\documentclass[10pt]{amsart}

\usepackage{amssymb}
\usepackage{amsmath}
\usepackage{amsfonts}
\usepackage[section]{placeins}
\usepackage{url}
\usepackage{microtype}
\usepackage{booktabs, multirow}
\usepackage[french, english]{babel}
\usepackage{amsmath,amssymb,amsthm, comment, mathrsfs, url, yfonts}
\usepackage{array}
\usepackage[colorlinks]{hyperref}
\hypersetup{citecolor=blue}
\usepackage{amscd}
\usepackage{todonotes}

\theoremstyle{plain}
\newtheorem{theorem}{Theorem}[subsection]
\newtheorem{proposition}[theorem]{Proposition}
\newtheorem{corollary}[theorem]{Corollary}
\newtheorem{lemma}[theorem]{Lemma}
\newtheorem{conjecture}[theorem]{Conjecture}
\newtheorem{ithm}{Theorem}

\newtheorem{iconj}{Conjecture}

\theoremstyle{definition}
\newtheorem{remark}[theorem]{Remark}
\newtheorem{definition}[theorem]{Def\/inition}

\newtheorem{caveat}[theorem]{Caveat}

\newcommand{\begpf}{\noindent{\bf Proof.}\enspace}
\newcommand{\epf}{{\ifhmode\unskip\nobreak\hfil\penalty50 \hskip1em
\else\nobreak\fi \nobreak\mbox{}\hfil\mbox{$\square$} \parfillskip=0pt
\finalhyphendemerits=0 \par\vskip5pt}}

\begin{document}
\title[A Serre weight conjecture for geometric Hilbert modular forms]{A Serre weight conjecture for geometric Hilbert modular forms in characteristic $p$}
\author{Fred Diamond}
\email{fred.diamond@kcl.ac.uk}
\address{Department of Mathematics, King's College London, Strand, London WC2R 2LS, United Kingdom}
\author{Shu Sasaki}
\email{s.sasaki.03@cantab.net}
\address{School of Mathematical Sciences, Queen Mary University of London, E1 4NS, UK}
\thanks{This research was partially supported by EPSRC Grant EP/L025302/1 and the Heilbronn Institute for Mathematical Research (FD),
and by EPSRC postdoctoral fellowship EP/G050511/1 and SFB/TR45 of DFG (SS)}
\begin{abstract}Let $p$ be a prime and $F$ a totally real field in which $p$ is unramified.
We consider mod $p$ Hilbert modular forms for $F$, defined as sections of
automorphic line bundles on Hilbert modular varieties of level prime to
$p$ in characteristic $p$. For a mod $p$ Hilbert modular Hecke eigenform of
arbitrary weight (without parity hypotheses), we associate a
two-dimensional representation of the absolute Galois group of $F$, and
we give a conjectural description of the set of weights of all eigenforms
from which it arises. This conjecture can be viewed as a ``geometric''
variant of the ``algebraic'' Serre weight conjecture of
Buzzard--Diamond--Jarvis, in the spirit of Edixhoven's variant of Serre's
original conjecture in the case $F = \mathbb{Q}$.  We develop techniques for studying
the set of weights giving rise to a fixed Galois representation, and prove
results in support of the conjecture, including cases of partial weight
one.
\end{abstract}

\maketitle

\section{introduction}

\subsection{The weight part of Serre's Conjecture}
Let $p$ be a rational prime.   Serre's Conjecture~\cite{Se}, now a theorem of Khare and Wintenberger~\cite{KW, KW2}
(completed by a result of Kisin~\cite{Ki2}) asserts that every odd, continuous, irreducible representation 
$\rho: \mathrm{Gal}(\overline{\mathbb{Q}}/\mathbb{Q})\rightarrow \mathrm{GL}_2(\overline{\mathbb{F}}_p)$
is {\em modular} in the sense that it is isomorphic to the
mod $p$ Galois representation associated to a modular eigenform.  Furthermore, Serre predicts
the minimal weight $k\geq 2$ such that $\rho$ arises from an eigenform of weight $k$
and level prime to $p$, the recipe for this minimal weight being in terms of the restriction of $\rho$ to an inertia
subgroup at $p$.   Under the assumption that $\rho$ is modular, the fact that it arises from an eigenform of Serre's predicted weight
was known prior to the work of Khare--Wintenberger (assuming $p>2$), and indeed this plays a crucial role in their
proof of Serre's Conjecture.  This fact, called the {\em weight part of Serre's Conjecture}, was proved by Edixhoven~\cite{E}
using the results of Gross~\cite{Gr} and Coleman--Voloch~\cite{CV} on companion forms.  Edixhoven also presents
(and proves for $p>2$\footnote{Both versions in the case $p=2$ should ultimately follow from the results of 
Khare--Wintenberger and Kisin, as explained in~\cite{Cal}.})
an alternative formulation, which predicts the minimal weight $k\ge 1$ such that $\rho$
arises from a mod $p$ eigenform of weight $k$ and level prime to $p$, where mod $p$ modular
forms are viewed as sections of certain line bundles on the reduction mod $p$ of a modular curve.   The qualitative
difference between the two versions of the conjecture stems from the fact that a mod $p$ modular form
of weight one does not necessarily lift to characteristic zero.

 There has been a significant amount of work towards generalising the original formulation
 of the weight part of Serre's Conjecture to other contexts where one has (or expects)
 Galois representations associated to automorphic forms.   This line of research was first
 developed by Ash and collaborators in the context of $\mathrm{GL}_n$ over $\mathbb{Q}$ 
 (in particular~\cite{AS}), and the most general formulation to date is due to Gee, Herzig
 and Savitt in~\cite{GHS}.   We refer the reader to the introduction of \cite{GHS} for a
 discussion of this history and valuable perspectives provided by representation theory,
 $p$-adic Hodge theory and the Breuil--M\'ezard Conjecture.
 
 An important setting for the development of generalisations of the weight part of Serre's
 Conjecture has been that of Hilbert modular forms, i.e., automorphic forms for $G = \mathrm{Res}_{F/\mathbb{Q}}\mathrm{GL}_2$ 
 where $F$ is a totally real field.  Work in this direction was initiated by Buzzard, Jarvis and one
 of the authors in~\cite{BDJ}, where a Serre weight conjecture is formulated under the assumption that
 $p$ is unramified in $F$.   For a totally odd, continuous, irreducible representation
 \begin{equation}\label{eqn:rho} \rho: \mathrm{Gal}(\overline{F}/F)\rightarrow \mathrm{GL}_2(\overline{\mathbb{F}}_p),\end{equation}
  there is a notion
 of $\rho$ being {\em modular of weight $V$}, where $V$ is an irreducible $\overline{\mathbb{F}}_p$-representation
 of $G(\mathbb{F}_p) = \mathrm{GL}_2(O_F/pO_F)$, where $O_F$ denotes the ring of integers of $F$.
 In this context, the generalisation of the weight part of Serre's Conjecture assumes that $\rho$ is modular
 of {\em some} weight, and predicts the set of {\em all} such weights in terms of the restriction of $\rho$ to inertia
 groups at primes over $p$.   This prediction can be viewed as a conjectural description of all pairs
 $(\tau_\infty,\tau_p)$ where $\tau_\infty$ (resp.~$\tau_p$) is a cohomological type at $\infty$
 (resp.~$K$-type at $p$) of an automorphic representation giving rise to $\rho$ (see~\cite[Prop.~2.10]{BDJ}).
 The conjecture was subsequently generalised in \cite{Sc,G2} to include the case where $p$ is ramified in $F$,
 and indeed proved under mild technical hypotheses (for $p>2$) 
 in a series of papers by Gee and collaborators culminating in \cite{GLS,GKi}, with an
 alternative endgame provided by Newton~\cite{N}.

 It is also natural to consider the problem of generalising Edixhoven's variant of the
 weight part of Serre's Conjecture, especially in view of the innovation due to Calegari--Geraghty~\cite{CG}
 on the Taylor--Wiles method for proving automorphy lifting theorems.
 By contrast with the original formulation of the weight part of Serre's Conjecture, there has been
 relatively little work in this direction.  The main aim of this paper is to formulate such a variant
 in the setting of Hilbert modular forms associated to a totally real field $F$ in which $p$ is
 unramified.  More precisely, for $\rho$ as in (\ref{eqn:rho}), we give a conjectural
 description of the set of all weights of mod $p$ Hilbert modular eigenforms giving rise to $\rho$,
where we view mod $p$ Hilbert modular forms as sections of certain
line bundles on the special fibre of a Hilbert modular variety.
Furthermore, we develop some tools for studying the set of possible weights,
and prove results towards the conjecture  in the first case that exhibits genuinely
new phenomena relative to the settings of~\cite{E} and~\cite{BDJ}.
A forthcoming paper of the authors generalises all the conjectures and
results of this article to the case where $p$ is ramified in $F$.
Another direction that warrants further research is consideration of
higher degree cohomology of automorphic bundles, and it would be
natural to pursue the problem in the context of more general Shimura varieties.

A key point we should make is that in the setting of classical modular forms, the enhancement
provided by Edixhoven's variant of Serre's Conjecture pertains essentially just to unramified-at-$p$ Galois representations and
weight one modular forms, but already in the Hilbert modular setting, a much richer tableau emerges
from the geometric variant, in terms of both the related $p$-adic Hodge theory and
arithmetic of automorphic forms.  We touch on this further  in the course of outlining the
contents of the paper below.

\subsection{Mod $p$ Hilbert modular forms and Galois representations}
The foundations for this paper have their roots in the work of Andreatta--Goren~\cite{AG}, which develops
the theory of mod $p$ Hilbert modular forms and partial Hasse invariants.  In particular, they use the
partial Hasse invariants to define the filtration, which we refer to instead as the {\em minimal weight},
of a mod $p$ Hilbert modular form.
However, the framework for~\cite{AG} is based on an alternate notion of Hilbert modular forms,
defined using Shimura varieties and automorphic forms associated to the reductive group $G^*$,
the preimage of $\mathbb{G}_m$ under $\det:G \to \mathrm{Res}_{F/\mathbb{Q}}\mathbb{G}_m$
where  $G =  \mathrm{Res}_{F/\mathbb{Q}}\mathrm{GL}_2$. We wish to work throughout with
automorphic forms with respect to $G$ itself, which are more amenable to the theory of
Hecke operators and associated Galois representations.  To this end we need to adapt
the setup of \cite{AG}.  

We begin by recalling the definition of Hilbert modular varieties
in \S\ref{sec:HMV} and Hilbert modular forms in~\S\ref{sec:hmfs} in our context.
For us, a weight will be a pair $(k,l) \in \mathbb{Z}^\Sigma \times \mathbb{Z}^\Sigma$,
where $\Sigma$ is the set of embeddings $F \to \mathbb{Q}$.
A fundamental observation is the absence in characteristic $p$ of the parity condition
on $k$ that appears in the usual definition of weights of Hilbert modular forms
(with respect to $G$, as opposed to $G^*$) in characteristic zero; our $(k,l)$ is
arbitrary.  In \S\ref{section:Hecke} we explain the construction of Hecke operators
in our setting, and in \S\ref{sec:Hasse} we recall (and adapt) the definition
of partial Hasse invariants from~\cite{AG}.

In \S\ref{sec:galois} we establish the existence of Galois representations associated to mod $p$
Hilbert modular eigenforms of arbitrary weight.  More precisely we prove (see Theorem~\ref{thm:galois}):
\begin{ithm} \label{ithm:galois}  If $f$ is a mod $p$ Hilbert modular eigenform of 
weight $(k,l)$ and level $U(\mathfrak{n})$ with $\mathfrak{n}$ prime to $p$, then
there is a Galois representation $\rho_f : G_F \to \mathrm{GL}_2(\overline{\mathbb{F}}_p)$
such that if $v\nmid \mathfrak{n}p$, then $\rho_f$ is unramified at $v$ and the characteristic
polynomial of $\rho_f(\mathrm{Frob}_v)$ is 
$X^2 - a_v  X + d_v \mathrm{Nm}_{F/\mathbb{Q}}(v)$,
where $T_v f = a_v f$ and $S_v f = d_v f$.
\end{ithm}
This was proved independently by Emerton--Reduzzi--Xiao~\cite{ERX} and Goldring--Koskivirta~\cite{GoK}
under parity hypotheses on $k$.  The new ingredient allowing us to treat arbitrary $(k,l)$ is to
use congruences to forms of level divisible by primes over $p$.  
This introduces a number of technical difficulties, the most critical of which
is overcome using a cohomological vanishing result proved in joint work
with Kassaei in~\cite{DKS}.

\subsection{A geometric Serre weight conjecture}
In \S\ref{sec:conj}, we introduce the notion of {\em geometric} modularity and formulate a conjecture that specifies
the set of weights for which a given $\rho$ is geometrically modular. A key point is that the geometric setting
allows for the notion of a minimal weight (among the possible $k$ for a fixed $l$) of eigenforms giving rise to $\rho$,
something not apparent in the framework of \cite{BDJ}.  An investigation of this phenomenon and its interaction
with properties of $\Theta$-cycles (under which $l$ varies) in this context led us to the expectation that this minimal
weight should lie in the cone $\Xi_{{\mathrm{min}}}^+ = \Xi_{{\mathrm{min}}} \cap \mathbb{Z}^\Sigma_{>0}$, where
$$\Xi_{{\mathrm{min}}}   := \{\, k \in \mathbb{Z}^\Sigma \, |\, \mbox{$pk_\tau \ge k_{\mathrm{Fr}^{-1}\circ\tau}$ for all $\tau \in \Sigma$}\,\},$$
and that geometric modularity of $\rho$ for weights in $\Xi_{{\mathrm{min}}}^+$
can be characterised using $p$-adic Hodge theory.
In particular, we make the following conjecture (see Conjecture~\ref{conj:geomweights} for a
stronger version, and Definitions~\ref{def:HTtype} and~\ref{def:weightlift} for conventions on Hodge--Tate weights):
\begin{iconj}
Suppose that $\rho: G_F \rightarrow {\mathrm{GL}}_2(\overline{\mathbb{F}}_p)$ is irreducible and
geometrically modular (of some weight), and let $l \in \mathbb{Z}^\Sigma$. Then
there exists $k_{\mathrm{min}} = k_{\mathrm{min}}(\rho,l) \in \Xi_{\mathrm{min}}^+$ such that the following
holds for all $k \in \Xi_{\mathrm{min}}^+$:
$$\begin{array}{c}
 \mbox{$\rho$ is geometrically modular of weight $(k,l)$  \; $\Longleftrightarrow$ \; $k \ge_{\mathrm{Ha}} k_{\mathrm{min}}$} \\ \ \\
\mbox{$\Longleftrightarrow$ \; $\rho|_{G_{F_v}}$  has a crystalline lift of weight $(k_\tau,l_\tau)_{\tau\in \Sigma_v}$ for all $v|p$.}\end{array}$$
\end{iconj}
The inequality in the statement means that $k - k_{\mathrm{min}}$ is a non-negative integral linear combination of weights of partial
Hasse invariants, and the existence of a weight $k_{\mathrm{min}}$ such that the second equivalence holds already
has deep consequences in $p$-adic Hodge theory which are not apparent from the framework of algebraic
Serre weight conjectures.  We also stress that one might have expected that, by analogy with the 
algebraic setting in \cite{BDJ}, the preceding conjecture held with $\Xi_{\mathrm{min}}^+$ replaced by
the set of totally positive weights.  Indeed we had been positing such a formulation to experts
on Serre weight conjectures, until R.~Bartlett showed us a counterexample to the resulting $p$-adic
Hodge-theoretic implications; this led us to examine the possible $\Theta$-cycles more closely and arrive
at the above version.   The role of $\Xi_{{\mathrm{min}}}$  in the conjecture in turn inspired the main result of \cite{DK},
thus answering the basic question posed in \cite{AG} of whether the minimal weight of a mod $p$ Hilbert modular
form is totally non-negative; in fact \cite[Cor.~1.2]{DK} establishes the stronger result that it lies in $\Xi_{\min}$.

In \S\ref{sec:conj} we also explain the relation with the Serre weight conjectures of~\cite{BDJ},
which can be viewed as specifying the weights $(k,l) \in \mathbb{Z}^\Sigma_{\ge 2} \times \mathbb{Z}^\Sigma$
(i.e., {\em algebraic weights}) for which $\rho$ is {\em algebraically} modular.
In particular Conjecture~\ref{conj:algvsgeomweights} predicts that
geometric and algebraic modularity of $\rho$ for a weight $(k,l)$ are equivalent if
$k \in \mathbb{Z}^\Sigma_{\ge 2} \cap \Xi_{{\mathrm{min}}}$.  The hypothesis
$k \in \mathbb{Z}^\Sigma_{\ge 2}$ is needed for the notion of algebraic modularity,
and it is not hard to see that it implies geometric modularity even without the assumption $k \in  \Xi_{{\mathrm{min}}}$,
at least if all $k_\tau$ have the same parity (see Proposition~\ref{prop:algvsgeomweights}).  On the other hand, the opposite implication
seems much more difficult, and its failure for $k \not\in \Xi_{{\mathrm{min}}}$ can be observed through
the optic of modular representation theory (see Remark~\ref{rmk:TobyQ}).  Furthermore for {\em algebraic} weights in  $\Xi_{{\mathrm{min}}}$,
the Breuil--M\'ezard Conjecture~\cite{BM} (or more precisely its generalisation in~\cite[Conj.~1.1.5]{GKi})
converts the $p$-adic Hodge-theoretic implications of Conjecture~\ref{conj:geomweights} into
non-trivial results on the mod $p$ representation theory of ${\mathrm{GL}}_2({\mathbb{F}}_{p^r})$
(see Wiersema's PhD thesis~\cite{HWPhD}).    We remark also that the interplay between
algebraic and geometric Serre weights is reflected in the geometry
of Hilbert modular varieties; this theme motivated the construction in~\cite{DKS} of a filtration
and Jacquet--Langlands relation for mod $p$ Hilbert modular forms of pro-$p$-Iwahori
level at $p$.

In the case $F = \mathbb{Q}$, the only non-algebraic weights with $k \in \Xi_{{\mathrm{min}}}^+$
have the form $(1,l)$, for which Conjecture~\ref{conj:geomweights2}\footnote{Strictly speaking,
we assume $F \neq \mathbb{Q}$ throughout the paper to allow for a more uniform exposition,
and because nothing new would be presented in the case $F= \mathbb{Q}$.  For $F =\mathbb{Q}$, the equivalence
between algebraic and geometric modularity for $k\ge 2$ is standard, and the analogue of Conjecture~\ref{conj:geomweights2}
reduces via the Breuil--M\'ezard Conjecture to Edixhoven's variant of the weight part of Serre's Conjecture,
hence is known (at least if $p > 2$).} reduces to the statement that $\rho$ is unramified at $p$ if and only if it arises from a mod $p$
eigenform of weight one.  For general $F$, the analogous statement relating (modular) $\rho$ unramified at $p$
to parallel weight one forms is established (under technical hypotheses) by work of Gee--Kassaei~\cite{GK} and
Dimitrov--Wiese~\cite{DW}; however there is a much richer range of possibilities for $\rho$ to arise (minimally)
from forms with non-algebraic minimal weights.   The first instance where this is apparent is for real quadratic fields $F$
in which $p$ is inert, and we investigate this in detail in \S\ref{sec:quadratic}.

\subsection{Partial $\Theta$-operators, $q$-expansions and the inert quadratic case}
We have already indicated how the perspective afforded by Conjectures~\ref{conj:geomweights}
and~\ref{conj:algvsgeomweights} leads to new results on the geometry of Shimura varieties and
mod $p$ automorphic forms, as in \cite{DK} and \cite{DKS}.  On the other hand, progress on these conjectures
evidently requires more geometric input than was needed for the proof of the algebraic Serre weight
conjecture of \cite{BDJ}.   To that end, we review and develop several useful general tools in
the context of Hilbert modular varieties, beginning with the theory of $\Theta$-operators in \S\ref{sec:Theta}.
We again proceed by adapting the treatment in \cite{AG}, but in doing so we introduce some new perspectives
which we feel simplify and clarify some aspects of their construction.  (See Remark~\ref{rmk:logdiffs}
and the proof of Theorem~\ref{thm:theta}.)  Furthermore the approach taken here leads to a substantial
improvement on the results in \cite{AG} on $\Theta$-operators when $p$ is ramified in $F$; see~\cite{PRtheta}.

In the last few sections, we make critical use of $q$-expansions.  Most of \S\ref{sec:qexpansions} is 
a straightforward application of standard methods and results describing $q$-expansions
and the effect on them of Hecke operators.  In addition to this, we construct partial Frobenius operators, whose
image we relate to the kernel of partial $\Theta$-operators in Theorem~\ref{thm:Phiv}; this argument is
(to our knowledge) new, and the result generalises a theorem of Katz~\cite{K2}.  These results are further
generalised in~\cite{PRtheta} to the case where $p$ is ramified in $F$.

In \S\ref{sec:eigenforms} we first prove various technical results on eigenforms and their $q$-expansions.
We then study the behaviour of the minimal weight for $\rho$ as $l$ varies (see Theorem~\ref{thm:thetarho}),
and prove that if an eigenform of algebraic weight is ordinary at a prime over $p$, then so is the associated
Galois representation (Theorem~\ref{thm:ordinary}).

Finally in \S\ref{sec:quadratic} we specialise to the inert quadratic case.  We first use results
from integral $p$-adic Hodge theory to describe those $\rho$ for which the 
(conjectural) minimal weight is not algebraic (i.e., has $k_\tau = 1$ for some $\tau$).
We then use the tools developed in the preceding sections to transfer
modularity results between algebraic and non-algebraic weights.
In particular we prove cases of Conjecture~\ref{conj:geomweights2} in the setting of
partial weight one, conditional on our conjectured equivalence between
algebraic and geometric modularity (see Theorem~\ref{thm:quadratic1}).
Since one direction of this equivalence is easy under a parity hypothesis,
we also obtain the following unconditional result (Theorem~\ref{thm:quadratic2}):
\begin{ithm}  Suppose that $[F:\mathbb{Q}] = 2$, $p$ is inert in $F$,
$3 \le k_0 \le p$, $k_0$ is odd and $\rho:G_F \to \mathrm{GL}_2(\overline{\mathbb{F}}_p)$
is irreducible and modular.  
If $\rho|_{G_{F_p}}$ has a crystalline lift of weight $((k_0,1),(0,0))$, then
$\rho$ is geometrically modular of weight $((k_0,1),(0,0))$.
\end{ithm}
Our method is indicative of a general strategy for transferring results
from the setting of algebraic Serre weight conjectures to their geometric variants;
further work in this direction is carried out in Wiersema's thesis~\cite{HWPhD}.

\subsection{Acknowledgements}  We are grateful to Payman Kassaei for many
valuable discussions related to this work.   We also thank Fabrizio Andreatta,
Eyal Goren, Kai-Wen Lan and David Savitt for helpful correspondence, and Robin Bartlett,
Toby Gee, David Helm and Vytautas Pa\v{s}k\=unas for useful conversations.
This research also benefitted from the authors' participation in the Workshop
on Serre's Modularity Conjecture at the University of Luxembourg in June 2015 (organised
by Mladen Dimitrov, Haluk \c{S}eng\"un, Gabor Wiese and Hwajong Yoo),
especially through discussions with Mladen Dimitrov and Jacques Tilouine.

We wish to thank Chris Birkbeck, Hanneke Wiersema and two anonymous referees for
carefully reading parts of this paper and suggesting several clarifications and minor corrections.
Finally, we are grateful to yet another referee for making us aware of
deficiencies of an earlier version of the introduction to this paper, and to Toby Gee
for helpful feedback leading to its improvement.

\section{Hilbert modular varieties}

\label{sec:HMV}

In this section we recall the definitions and basic properties of the models for Hilbert modular
varieties used throughout the paper.

\subsection{General notation}
Let $p$ be a fixed rational prime. Let $F$ be a totally real field in which $p$ is unramified. 
We let $O_F$ denote the ring of integers of $F$, and $O_{F,\ell} = O_F\otimes \mathbb{Z}_\ell$
for any prime $\ell$.

Since this paper offers nothing new in the case $F = \mathbb{Q}$ (relative to \cite{E}), we will assume throughout that $F \neq \mathbb{Q}$ in order to avoid complications arising from consideration of the cusps.

Let $\mathfrak{d}=\mathfrak{d}_{F/\mathbb{Q}}$ denote the different of $F$ over $\mathbb{Q}$. Fix algebraic closures $\overline{\mathbb{Q}}, \overline{\mathbb{Q}}_p$ of $\mathbb{Q}$ and $\mathbb{Q}_p$ respectively, and fix embeddings of $\overline{\mathbb{Q}}$ into $\overline{\mathbb{Q}}_p$ and $\mathbb{C}$.

Let $\Sigma$ denote the set embeddings of $F$ into $\overline{\mathbb{Q}}$. 
Let $L$ denote a finite extension of $\mathbb{Q}_p$ in $\overline{\mathbb{Q}}_p$ containing the image of every embedding in $\Sigma$, $\mathcal{O}$ its ring of integers, $\pi$ a uniformiser and $E = \mathcal{O}/\pi$ its residue field.
We identify $\Sigma$ with the set of embeddings of $F$ into $L$ (and hence of $O_F$
into $\mathcal{O}$), as well as the set of embeddings of $F$ into $\mathbb{R}$.

If $T$ is a subset of $F_\infty = F\otimes \mathbb{R} \cong \prod_{\tau \in \Sigma} \mathbb{R}$, we let
$T_+$ the set of totally positive elements in $T$.

\subsection{Hilbert modular varieties of level $N$}
\begin{definition} For a fractional ideal $J$ of $F$ and an integer $N\ge 3$, let $\mathcal{M}_{J, N}$ denote the functor which sends an $\mathcal{O}$-scheme $S$ to the set of isomorphism classes of data $(A, i, \lambda, \eta)$ comprising
\begin{itemize}

\item an abelian scheme $A/S$ of relative dimension $[F:\mathbb{Q}]$, 

\item a ring homomorphism $\iota: O_F\rightarrow \mathrm{End}(A/S)$, 

\item an $O_F$-linear isomorphism $\lambda: (J,J_+) \simeq (\mathrm{Sym}(A/S),\mathrm{Pol}(A/S))$
such that the induced map $A\otimes_{O_F} J \rightarrow A^\vee$ is an isomorphism,
where $\mathrm{Sym}(A/S)$ (resp. $\mathrm{Pol}(A/S)$) denotes the \'etale sheaf whose sections are symmetric $O_F$-linear morphisms (resp. polarisations) $A\rightarrow A^\vee$,

\item an $O_F$-linear isomorphism $\eta: (O_F/N)^2 \simeq A[N]$.
   
\end{itemize}
We call such a quadruple a {\em $J$-polarised Hilbert-Blumenthal abelian variety with level $N$ structure} (or simply an {\em HBAV} when $J$ and $N$ are fixed) over $S$.
\end{definition}

The functor $\mathcal{M}_{J,N}$ is representable by a smooth $\mathcal{O}$-scheme,
which we shall denote $Y_{J, N}$; see \cite[Thm.~2.2]{DP} and the discussion before
it, from which it also follows (using for example \cite[Thm.~1.4]{C}) that $Y_{J,N}$
is quasi-projective over $\mathcal{O}$.

Let $Z_{J,N}$ denote the finite $\mathcal{O}$-scheme representing
$O_F$-linear isomorphisms $J/NJ \simeq \mathfrak{d}^{-1}\otimes \mu_N$.
If $(A,\iota,\lambda,\eta)$ is an HBAV over $S$, then $\lambda \otimes \wedge^2\eta$
defines an isomorphism
$$J/NJ = J \otimes_{O_F} \wedge^2_{O_F}(O_F/N)^2 \simeq \mathrm{Sym}(A/S)\otimes_{O_F}\wedge^2_{O_F} A[N],$$
where $A[N]$ is viewed as an \'etale sheaf on $S$.  Composing with the isomorphisms
$$\mathrm{Sym}(A/S)\otimes_{O_F} \wedge^2_{O_F} A[N] \simeq \mathrm{Hom}(O_F,\mu_N)
 \simeq \mathfrak{d}^{-1}\otimes \mu_N$$
 induced by the Weil pairing and the trace pairing thus gives an element of $Z_{J,N}(S)$.
 In particular taking $S = Y_{J,N}$ and the universal HBAV over it, we obtain a canonical morphism $Y_{J,N} \to Z_{J,N}$
 with geometrically connected fibres.

\subsection{Unit action on polarisations}
The group $O_{F, +}^\times $ of totally positive units in $O_F$ acts on $Y_{J, N}$  by $\nu$ in $O_{F, +}^\times$ sending $(A, \iota, \lambda, \eta)\in Y_{J, N}(S)$ for every $\mathcal{O}$-scheme $S$ to $(A, \iota, \nu\lambda, \eta)\in Y_{J, N}(S)$.  Similarly $u\in \mathrm{GL}_2(O_F/N O_F)$ acts by sending  $(A, \iota, \lambda, \eta)$ to  $(A, \iota, \lambda, \eta \circ r_{u^{-1}})$ where $r_{u^{-1}}$ denotes right multiplication by $u^{-1}$, thus defining a right action of $\mathrm{GL}_2(O_F/N O_F)$, and hence of $\mathrm{GL}_2(\widehat{O}_F)$ on $Y_{J,N}$ through the projection $\mathrm{GL}_2(\widehat{O}_F)\rightarrow \mathrm{GL}_2(O_F/NO_F)$ where $\widehat{O}_F$ denotes the profinite completion of $O_F$. If $\mu \in O_F^\times$, then the action of $\mu^2 \in O_{F,+}^\times$ on $Y_{J,N}$ coincides with that of $\mu^{-1} I_2 \in \mathrm{GL}_2(\widehat{O}_F)$ (where $I_2$ denote the 2-by-2 identity matrix).

\subsection{Adelic action on level structures}
Now let $U$ be an open compact subgroup of $\mathrm{Res}_{F/\mathbb{Q}}\mathrm{GL}_2(\widehat{\mathbb{Z}}) \simeq \mathrm{GL}_2(\widehat{O}_F)$ containing $\mathrm{GL}_2(O_{F,p})$.  Choose an integer $N \ge 3$ such that $N$ is not divisible by $p$ and $U(N) \subset U$, where $U(N):= \ker( \mathrm{GL}_2(\widehat{O}_F) \to  \mathrm{GL}_2(O_F/NO_F)$. Then the action of $O_{F,+}^\times \times \mathrm{GL}_2(\widehat{O}_F)$ induces one on $Y_{J,N}$ of the finite group
$$G_{U,N} := (O_{F,+}^\times \times U) / \{\,(\mu^2,u) \,|\, \mu \in O_F^\times, u \in U, u \equiv \mu I \bmod N\,\}.$$
Note that the action of  $(\nu,u) \in O_{F,+}^\times \times \mathrm{GL}_2(\widehat{O}_F)$ on $Y_{J,N}$
is compatible with the natural action on $Z_{J,N}$ defined by multiplication by $\nu\det(u)^{-1}$.

We will show that if $U$ is sufficiently small, then $G_{U,N}$ acts freely on $Y_{J,N}$.  To make
this precise, let $\mathcal{P}_F$ denote the set of primes $r$ in $\mathbb{Q}$ such that the maximal totally real subfield $\mathbb{Q}(\mu_r)^+ $ of $\mathbb{Q}(\mu_r)$ is contained in $F$, and let $\mathcal{C}_F$ denote the set of quadratic CM-extensions $K/F$ (in a fixed algebraic closure of $F$) such that either:
\begin{itemize}
\item $K = F(\mu_r)$ for some odd prime $r \in \mathcal{P}_F$, or
\item $K = F(\sqrt{\beta})$ for some $\beta \in O_F^\times$.
\end{itemize}
Note that the sets $\mathcal{P}_F$ and $\mathcal{C}_F$ are finite.

For an ideal $\mathfrak{n}$ of $O_F$, we define the following open compact subgroups of $\mathrm{GL}_2(\widehat{O}_F)$:
$$\begin{array}{rcl}
U_0(\mathfrak{n}) &:= &\left\{\,\left.\left(\begin{array}{cc} a&b \\ c & d \end{array}\right) \in \mathrm{GL}_2(\widehat{O}_F)\,\right|\,c \in \mathfrak{n}\widehat{O}_F\,\right\}; \\&&\\
U_1(\mathfrak{n}) &:= &\left\{\,\left.\left(\begin{array}{cc} a&b \\ c & d \end{array}\right) \in U_0(\mathfrak{n})\,\right|\,d - 1 \in \mathfrak{n}\widehat{O}_F\,\right\}; \\&&\\
{}^1U_1(\mathfrak{n}) &:= &\left\{\,\left.\left(\begin{array}{cc} a&b \\ c & d \end{array}\right) \in U_1(\mathfrak{n})\,\right|\,a - 1 \in \mathfrak{n}\widehat{O}_F\,\right\}.
\end{array}$$

\begin{lemma} \label{lemma:neat} 
Suppose that one of the following holds:
\begin{itemize}
\item $U \subset {}^1U_1(\mathfrak{n})$ for some $\mathfrak{n}$ such that if $r \in \mathcal{P}_F$, then $\mathfrak{n}$ does not contain $\mathfrak{r}O_F$ where $\mathfrak{r}$ is the prime over $r$ in $\mathbb{Q}(\mu_r)^+$, or
\item $U \subset U_0(\mathfrak{n})$ for some $\mathfrak{n}$ such that if $\mu_r \subset K$ and $K  \in \mathcal{C}_F$, then $\mathfrak{n} \subset \mathfrak{q}$ for some prime $\mathfrak{q}$ of $F$ inert in $K$ and not dividing $r$.
\end{itemize}
Then $G_{U,N}$ acts freely on $Y_{J,N}$.
\end{lemma}

\begpf  For $G_{U,N}$ to act freely on $Y_{U,N}$ means that the morphism $G_{U,N} \times Y_{U,N} \to Y_{U,N} \times_{\mathcal{O}} Y_{U,N}$ defined by $(g,x) \mapsto (gx,x)$ is a closed
immersion.  Since this morphism is finite and the fibre over every closed point is reduced, it suffices to prove that for every geometric point $x \in Y_{U,N}(S)$, the map $G_{U,N} \to Y_{U,N}(S)$ defined by $g  \mapsto gx$ is injective, i.e., that the stabiliser of $x$ in $G_{U,N}$ is trivial.

Suppose then that $(A,\iota,\lambda,\eta)$ is an HBAV over an algebraically closed field, and that $(\nu,u) \in O_{F,+}^\times \times U$ is such that $(A,\iota,\nu \lambda, \eta\circ r_{u^{-1}})$ is isomorphic to  $(A,\iota,\lambda, \eta)$.  This means that there is an automorphism $\alpha$ of $A$ such that $\alpha$ commutes with the action of $O_F$ and satisfies $\alpha\circ \eta = \eta \circ r_{u^{-1}}$ and $\lambda(j) =\alpha^\vee\circ \lambda(\nu j) \circ \alpha$ for $j \in J$.

We wish to prove that $\alpha = \iota(\mu)$ for some $\mu \in O_F^\times$.  Suppose this is not the case. Viewing $F$ as a subfield of $\mathrm{End}^0(A) = \mathbb{Q} \otimes \mathrm{End}(A)$ via $\iota$, it follows from the classification of endomorphism algebras of abelian varieties that $F(\alpha)$ is a quadratic CM-extension $K/F$. Since $\alpha$ is an automorphism, it is a unit in an order in $K$, so $\alpha \in O_K^\times$. Since $O_F^\times$ and $O_K^\times$ have the same rank and $\alpha \notin O_F^\times$, we have $\alpha^n \in O_F^\times$ for some $n>0$; replacing $\alpha$ by a power, we may assume $n$ is a prime $r$.  Since $\alpha^r \in F$ and $K = F[\alpha]$ is Galois over $F$, it follows that $\zeta_r \in K$, and hence either $K = F(\mu_r)$ or $r = 2$. In either case we conclude that $r \in \mathcal{P}_F$, $\mu_r \subset K$ and $K \in \mathcal{C}_F$.

Now let $f(X)$ denote the minimal polynomial of $\alpha$ over $F$.  Note that since $\alpha^r \in O_F^\times$, we have
\begin{equation}\label{charpoly} f(X) = (X - \alpha)(X-\zeta_r\alpha)= X^2 - (1+\zeta_r)\alpha X + \zeta_r\alpha^2\end{equation}
for some $\zeta_r \in \mu_r$. For each prime $\ell$ of $F$ not dividing $p$, the $\ell$-adic Tate module $T_\ell(A)$ is free of rank two over $O_{F,\ell}$ and is annihilated by $f(\alpha)$, so $f(X)$ is in fact the characteristic polynomial of $\alpha$ on $T_\ell(A)$. It follows that $f(X)$ is the characteristic polynomial of $\alpha$ on $A[N]$, and hence also the characteristic polynomial of $u$ on $(O_F/N)^2$.

Suppose now that $U$ is as in the first bullet in the statement of the lemma,  Since $U \subset {}^1U_1(\mathfrak{n})$, the characteristic polynomial of $u$ is $(X-1)^2 \bmod \mathfrak{n}$. Comparing with (\ref{charpoly}), we see that $(1+\zeta_r)\alpha \equiv 2 \bmod \mathfrak{n}$ and $\zeta_r\alpha^2 \equiv 1 \bmod \mathfrak{n}$. If $r = 2$, this implies $2\in \mathfrak{n}$, contradicting the hypothesis on $\mathfrak{n}$. If $r$ is odd, this implies $\zeta_r\alpha^2(\zeta_r - \zeta_r^{-1})^2 \in \mathfrak{n}$; since $\zeta_r\alpha^2 \in O_F^\times$ and $(\zeta_r - \zeta_r^{-1})^2$ generates $\mathfrak{r}$, this also contradicts the hypothesis on $\mathfrak{n}$.

Suppose now that $U$ is as in the second bullet of the statement. Then there is a prime $\mathfrak{q}$  dividing $\mathfrak{n}$ such that $\mathfrak{q}$ is inert in $K$ and does not divide $r$. Since the discriminant of $f(X)$ is only divisible by primes over $r$, we have $O_{K,\mathfrak{q}} = O_{F,\mathfrak{q}}[\alpha]$, so $f(X)$ is irreducible modulo $\mathfrak{q}$.  On the other hand, since $u \in U_0(\mathfrak{q})$ its characteristic polynomial factors over $O_F/\mathfrak{q}$, and we again obtain a contradiction.

We have now shown that $\alpha = \iota(\mu)$ for some $\mu \in O_F^\times$.  It follows that $u^{-1} \equiv \mu I \bmod N$, and that $\nu = \mu^{-2}$. Therefore the image of $(\nu,u)$ in $G_{U,N}$ is trivial, as required. \epf 

\begin{caveat}\label{def:neat} Unless otherwise indicated, we assume throughout the paper that the 
open compact subgroup $U$ of $\mathrm{GL}_2(\widehat{O}_F)$ contains $\mathrm{GL}_2(O_{F,p})$ and
is sufficiently small that the conclusion of Lemma \ref{lemma:neat} holds for some, hence all, $N\ge 3$
such that $U(N) \subset U$.
\end{caveat}

\subsection{Hilbert modular varieties of level $U$}
We fix a set $T$ of representatives $t$ in $(\mathbb{A}_F^{\infty})^\times$ for the strict ideal class group $(\mathbb{A}_F^{\infty})^\times/F_+^\times\widehat{O}_F^\times \cong
\mathbb{A}_F^{\times}/F^\times\widehat{O}_F^\times F_{\infty,+}^\times$, and let $J_t$ denote the corresponding fractional ideal of $F$.
We assume the representatives $t$ are chosen
so that the $J_t$ are prime to $p$; i.e., that $t_p \in O_{F,p}^\times$ for each $t \in T$.

Since $Y_{J_t,N}$ is quasi-projective over $\mathcal{O}$, the quotient $Y_{J_t,N}/G_{U,N}$ is representable by a scheme over $\mathcal{O}$ (by \cite[Prop.V.1.8]{SGA1}), and we define
$$Y_U =  \coprod_{t \in T}  Y_{J_t, N}/G_{U,N}.$$
Then $Y_U$ is smooth over $\mathcal{O}$ and the  projection $\coprod_{t \in T}  Y_{J_t, N} \to Y_U$ is Galois and \'etale with Galois group $G_{U,N}$
(in view of Lemma~\ref{lemma:neat} and Caveat~\ref{def:neat}).
Moreover $Y_U$ is defined over $\mathcal{O}\cap\overline{\mathbb{Q}}$ and is independent of the choices of $N$ and $T$.  

\subsection{Components}
\label{subsec:components}
Let $Z_U = \coprod_{t\in T} Z_{J_t,N}/G_{U,N}$ with $(\nu,u) \in G_{U,N}$
acting by multiplication by $\nu\det(u)^{-1}$, so if 
$\mu_N(\overline{\mathbb{Q}}) \subset \mathcal{O}$, then
$Z_U(\mathcal{O})$ can be identified with the set of (geometrically)
connected components of $Y_U$.  Fixing a generator $\zeta_N$ for
$\mathfrak{d}^{-1} \otimes \mu_N(\mathcal{O})$ as an $O_F$-module,
we obtain a bijection 
\begin{equation}\label{eqn:components}
(\mathbb{A}_F^\infty)^\times/F_+^\times \det(U) \simeq Z_U(\mathcal{O})
\end{equation}
by sending $x F_+^\times \det(U)$ to the $G_{U,N}$-orbit of the isomorphism
$J_t/NJ_t \simeq \mathfrak{d}^{-1}\otimes \mu_N(\mathcal{O})$ sending the class
of $(x\alpha)^{-1} \in J_t \otimes_{O_F}\widehat{O}_F$ to $\zeta_N$, where 
$t \in T$ and $\alpha \in F_+^\times$ (unique up to multiplication by an
element of $O_{F,+}^\times$) are chosen so that $x^{-1} \in t \alpha \widehat{O}_F^\times$.

\subsection{Complex points}
We recall that $Y_U$ is defined over $\mathcal{O} \cap \overline{\mathbb{Q}}$, and 
a standard construction yields an isomorphism 
\begin{equation}\label{eqn:chmv} \mathrm{GL}_2(F) \backslash 
     \mathrm{GL}_2(\mathbb{A}_F)/UU_\infty  \simeq Y_U(\mathbb{C})\end{equation}
where $U_\infty = \prod_{\tau \in \Sigma} \mathrm{SO}_2(\mathbb{R})\mathbb{R}^\times \subset \prod_{\tau\in\Sigma}\mathrm{GL}_2(\mathbb{R}) = \mathrm{GL}_2(F_\infty)$,
allowing us to view $Y_U$ as a model for the Hilbert modular variety of level $U$.
More precisely, by the Strong Approximation Theorem, any double coset as in (\ref{eqn:chmv})
can be written in the form $\mathrm{GL}_2(F) g_\infty \mathrm{diag}(1,x)UU_\infty$
for some $g_\infty \in \mathrm{GL}_2(F_\infty)$,
$x \in (\mathbb{A}_F^\infty)^\times$, such that $\det(g_\infty) \in F_{\infty,+}^\times$
and $x\widehat{O}_F = J\mathfrak{d}\widehat{O}_F$ for some $J$.
Such a double coset corresponds under (\ref{eqn:chmv}) to the $G_{U,N}$-orbit of
the HBAV over $\mathbb{C}$ defined by
$$\mathbb{C}\otimes O_F / (g_\infty(z_0)O_F \oplus  (J\mathfrak{d})^{-1})$$
with the evident $O_F$-action, isomorphism 
$\lambda: (J,J_+) \simeq (\mathrm{Sym}(A/S),\mathrm{Pol}(A/S))$ defined so that
$\lambda(\alpha)$ corresponds to the Hermitian form
$\mathrm{tr}_{F/\mathbb{Q}}(\alpha s \bar{t}/\mathrm{Im}(g_\infty(z_0))$, and level
$N$-structure defined by $(a,b) \mapsto (a g_\infty(z_0) + bx^{-1})/N$, where
$z_0 = i \otimes 1 \in \mathbb{C} \otimes O_F$.

\section{Hilbert modular forms}

\label{sec:hmfs}

In this section we recall the definition of Hilbert modular forms as sections of certain line
bundles on Hilbert modular varieties.

\subsection{Automorphic line bundles}
The condition that $A\otimes_{O_F} J \rightarrow A^\vee$ is an isomorphism (in the definition of an HBAV) is called the ``Deligne--Pappas'' condition.  
Our assumption that $p$ is unramified in $F$ ensures its equivalence with the
``Rapoport condition'' that $\mathrm{Lie}(A/S)$
is, locally on $S$, free of rank one over $O_F \otimes \mathcal{O}_S$ (\cite[Cor.~2.9]{DP}),
 and hence so is its $\mathcal{O}_S$-dual $e^*\Omega^1_{A/S} \simeq s_*\Omega^1_{A/S}$, where $s:A \to S$ is the structure morphism and $e:S\to A$ is the identity section. Since $O_F\otimes \mathcal{O}_S \simeq \bigoplus_{\tau\in\Sigma} \mathcal{O}_S$ as a coherent sheaf of $\mathcal{O}_S$-algebras, we may accordingly decompose $s_*\Omega^1_{A/S}$ as a direct sum of line bundles on $S$. Applying this to the universal HBAV $A_{J,N}$ over $Y_{J,N}$, we obtain a decomposition
 $s_*\Omega^1_{A_{J,N}/Y_{J,N}} = \bigoplus_{\tau\in\Sigma} \omega_\tau$ where each $\omega_\tau$ is a line bundle on $S$. For a  tuple $k=(k_\tau)_{\tau\in \Sigma} \in \mathbb{Z}^\Sigma$, we let $\omega^{\otimes k}$ denote the line bundle $\bigotimes_\tau \omega_\tau^{\otimes k_\tau}$ on $Y_{J,N}$.
 \begin{remark}
 Note that the definition of $\omega^{\otimes k}$ makes sense ``integrally'' because $p$ is assumed to be unramified in $F$ so that the Rapoport condition is satisfied; 
 in the ramified case, one can instead proceed as in \cite{S} using models for Hilbert modular varieties defined by Pappas and Rapoport in~\cite{PR}.
 \end{remark}

Since $\mathcal{H}^1_{\mathrm{DR}}(A/S) = R^1s_*\Omega^\bullet_{A/S}$ is locally free of rank two over $O_F \otimes \mathcal{O}_S$ (by \cite[Lem.~1.3]{R}) sitting in the exact sequence 
$$0\rightarrow s_\ast \Omega^1_{A/S}\rightarrow \mathcal{H}^1_{\mathrm{DR}}(A/S) \rightarrow R^1s_\ast\mathcal{O}_A\rightarrow 0$$
of locally free modules over $O_F\otimes\mathcal{O}_S$ (given by the Hodge--de Rham spectral sequence), 
\begin{equation}\label{wedge de Rham} \wedge^2_{O_F \otimes \mathcal{O}_S}\mathcal{H}^1_{\mathrm{DR}}(A/S) \simeq 
s_*\Omega^1_{A/S} \otimes_{_{O_F\otimes\mathcal{O}_S}} R^1s_*\mathcal{O}_A\end{equation}
is locally free of rank one over $O_F \otimes \mathcal{O}_S$ and similarly decomposes as a direct sum of line bundles indexed by $\tau \in \Sigma$. We let $\delta_\tau$ denote the line bundles so obtained from the universal HBAV over $S=Y_{J,N}$, and for a tuple
$l=(l_\tau)_{\tau\in \Sigma} \in \mathbb{Z}^\Sigma$, we let $\delta^{\otimes l}$ denote the line bundle $\bigotimes_\tau \delta_\tau^{\otimes l}$. Finally we let $\mathcal{L}_{J,N}^{k,l}$ denote the
line bundle $\omega^{\otimes k}\otimes\delta^{\otimes l}$.

Recall that we defined the action of $O_{F, +}^\times \times \mathrm{GL}_2(\widehat{O}_F) $ on $S = Y_{J, N}$  by requiring the pull-back via $(\nu,u)$ of the universal HBAV $(A, \iota, \lambda, \eta)$ to be isomorphic to $(A,\iota,\nu\lambda,\eta\circ r_{u^{-1}})$; we let $\alpha_{\nu,u} : A \to (\nu,u)^*A$ be the unique such isomorphism.  Note that $((\nu,u)^*\alpha_{\nu',u'}) \circ \alpha_{\nu,u} = \alpha_{\nu\nu',uu'}$ for $(\nu,u),(\nu',u') \in  O_{F, +}^\times \times \mathrm{GL}_2(\widehat{O}_F)$ (where we identify $(\nu,u)^*\circ(\nu',u')^*$ with $(\nu\nu',uu')^*$ via the natural isomorphism resulting from the equality $(\nu',u')\circ(\nu,u) = (\nu\nu',uu')$).  It follows that the induced $O_F\otimes O_S$-linear isomorphisms
$$\alpha_{\nu,u}^*:  (\nu,u)^*(s_*\Omega^1_{A/S}) \to s_*\Omega^1_{A/S},\qquad
(\nu,u)^*(R^1s_*\Omega^\bullet_{A/S}) \to R^1s_*\Omega^\bullet_{A/S}$$
satisfy the relation $\alpha_{\nu,u}^*\circ(\nu,u)^*\alpha_{\nu',u'}^* = \alpha_{\nu\nu',uu'}^*$. We thus obtain an action of the group $O_{F, +}^\times \times \mathrm{GL}_2(\widehat{O}_F) $ on the sheaves 
$s_*\Omega^1_{A/S}$ and $R^1s_*\Omega^\bullet_{A/S}$, and hence on the line bundles $\mathcal{L}_{J,N}^{k,l}$, compatible with its action on $Y_{J,N}$.  

Recall that if $\mu \in O_F^\times$, then $(\mu^2,\mu I_2)$ acts trivially on $Y_{J,N}$. In this case the isomorphism $\alpha_{\mu^2,\mu I_2}$ is given by $\iota(\mu)$, and it follows that the induced action on $\mathcal{L}_{J,N}^{k,l}$ is multiplication by the element $\mu^{k+2l} := \prod_{\tau} \tau(\mu)^{k_\tau+2l_\tau}$. In particular if $k_\tau + 2l_\tau$ is an integer $w$ independent of $\tau$, then $\mu^{k+2l} = 
\mathrm{Nm}_{F/\mathbb{Q}}(\mu)^w$.  Thus if %
$w$ is even, then the action of $O_{F,+}^\times \times U$ on $\mathcal{L}_{J,N}^{k,l}$ factors through $G_{U,N}$ and hence defines descent data;
we let $\mathcal{L}^{k,l}_U$ the resulting line bundle on $Y_U$ (given by \cite[Cor.VIII.1.3]{SGA1}).  The same holds if $w$ is odd and $\mathrm{Nm}_{F/\mathbb{Q}}(\mu) = 1$ for all $\mu \in O_F^\times \cap U$. Note that the line bundle $\mathcal{L}_U^{k,l}$ is independent of the choice of $N$.

\subsection{Hilbert modular forms}

\begin{definition}\label{def:paritious} For two tuples $k$ and $l$ above, we say $(k,l)$ is {\em paritious} if $k_\tau+2l_\tau$ is independent of $\tau$.  For such $(k,l)$, we call an element of $H^0(Y_U, \mathcal{L}_U^{k,l})$ a {\em Hilbert modular form} of {\em weight} $(k, l)$ and of {\em level} $U$ (where in addition to Caveat~\ref{def:neat}, we assume that $\mathrm{Nm}_{F/\mathbb{Q}}(\mu) = 1$ for all $\mu \in O_F^\times \cap U$ if $k_\tau + 2 l_\tau$ is odd).\end{definition}

We now make an observation critical to our consideration of weights of mod $p$ Hilbert modular forms. Let $\overline{Y}_{J,N}$ denote the special fibre of $Y_{J,N}$, and similarly let $\overline{\mathcal{L}}_{J,N}^{k,l}$ denote the pull-back of $\mathcal{L}_{J,N}^{k,l}$ to $\overline{Y}_{J,N}$.  If $\mu^{k+2l} \equiv 1 \bmod \pi$ for all $\mu \in O_F^\times \cap U$, then the action of $O_{F,+}^\times \times U$ on $\overline{\mathcal{L}}_{J,N}^{k,l}$ factors through $G_{U,N}$, and hence defines descent data, giving rise to a line bundle $\overline{\mathcal{L}}^{k,l}_U$ on the special fibre $\overline{Y}_U$ of $Y_U$ (again independent of the choice of $N$).  If $(k,l)$ is paritious, then this is simply the pull-back of $\mathcal{L}_U^{k,l}$ to $\overline{Y}_U$, but the line bundles $\overline{\mathcal{L}}^{k,l}_U$ may be defined even if $(k,l)$ is not paritious.
In particular if $O_F^\times \cap U$ is contained in the kernel of reduction modulo $p$, then $\mathcal{L}_U^{k,l}$ is defined for {\em all} pairs $(k,l)$.  This holds for example if $U \subset U_1(\mathfrak{n})$ for some ideal $\mathfrak{n}$ such that the kernel of
$O_F^\times \to (O_F/\mathfrak{n})^\times$ is contained in the kernel of $O_F^\times \to (O_F/p)^\times$.

More generally for any $\mathcal{O}$-algebra $R$ in which the image of $\mu^{k+2l}$ is
trivial for all $\mu \in O_F^\times \cap U$, we obtain a line bundle  $\mathcal{L}_{U,R}^{k,l}$
on $Y_{U,R} = Y_U \times_{\mathcal{O}} R$ by descent from the pull-back of the line bundles
$\mathcal{L}^{k,l}_{J,N}$.

\begin{definition} If $U$, $k$, $l$ and $R$ are such that the image of $\mu^{k+2l}$ in $R$ is trivial for all $\mu \in O_F^\times \cap U$, then we call an element of $H^0(Y_{U,R}, \overline{\mathcal{L}}_{U,R}^{k,l})$ a {\em Hilbert modular form over $R$}
of {\em weight} $(k, l)$ and {\em level} $U$, and we write $M_{k,l}(U;R)$ for the $R$-module
of  such forms.  If $R=E$, then we call such a form
a {\em mod $p$ Hilbert modular form} (of {\em weight} $(k, l)$ and of {\em level} $U$).\end{definition}

\begin{definition}
\label{def:pneat}
We say that $U$ is {\em $p$-neat} if $O_F^\times \cap U$ is contained in the kernel of reduction modulo $p$ (in addition to $U$ being sufficiently small in the sense of Caveat~\ref{def:neat}). \end{definition} 

\subsection{The Koecher Principle}

The Koecher Principle implies that $M_{k,l}(U;R)$ is a finitely generated $R$-module (assuming an $\mathcal{O}$-algebra $R$
is Noetherian), and that $M_{0,0}(U;R) =  H^0(Y_{U,R}, \mathcal{O}_{Y_{U,R}})$ is the set of locally constant
functions on $Y_{U,R}$.  Both of these assertions follow from the analogous ones with $Y_U$ replaced
by $Y_{J,N}$, proved by Rapoport. (The case $J=O_F$ is treated by Prop.~4.9 and the discussion preceding
Prop.~6.11 of \cite{R}, and the modifications needed for the case of arbitrary $J$ are given in \cite{C2};
see also \cite[Thm.~8.3]{Di} and \cite[Thm.~7.1]{DT} for variants with different level structure and descent data
in place.)

\subsection{Canonical trivialisations} \label{subsec:cantriv}
We observe that the sheaves $\delta^{\otimes l}$ on $Y_{J,N}$ are in fact free (not just locally so).
Indeed if $A$ is the universal HBAV over $S=Y_{J,N}$, then we have a sequence of canonical
isomorphisms:
\begin{equation}\label{eqn:duality}
\begin{array}{c}
R^1s_*\mathcal{O}_A \simeq \mathrm{Lie}(A^\vee) \simeq \mathrm{Lie}(A) \otimes_{O_F} J
 \simeq \mathcal{H}om_{\mathcal{O}_S}(s_*\Omega^1_{A/S}, \mathcal{O}_S)\otimes_{O_F} J\\
  \simeq \mathcal{H}om_{O_F\otimes\mathcal{O}_S}(s_*\Omega^1_{A/S}, J\mathfrak{d}^{-1}\otimes \mathcal{O}_S),
\end{array}
\end{equation}
  from which it follows that $\wedge^2_{O_F \otimes \mathcal{O}_S}\mathcal{H}^1_{\mathrm{DR}}(A/S)
  \cong s_*\Omega^1_{A/S} \otimes_{O_F \otimes \mathcal{O}_S} R^1s_*\mathcal{O}_A$
  is canonically isomorphic to $J\mathfrak{d}^{-1}\otimes \mathcal{O}_S$, which is free of rank one
  over $O_F\otimes \mathcal{O}_S$. Therefore each $\delta_\tau$ is free of rank one over $\mathcal{O}_S$, and hence so are the sheaves $\delta^{\otimes l}$.
  
Under the action of  $(\nu,u) \in O_{F,+}^\times \times U$ on $Y_{J,N}$, one finds that
the canonical isomorphism $\psi$ from $\wedge^2_{O_F \otimes \mathcal{O}_S}\mathcal{H}^1_{\mathrm{DR}}(A/S)$ to $J\mathfrak{d}^{-1}\otimes \mathcal{O}_S$ is multiplied by
$\nu$ (in the sense that $(\nu,u)^*\psi = (\nu\otimes 1)\psi \circ \alpha_{\nu,u}^*$).  Therefore the
action of $(\nu,u)$ on the resulting trivialisation of $\mathcal{L}_{J,N}^{0,l} = \delta^{\otimes l}$ is
multiplication by $\nu^l$. In particular, if $(0,l)$ is paritious (i.e., $l_\tau$ is independent of $\tau$),
then $\nu^l = 1$ so the trivialisation of $\mathcal{L}_{J,N}^{0,l}$ on $Y_{J,N}$ is invariant under $G_{U,N}$, hence descends to one on $Y_U$.
Similarly if $\nu^l \equiv 1 \bmod \pi$ for all $\nu \in O_{F,+}^\times$, then the canonical trivialisation of 
$\overline{\mathcal{L}}_{J,N}^{0,l}$ on $\overline{Y}_{J,N}$ descends to one on $\overline{Y}_U$.

\subsection{Complex Hilbert modular forms}
If $(k,l)$ is paritious, then under the identification (\ref{eqn:chmv}), the line bundle $\mathcal{L}_U^{k,l}$
gives the usual automorphic line bundle whose sections are classical Hilbert modular forms of weight $(k,l)$ and level $U$.  More precisely, $\mathcal{L}_U^{k,l}$ is defined over $\mathcal{O}\cap \overline{\mathbb{Q}}$, and its fibre at the point $y_{g_\infty,x} \in Y_U(\mathbb{C})$ corresponding to the double coset $\mathrm{GL}_2(F)g_\infty \mathrm{diag}(1,x) UU_\infty$ has basis $e^{k,l} = \otimes_\tau (ds_\tau^{\otimes k_\tau} \otimes h_\tau^{\otimes l_\tau})$, where $s= (s_\tau)_{\tau\in\Sigma}$ are the
coordinates on $\mathbb{C}\otimes F \cong \mathbb{C}^\Sigma$ and $h_\tau$ is the basis for $\delta_\tau$
given by the trivialisation defined above.  For $\phi \in M_{k,l}(U;\mathbb{C})$, we define the function
$f_\phi: \mathrm{GL}_2(\mathbb{A}_F) \to \mathbb{C}$ so that
$$y_{g_\infty,x}^* \phi = ||x||^{-1}\det(g_\infty)^{l-1}j(g_\infty,z_0)^k f_\phi(\gamma g_\infty u) e^{k,l}\mbox{\ for all $\gamma \in \mathrm{GL}_2(F)$, $u \in U$,}$$
where $j(g_\infty,z) = cz +d$ for $g_\infty = \begin{pmatrix} a & b \\ c & d \end{pmatrix} \in 
\mathrm{GL}_2^+(\mathbb{R})$ and $z$ in the complex upper-half plane $\mathfrak{H}$, and the
exponents $k$ and $l-1$ denote products over the embeddings $\tau \in \Sigma$.
Then $\phi \mapsto  f_\phi$ defines an isomorphism 
$M_{k,l}(U;\mathbb{C}) \simeq A_{k,l}(U)$, where $A_{k,l}(U)$ is the set of functions 
$f: \mathrm{GL}_2(\mathbb{A}_F) \to \mathbb{C}$ such that:
\begin{itemize}
\item $f(\gamma hu) = \det(u_\infty)^{1-l} j(u_\infty,i)^{-k} f(h)$ for all $\gamma\in \mathrm{GL}_2(F)$, 
$h \in \mathrm{GL}_2(\mathbb{A}_F)$ and $u\in UU_\infty$;
\item $f_h(g_\infty(z_0)) = \det(g_\infty)^{l-1} j(g_\infty,i)^k f(h g_\infty)$ is holomorphic on $\mathfrak{H}^\Sigma$ for all $h \in \mathrm{GL}_2(\mathbb{A}_F^\infty)$.
\end{itemize}
Note also that $f\mapsto (f_{\mathrm{diag}(1,x)})$ defines an isomorphism $A_{k,l}(U) \simeq \oplus_x M_k(\Gamma_{U,x})$, where $x$ runs over a set of representatives of $F^\times \backslash \mathbb{A}_F^\times / \det(U) F_{\infty,+}^\times$, $\Gamma_{U,x} = \mathrm{GL}_2^+(F) \cap \mathrm{diag}(1,x) U\mathrm{diag}(1,x)^{-1}$ and $M_k(\Gamma)$ denotes the set of holomorphic functions
$\varphi:\mathfrak{H}^\Sigma \to \mathbb{C}$ such that $\varphi(\gamma(z)) = 
\det(\gamma)^{-k/2}j(\gamma,z)^k \varphi(z)$ for all $\gamma \in \Gamma$.

\subsection{Forms of weight $(0,l)$ in characteristic $p$}

Let us now return to characteristic $p$ and
give sufficient hypotheses for the sheaf $\overline{\mathcal{L}}_U^{0,l}$ on the special fibre
$\overline{Y}_U$ to be globally free, even when $(0,l)$ is not paritious.  Suppose that
$\mu_N(\overline{\mathbb{Q}}) \subset \mathcal{O}$, so the geometric components
of $Y_{J,N}$ are defined over $\mathcal{O}$. Recall that the set of geometric components
is in bijection with $Z_{J,N}(\mathcal{O})$, with $(\nu,u)$ acting by $\nu\det(u)^{-1}$, so the
stabiliser of each component of $Y_{J,N}$ is
$\{\,(\nu,u) \in O_{F,+}^\times \times U\,|\,\nu\equiv \det(u)\bmod N\,\}$.
Letting $H_{U,N}$ denote the corresponding subgroup of $G_{U,N}$, we see that
if $\nu^l \equiv 1 \bmod \pi$ for each $\nu \in O_{F,+}^\times \cap \det(U)$, then
the trivialisation of $\overline{\mathcal{L}}_{J,N}^{0,l}$ on $\overline{Y}_{J,N}$ is invariant
under $H_{U,N}$, so descends to the quotient $\overline{Y}_{J,N}/H_{U,N}$.
Note that this hypothesis also implies
that $\mu^{2l} \equiv 1 \bmod \pi$ for all $\mu \in O_F^\times \cap U$, so that
$\overline{\mathcal{L}}_{J,N}^{0,l}$ descends to $\overline{Y}_U$; since the projection
from $\coprod\overline{Y}_{J,N}/H_{U,N}$ is an isomorphism on each connected component,
it follows that $\overline{\mathcal{L}}_U^{0,l}$ is (globally) free on $\overline{Y}_U$.

We record this as follows (recall that $Z_U$ is defined in see \S\ref{subsec:components}):
\begin{proposition} \label{prop:twist}
Suppose that $\mu_N(\overline{\mathbb{Q}}) \subset \mathcal{O}$
for some $N$ prime to $p$ such that $U(N) \subset U$.
If $\nu^l \equiv 1 \bmod \pi$ for all $\nu \in O_{F,+}^\times \cap \det(U)$, then
the sheaf $\overline{\mathcal{L}}_U^{0,l}$ on $\overline{Y}_U$ is
(non-canonically) isomorphic to $\mathcal{O}_{\overline{Y}_U}$, and $M_{0,l}(U;E)$ to
the space of functions $Z_U(\mathcal{O}) \to E$.
\end{proposition}

Note that the hypotheses of the proposition are satisfied for {\em all} $l \in \mathbb{Z}^\Sigma$
if  $O_{F,+}^\times \cap \det(U)$ is contained in the kernel of reduction modulo $p$.
This holds for example if  $U \subset {}^1U_1(\mathfrak{n})$ for some ideal $\mathfrak{n}$
such that the kernel of $O_F^\times \to (O_F/\mathfrak{n})^\times$ is contained in the kernel
of $O_F^\times \to (O_F/p)^\times$.  In this case $U$ is also $p$-neat, so the sheaves
$\overline{\mathcal{L}}_U^{k,l}$ are defined for all pairs $(k,l)$, and the spaces of mod $p$ Hilbert modular forms $H^0(\overline{Y}_U, \overline{\mathcal{L}}_U^{k,l})$ for fixed $k$ and varying $l$
are (non-canonically) isomorphic.

\section{Hecke operators} \label{section:Hecke} 

In this section, we define Hecke operators geometrically on spaces of mod $p$ Hilbert modular forms. %

\subsection{Adelic action on Hilbert modular varieties}
\label{subsec:adelicact1}

Suppose that $U_1$ and $U_2$ are open compact subgroups of $\mathrm{GL}_2(\widehat{O}_F)$;
we assume as usual that Caveat~\ref{def:neat} holds, so $U_1$ and $U_2$ contain $\mathrm{GL}_2(O_{F,p})$
and are sufficiently small in the sense that the conclusion of Lemma~\ref{lemma:neat} holds.

Suppose that $g \in \mathrm{GL}_2(\mathbb{A}_F^\infty) = \mathrm{GL}_2(\widehat{O}_F\otimes\mathbb{Q})$ with $g_p \in \mathrm{GL}_2(O_{F,p}) $ and $g^{-1}U_1 g \subset U_2$.  We now proceed to define a morphism
$\rho_g: Y_{U_1} \to Y_{U_2}$ which corresponds to the one defined by right multiplication by $g$ on the associated Hilbert modular varieties; i.e. on complex points it is given by 
$\mathrm{GL}_2(F)xU_1U_\infty \mapsto \mathrm{GL}_2(F)xgU_2U_\infty$.

We first choose:
\begin{itemize}
\item $\alpha \in O_F$ such that $\alpha g \in M_2(\widehat{O}_F)$ and $\alpha \in O_{F,p}^\times$;
\item $N_2$ prime to $p$ such that $U(N_2) \subset U_2$;
\item $N_1$ prime to $p$ such that $U(N_1) \subset U_1$ and 
         $(\alpha g)^{-1}N_1/N_2 \in M_2(\widehat{O}_F)$.
\end{itemize}

We will define a morphism $\tilde{\rho}_g:\coprod Y_{J,N_1} \to \coprod Y_{J,N_2}$ whose composite
with the projection to $Y_{U_2}$ factors through $Y_{U_1}$, yielding the desired morphism
$\rho_g:Y_{U_1} \to Y_{U_2}$, independent of the above choices of $\alpha$, $N_1$ and $N_2$.

We first note that the conditions above imply that $N_2|N_1$,
$g^{-1}U(N_1)g \subset U(N_2)$, and (right) multiplication
by $(\alpha g)^{-1}N_1/N_2$ induces an injective $O_F$-linear map
$$j: (O_F/N_2)^2 \longrightarrow (O_F/N_1)^2 /(O_F/N_1)^2 \cdot(\alpha g)^{-1}N_1 .$$

Let $(A_1,\iota_1,\lambda_1,\eta_1)$ denote the universal HBAV over $S=Y_{J_1,N_1}$ where
$J_1 = J_{t_1}$ for some $t_1 \in T$, and let
$A_1' = A_1/\eta_1(C)$ where $C =(O_F/N_1)^2\cdot(\alpha g)^{-1}N_1 $.  Then $A_1'$ inherits
an $O_F$-action $\iota_1'$ from $A_1$, and $\eta_1$ induces an $O_F$-linear closed immersion
$$(O_F/N_1)^2 / (O_F/N_1)^2\cdot(\alpha g)^{-1}N_1 \longrightarrow A_1'$$
whose composite with $j$ defines an isomorphism $\eta_1' : (O_F/N_2)^2 \to A_1'[N_2]$.

Now consider the injective $O_F$-linear map $\pi^*: \mathrm{Sym}(A_1'/S) \to \mathrm{Sym}(A_1/S)$
defined by $f \mapsto \pi^\vee \circ f \circ \pi$, where $\pi$ is the natural projection $A_1 \to A_1'$.
\begin{lemma} \label{lemma:pol} The image of $\pi^*$ is
$(\det (\alpha g))\mathrm{Sym}(A_1/S)$ where $(\det(\alpha g))$ denotes
the ideal $O_F \cap \det(\alpha g)\widehat{O}_F$ of $O_F$.
\end{lemma}
\begpf Note that since $\mathrm{Sym}(A_1/S)$ is an invertible $O_F$-module, the image
of $\pi^*$ is (locally on $S$) of the form $I\mathrm{Sym}(A_1/S)$ for some ideal $I$ of $O_F$,
non-zero since $\pi^*$ is injective.  Moreover since $\ker(\pi) \subset A_1[N_1]$, there is an isogeny
$\phi:A_1' \to A_1$ such that $\phi\circ \pi$ is multiplication by $N_1$; since $\pi^*\circ \phi^*$ is multiplication
by $N_1^2$, it follows that $N_1^2 \in I$, so $I$ can only be divisible by primes dividing $N_1$.

We now determine $I \otimes \mathbb{Z}_\ell$ for each prime $\ell | N_1$. 
Note in particular that $\ell \neq p$, so $\ell$ is invertible in $\mathcal{O}_S$.
Consider the commutative
diagram:
$$\begin{array}{ccccc}
\mathrm{Sym}(A_1'/S) \otimes \mathbb{Z}_\ell  && \longrightarrow && \mathrm{Sym}(A_1/S) \otimes \mathbb{Z}_\ell \\ &&&& \\ \downarrow \wr &&&&  \downarrow \wr \\&&&&\\
\mathrm{Hom}_{\mathbb{Z}_\ell} (\wedge^2_{O_{F,\ell}}T_\ell(A_1'),\mathbb{Z}_\ell(1))
&&\longrightarrow &&
\mathrm{Hom}_{\mathbb{Z}_\ell} (\wedge^2_{O_{F,\ell}}T_\ell(A_1),\mathbb{Z}_\ell(1)),
\end{array}$$
of $O_F$-linear maps of $\ell$-adic sheaves on $S$, where the top map is 
$\pi^*\otimes \mathbb{Z}_\ell$, the vertical isomorphisms are induced by the Weil pairings,
and the bottom map is given by the map $T_\ell(\pi): T_\ell(A_1) \to T_\ell(A_1')$ on
$\ell$-adic Tate modules induced by $\pi$.
The cokernel of $\pi^* \otimes \mathbb{Z}_\ell$ is therefore isomorphic to that of the bottom map,
which in turn is isomorphic to $\mathrm{Hom}_{\mathbb{Z}_\ell} (M_\ell,\mathbb{Q}_\ell/\mathbb{Z}_\ell(1))$, where $M_\ell$ is the cokernel of $\wedge^2_{O_{F,\ell}}T_\ell(\pi)$.
Since the $\ell$-adic sheaves $T_\ell(A_1)$ and $T_\ell(A_1')$ are locally free of rank two over
$O_{F,\ell}$ and the cokernel of $T_\ell(\pi)$ is isomorphic to 
$$\ker(\pi)\otimes \mathbb{Z}_\ell \cong C \otimes \mathbb{Z}_\ell \cong 
 O_{F,\ell}^2\cdot(\alpha g)^{-1} /O_{F,\ell}^2,$$
 it follows that $M_\ell$ is isomorphic to $O_{F,\ell}/\det(\alpha g)O_{F,\ell}$.
 
 We have now shown that the cokernel of $\pi^*\otimes \mathbb{Z}_\ell$ is (\'etale locally)
 isomorphic to $O_{F,\ell}/\det(\alpha g)O_{F,\ell}$ for all $\ell$.
Since the cokernel of $\pi^*$ is also \'etale locally isomorphic $O_F/I$, it follows that
$O_F/I$ is isomorphic to $\widehat{O}_F/\det(\alpha g)\widehat{O}_F$, and hence
that $I= (\det(\alpha g))$.
\epf

It follows from the lemma that $\lambda_1:J_1 \simeq \mathrm{Sym}(A_1/S)$ restricts
to an isomorphism $IJ_1 \to \pi^*\mathrm{Sym}(A_1'/S)$, where $I =  (\det(\alpha g))$.
Moreover since $f$ is a section of $\mathrm{Pol}(A_1'/S)$ if and only if $\pi^*f$ is
a section of $\mathrm{Pol}(A_1/S)$, we see that $\lambda_1$ further restricts to an
isomorphism $(IJ_1)_+ \to \pi^*\mathrm{Pol}(A_1'/S)$.
Now let $J_2= J_{t_2}$ where $t_2 \in T$ is the fixed representative of $IJ_1$ in the strict
class group of $F$, and choose an element $\beta \in F_+^\times$ such that
$\beta J_2 = IJ_1$.  Thus $\beta$ is uniquely determined up to $O_{F,+}^\times$,
and the composite of $\lambda_1\circ \beta$ with the inverse of $\pi^*$ yields
an isomorphism
\begin{eqnarray*}
(J_2,(J_2)_+) \simeq (IJ_1, (IJ_1)_+)& \simeq &\pi^*(\mathrm{Sym}(A_1'/S), \mathrm{Pol}(A_1'/S))\\
 &\simeq& (\mathrm{Sym}(A_1'/S), \mathrm{Pol}(A_1'/S))\end{eqnarray*}
which we denote by $\lambda_1'$. 
 
Finally we note that since $A$ satisfies the Deligne--Pappas condition, so does $A'$.
This follows for example from the commutative diagram:
$$\begin{array}{ccccc}
A_1 \otimes_{O_F} IJ_1 & \to&  A_1 \otimes_{O_F} J_1 &\simeq& A_1^\vee\\
&&&&\\
{\scriptstyle{\pi\otimes 1}}\downarrow &&&& \uparrow{\scriptstyle{\pi^\vee}}\\
&&&& \\
 A_1' \otimes_{O_F} IJ_1 &\multicolumn{3}{c}{\longrightarrow}& (A_1')^\vee,
\end{array}$$
and the observation that the top left map is an isogeny with kernel $A_1[I]\otimes_{O_F} IJ_1$,
hence (constant) degree $|O_F/I|^2$, while $\deg(\pi\otimes 1) = \deg(\pi^\vee) = \deg(\pi) = |O_F/I|$, so the bottom map must be an isomorphism.

Now $(A_1',\iota_1',\lambda_1',\eta_1')$ is a $J_2$-polarised HBAV with level $N_2$ structure
over $Y_{J_1,N_1}$, so corresponds to a morphism $Y_{J_1,N_1} \to Y_{J_2,N_2}$ such
that the pull-back of the universal HBAV over $Y_{J_2,N_2}$ is $(A_1',\iota_1',\lambda_1',\eta_1')$.
Taking the union over $t_1 \in T$ yields the desired morphism $\tilde{\rho}_g:\coprod Y_{J,N_1} \to \coprod Y_{J,N_2}$.

It is straightforward to check that the composite of $\tilde{\rho}_g$
with the projection to $Y_{U_2}$ is independent of the choices of $\alpha$, $N_2$ and $\beta$,
and indeed of $N_1$ in the sense that if $N_1$ is replaced
by a multiple $N$, then the resulting morphism is obtained by composing with the natural
projection $\coprod Y_{J,N} \to \coprod Y_{J,N_1}$.  (The only non-trivial point is that if
$\alpha$ is replaced by a multiple $\delta\alpha$, then the resulting $J_2$-polarised HBAV
with level $N_2$ structure on $Y_{J_1,N_1}$ is isomorphic to the original 
$(A_1',\iota_1',\lambda_1',\eta_1')$ via the
map induced by $\iota_1(\delta)$.)  Moreover the resulting morphism to $Y_{U_2}$ is invariant
under the action of $G_{U_1,N_1}$ on  $\coprod Y_{J,N_1}$ (indeed we have
 $\tilde{\rho}_g \circ (\mu,u) = (\mu,g^{-1}ug) \circ \tilde{\rho}_g$ for all $(\mu,u) \in O_{F,+}^\times
 \times U_1$ on each $Y_{J_1,N_1}$ for any choice of $\beta$ as above),
hence factors through $Y_{U_1}$,
yielding the desired morphism $\rho_g: Y_{U_1} \to Y_{U_2}$.

Suppose that $U_1$, $U_2$ and $U_3$ are open compact subgroups of $\mathrm{GL}_2(\widehat{O}_F)$
with $g_1,g_2 \in \mathrm{GL}_2(\mathbb{A}_F^\infty)$ as above satisfying  $g_1^{-1}U_1 g_1
 \subset U_2$
and $g_2^{-1} U_2 g_2 \subset U_3$, so that $\rho_{g_1}: Y_{U_1} \to Y_{U_2}$ and 
$\rho_{g_2}: Y_{U_2} \to Y_{U_3}$ are defined.  Note that choosing $\alpha_2$, $N_2$
and $N_3$ to define $\rho_{g_2}$, and then $\alpha_1$, $N_1$ and (the same) $N_2$ to define
$\rho_{g_1}$, we may use $\alpha_1\alpha_2$, $N_1$ and $N_3$ to define $\rho_{g_1g_2}$.
Let  $(A_i,\iota_i,\lambda_i,\eta_i)$ denote the universal HBAV over $Y_{N_i,J_i}$ for $i=1,2,3$,
where $J_i = J_{t_i}$ for $t_i \in T$ such that $t_{i+1}$ represents the class of
$(\det (\alpha_i g_i))J_i$ for $i=1,2$.  The above construction of $\tilde{\rho}_{g_i}$
then yields a $J_{i+1}$-polarised abelian variety $(A_i',\iota_i',\lambda_i',\eta_i')$
with level $N_{i+1}$ structure over $Y_{J_i,N_i}$, where $A_i' = A_i/\eta_i(C_i)$ with
$C_i = (O_F/N_i)^2\cdot (\alpha_i g_i)^{-1}N_i$. It is straightforward to check that
the pull-back via $\tilde{\rho}_{g_1}$ of $(A_2',\iota_2',\lambda_2',\eta_2')$ is isomorphic 
to a $J_3$-polarised HBAV with level $N_3$-structure defining $\tilde{\rho}_{g_1g_2}$,
so that we may take $\tilde{\rho}_{g_1g_2} = \tilde{\rho}_{g_2} \circ \tilde{\rho}_{g_1}$
and conclude that $\rho_{g_1g_2} = \rho_{g_2} \circ \rho_{g_1}$.  

\subsection{Adelic action on Hilbert modular forms}\label{subsec:adelicaction}

We revert to the original setting of \S\ref{subsec:adelicact1}, with $g$, $U_1$ and $U_2$
satisfying $g^{-1}U_1g \subset U_2$, and use the notation in the
definition of $\rho_g$ (and in particular a choice of $N_1$, $N_2$, $\alpha$ and $\beta$),
but writing $S_i = Y_{J_i,N_i}$ for $i=1,2$ and $s_i:A_i \to S_i$ and $s_1':A_1' \to S_1$
for the structural morphisms.  We let $\pi_\alpha$ denote the canonical projection
$A_1 \to A_1' \simeq \tilde{\rho}_g^* A_2$; the dependence on $\alpha$ is such that
if $\delta \in O_F \cap O_{F,p}^\times$ (and $N_1$ is such that $(\delta\alpha g)^{-1}N_1/N_2
\in M_2(\widehat{O}_F)$), then $\pi_{\delta\alpha} = i(\delta)\pi_\alpha$.
It follows that the $O_F\otimes \mathcal{O}_{S_1}$-linear morphisms 
\begin{equation}
\label{eqn:action}
\begin{array}{ccccc}
 \tilde{\rho}_g^* s_{2,*} \Omega^1_{A_2/S_2}& \simeq &s_{1,*}' \Omega^1_{A_1'/S_1}
 & \to   & s_{1,*}\Omega^1_{A_1/S_1},\\&&&&\\
  \tilde{\rho}_g^* R^1s_{2,*} \Omega^\bullet_{A_2/S_2} &\simeq &R^1s_{1,*}' \Omega^\bullet_{A_1'/S_1}
&  \to  & R^1s_{1,*}\Omega^\bullet_{A_1/S_1}\end{array}\end{equation}
induced by  $(\alpha\otimes 1)^{-1}\pi_\alpha^*$ are independent of the choice of $\alpha$
(as well as $N_2$ and $\beta$, and even $N_1$ in the sense of compatibility with pull-back
by the natural projection). Note that these are in fact isomorphisms since the degree of
the isogeny $\pi_\alpha$ is invertible in $\mathcal{O}_S$.  Furthermore the commutativity of the diagram:
  $$\begin{array}{ccccc}
  A_1 &\multicolumn{3}{c}{\stackrel{\alpha_{\nu,u}}{\longrightarrow}} & (\nu,u)^*A_1 \\&&&&\\
 {\scriptstyle{\pi_\alpha}} \downarrow &&&& \downarrow {\scriptstyle{(\nu,u)^*\pi_\alpha} }\\
 &&&&\\
 \tilde{\rho}_g^*A_2 & \stackrel{\tilde{\rho}_g^*(\alpha_{\nu,g^{-1}ug})}{\longrightarrow} &
 \tilde{\rho}_g^*(\nu,g^{-1}ug)^*A_2 & \simeq & (\nu,u)^*\tilde{\rho}_g^* A_2 \end{array}$$ 
 implies that the isomorphisms in (\ref{eqn:action}) are compatible with the action of 
 $G_{U_1,N_1}$ (where $G_{U_1,N_1}$ acts on the sources via the homomorphism
 $(\nu,u) \mapsto (\nu,g^{-1}ug)$ to $G_{U_2,N_2}$ and pull-back by $\tilde{\rho}_g^*$).
 It follows that the same is true for the $\mathcal{O}_{S_1}$-linear isomorphisms 
 $\tilde{\rho}_g^*\mathcal{L}_{J_2,N_2}^{k,l} \stackrel{\sim}{\longrightarrow} \mathcal{L}_{J_1,N_1}^{k,l}$
 induced by those in (\ref{eqn:action}) for $k,l \in \mathbb{Z}^\Sigma$, which therefore
descend to define isomorphisms 
\begin{equation}\label{eqn:pullback}
\rho_g^*\mathcal{L}^{k,l}_{U_2,R} \stackrel{\sim}{\longrightarrow} \mathcal{L}^{k,l}_{U_1,R}\end{equation}
for any $\mathcal{O}$-algebra $R$ in which the image of $\mu^{k+2l}$ is trivial
for all $\mu \in O_F^\times \cap U_2$ (and hence all $\mu \in O_F^\times\cap U_1$).
We thus obtain an $R$-linear map $[U_1 g U_2]: M_{k,l}(U_2;R) \to M_{k,l}(U_1;R)$
defined as the product of $||\det(g)|| = \mathrm{Nm}_{F/\mathbb{Q}}(\det g)^{-1}$ with the
composite:
$$H^0(Y_{U_2,R},\mathcal{L}^{k,l}_{U_2,R}) 
\longrightarrow H^0(Y_{U_1,R}, \rho_g^*\mathcal{L}^{k,l}_{U_2,R})
\stackrel{\sim}{\longrightarrow} H^0(Y_{U_1,R},\mathcal{L}^{k,l}_{U_1,R}).$$

Returning now to the setting where $U_1$, $U_2$ and $U_3$ are open compact subgroups of $\mathrm{GL}_2(\widehat{O}_F)$ and 
$g_1,g_2 \in \mathrm{GL}_2(\mathbb{A}_F^\infty)$ are such that $g_1^{-1}U_1 g_1 \subset U_2$
and $g_2^{-1} U_2 g_2 \subset U_3$, we find that the composite:
$$\label{eqn:composite}
 A_1 \stackrel{\pi_{\alpha_1}}{\longrightarrow} \tilde{\rho}_{g_1}^* A_2 
 \stackrel{\tilde{\rho}_{g_1}^*\pi_{\alpha_2}}{\longrightarrow}  \tilde{\rho}_{g_1}^*\tilde{\rho}_{g_2}^* A_3
 \simeq \tilde{\rho}_{g_1g_2}^* A_3$$
 is $\pi_{\alpha_1\alpha_2}$.  This in turn implies that the composite
$$\rho_{g_1g_2}^*\mathcal{L}_{U_3,R}^{k,l} \simeq
     \rho_{g_1}^*\rho_{g_2}^*\mathcal{L}_{U_3,R}^{k,l}
      \stackrel{\sim}{\longrightarrow} \rho_{g_1}^*\mathcal{L}_{U_2,R}^{k,l} 
      \stackrel{\sim}{\longrightarrow} \mathcal{L}_{U_1,R}^{k,l},$$
is the isomorphism in (\ref{eqn:pullback}) used to define $ [U_1 g_1g_2 U_3]$, which therefore
coincides with $[U_1 g_1U_2] \circ [U_2 g_2 U_3]$.

For $R=\mathcal{O}$ and $(k,l)$ paritious, we thus obtain an action of the group
$\{\,g\in \mathrm{GL}_2(\mathbb{A}_F^\infty)\,|\, g_p \in \mathrm{GL}_2(O_{F,p})\,\}$ on
$$M_{k,l}(\mathcal{O}) :=   \varinjlim M_{k,l}(U;\mathcal{O}),$$
where the direct limit is over all sufficiently small open compact subgroups
$U$ of $\mathrm{GL}_2(\mathbb{A}_F^\infty)$ containing $\mathrm{GL}_2(O_{F,p})$.
Similarly we have an action on $M_{k,l}(\mathbb{C}) :=   \varinjlim M_{k,l}(U;\mathbb{C})$, which
is compatible by extension of scalars with the one just defined on $M_{k,l}(\mathcal{O})$.  One can check
that the action is also compatible under the isomorphisms $M_{k,l}(U;\mathbb{C}) \simeq A_{k,l}(U)$
with the usual action defined by right multiplication on the space of automorphic forms
$A_{k,l} :=   \varinjlim A_{k,l}(U)$.

Recall that for $R = E$ and arbitrary $(k,l)$, the space $M_{k,l}(U;E)$ is
defined for sufficiently small $U$ (for example $p$-neat as in Definition~\ref{def:pneat}), so we may similarly define
$$M_{k,l}(E) :=   \varinjlim M_{k,l}(U;E).$$
Then $M_{k,l}(E)$ is a smooth admissible representation of
$\{\,g\in \mathrm{GL}_2(\mathbb{A}_F^\infty)\,|\, g_p \in \mathrm{GL}_2(O_{F,p})\,\}$
over $E$, and we recover $M_{k,l}(U;E) = M_{k,l}(E)^U$ for sufficiently small $U$
containing $\mathrm{GL}_2(O_{F,p})$.  (Note that $M_{k,l}(E)^U = 0$ if 
$\overline{\mu}^{k+2l}\neq 1$ for some $\mu \in U \cap O_F^\times$.)

We may similarly define $M_{k,l}(R)$ for any $(k,l)$ and $R$ in which $p$ is nilpotent.
We again have $M_{k,l}(U;R) = M_{k,l}(R)^U$ for sufficiently small $U$ (indeed for any
$U$ for which we have already defined $M_{k,l}(U;R)$), so we may define $M_{k,l}(U;R)$
to be $M_{k,l}(R)^U$ for {\em any} open compact subgroup $U$ of 
$\mathrm{GL}_2(\mathbb{A}_F^\infty)$ containing $\mathrm{GL}_2(O_{F,p})$.
Note then that $M_{k,l}(U;R) = 0$ if $\overline{\mu}^{k+2l}\neq 1$ for some
$\mu \in U \cap O_F^\times$, but not necessarily under the weaker assumption (if $pR \neq 0$)
that $\mu^{k+2l}$ has non-trivial image in $R$ for some $\mu \in U \cap O_F^\times$.
We shall restrict our attention however to the case $R=E$.

\subsection{Hecke operators}
Suppose now that $U_1$ and $U_2$ are open compact subgroups of 
$\mathrm{GL}_2(\mathbb{A}_F^\infty)$ containing $\mathrm{GL}_2(O_{F,p})$ and
that $g$ is an element of $\mathrm{GL}_2(\mathbb{A}_F^\infty)$ 
such that $g_p \in \mathrm{GL}_2(O_{F,p})$.  We may then define the double coset
operator 
$$[U_1 g U_2]: M_{k,l}(U_2;E) \to M_{k,l}(U_1;E)$$
 to be the map 
$f \mapsto \sum_{i\in I} g_i f$ where  $U_1 g U_2 = \coprod_{i\in I} g_i U_2$.
It is straightforward to check that the map is independent of the choice of 
representatives $g_i$, that the image is indeed in $M_{k,l}(U_1;E)$, and that
the definition agrees with the one already made when $U_1$ and $U_2$ are
sufficiently small and $g^{-1}U_1 g \subset U_2$.  (Recall that a normalising
factor of $||\det(g)||$ is incorporated into the definition of the action.)

If $U_1$ and $U_2$ are sufficiently small we may reinterpret $[U_1 g U_2]$
in the usual way using trace morphisms as follows.  Letting $U_1' = U_1 \cap g U_2 g^{-1}$,
we have that $g^{-1} U_1' g \subset U_2$, so that $[U_1 g U_2] = [U_1 1U_1'] \circ [U_1'gU_2]$
and $[U_1'gU_2]$ is the composite
$$H^0(\overline{Y}_{U_2},\overline{\mathcal{L}}_{U_1'}^{k,l}) 
\to H^0(\overline{Y}_{U_1'}, \rho_g^*\overline{\mathcal{L}}_{U_2}^{k,l}) 
\to H^0(\overline{Y}_{U_1'}, \overline{\mathcal{L}}_{U_1'}^{k,l})$$
where the second map is $||\det(g)||$ times the one induced from (\ref{eqn:pullback}).
On the other hand $[U_1 1 U_1']$ is precisely the composite
$$H^0(\overline{Y}_{U_1'},\overline{\mathcal{L}}_{U_1'}^{k,l}) 
\to H^0(\overline{Y}_{U_1'}, \rho_1^*\overline{\mathcal{L}}_{U_1}^{k,l}) 
\to H^0(\overline{Y}_{U_1}, \overline{\mathcal{L}}_{U_1}^{k,l}),$$
where the first map is given by the inverse of the isomorphism
$\rho_1^* \mathcal{L}_{U_1}^{k,l} \to \mathcal{L}_{U_1'}^{k,l}$
(from (\ref{eqn:pullback})), and the 
last map is the trace times the index of $U_1' \cap O_F^\times$ in
$U_1 \cap O_F^\times$.  

For primes $v$ of $F$ such that $v \nmid p$ and $\mathrm{GL}_2(O_{F,v}) \subset U$,
we define the Hecke operators
\begin{equation}\label{eqn:TvSv} T_v := \left[ U \left(\begin{array}{cc} 1 & 0 \\ 0 & \varpi_v \end{array}\right) U \right]
\quad\mbox{and}\quad
 S_v := \left[ U \left(\begin{array}{cc} \varpi_v & 0 \\ 0 & \varpi_v \end{array}\right) U \right]\end{equation}
on $M_{k,l}(U;E)$, where $\varpi_v$ is a uniformiser of $O_{F,v}$.
These operators are independent of the choice of $\varpi_v$, and commute with
each other (for varying $v$).  Note that under the above interpretation via the trace map
(for sufficiently small $U$), we have $U_1' = U \cap U_0(v)$ and $U_1' \cap O_F^\times
= U \cap O_F^\times$, so that $T_v$ can be written as $\mathrm{Nm}_{F/\mathbb{Q}}(v)^{-1}$
times the composite
$$H^0(\overline{Y}_{U},\overline{\mathcal{L}}_{U}^{k,l}) 
\to H^0(\overline{Y}_{U'}, \rho_g^*\overline{\mathcal{L}}_{U}^{k,l}) 
\to H^0(\overline{Y}_{U'}, \rho_1^*\overline{\mathcal{L}}_{U}^{k,l}) 
\to H^0(\overline{Y}_{U}, \overline{\mathcal{L}}_{U}^{k,l}),$$
where $U' = U\cap U_0(v)$, the first map is the natural pull-back,
the second map is induced by the maps 
$\rho_g^*\overline{\mathcal{L}}_{U}^{k,l} \to 
\overline{\mathcal{L}}_{U'}^{k,l} \simeq
\rho_1^*\overline{\mathcal{L}}_{U}^{k,l}$
of (\ref{eqn:pullback}), and the last map is the trace.  We remark also that if $(k,l)$
is paritious, then the above definitions with $E$ replaced by $\mathcal{O}$ gives
Hecke operators compatible with the usual ones denoted $T_v$ and $S_v$ on
the corresponding spaces of automorphic forms.

\subsection{Adelic action on components}
We will describe below the action of the group
$\{\,g\in \mathrm{GL}_2(\mathbb{A}_F^\infty)\,|\, g_p \in \mathrm{GL}_2(O_{F,p})\,\}$ on
the spaces $M_{0,l}(E)$, but first we consider the right action via $\rho_g$ on geometric
components. More precisely, suppose as usual that
$g^{-1}U_1 g \subset U_2$ and $N_1$, $N_2$ and $\alpha$ are as in the definition
of $\rho_g$; assume moreover that $\mu_{N_1}(\overline{\mathbb{Q}}) \subset \mathcal{O}$
and consider the map $Z_{U_1}(\mathcal{O}) \to Z_{U_2}(\mathcal{O})$ induced by
$\rho_g$ (where $Z_{U_i}$ was defined in \S\ref{subsec:components}).  
Maintaining the notation in the construction of $\rho_g$, one finds that
the commutativity of the diagram in the proof of Lemma~\ref{lemma:pol}
implies that of
$$\begin{scriptsize}\begin{array}{ccccccc}
J_1\otimes_{O_F}\widehat{O}_F &
\stackrel{\sim}{\longrightarrow} &
\mathrm{Sym}(A_1/S_1)\otimes_{O_F}\widehat{O}_F &
\twoheadrightarrow&
\mathrm{Hom}(\wedge^2_{O_F}A_1[N_1],\mu_{N_1}) &
\stackrel{\sim}{\longrightarrow} &
\mathrm{Hom}(O_F/N_1,\mu_{N_1}) \\&&&&&&\\
\wr\uparrow&&\wr\uparrow&&\downarrow&&\downarrow\\&&&&&&\\
J_2\otimes_{O_F}\widehat{O}_F &
\stackrel{\sim}{\longrightarrow} &
\mathrm{Sym}(A_1'/S_1)\otimes_{O_F}\widehat{O}_F &
\twoheadrightarrow&
\mathrm{Hom}(\wedge^2_{O_F}A_1'[N_2],\mu_{N_2}) &
\stackrel{\sim}{\longrightarrow} &
\mathrm{Hom}(O_F/N_2,\mu_{N_2}),\end{array}\end{scriptsize}$$
where the horizontal arrows of the top (resp.~bottom) row
are induced (from left to right) by $\lambda_1$ (resp.~$\lambda_1'$),
the  Weil pairing on $A_1$ (resp.~$A_1'$), and $\eta_1$
(resp.~$\eta_1'$), the first vertical arrow by
$\beta\det(\alpha g)^{-1}$, the second by $\det(\alpha g)^{-1}\pi^*$,
the third by the surjections $\wedge^2_{O_F}A_1[N_1] \to \wedge^2_{O_F}A_1'[N_2]$
(arising from the isomorphisms $\det(\alpha g_\ell)^{-1}\wedge^2_{O_{F,\ell}}T_\ell(\pi)$
for $\ell|N_2$) and
$\cdot^{N_1/N_2}: \mu_{N_1} \to \mu_{N_2}$, and the last by the natural
projection and $\cdot^{N_1/N_2}$.
It follows that if $\zeta \in Z_{J_1,N_1}(\mathcal{O})$ (i.e., $\zeta:J_1/N_1 
\stackrel{\sim}{\to} \mathfrak{d}^{-1}\otimes \mu_{N_1}(\mathcal{O})$) then
$\tilde{\rho}_g(\zeta)$ is the isomorphism
$J_2/N_2 \simeq \mathfrak{d}^{-1}\otimes \mu_{N_2}(\mathcal{O})$ induced
by $x\mapsto \zeta(\beta\alpha^{-2}\det(g)^{-1}x)^{N_1/N_2}$.
It follows in turn that the map $Z_{U_1}(\mathcal{O}) \to Z_{U_2}(\mathcal{O})$ induced by
$\rho_g$ corresponds to multiplication by $\det(g)^{-1}$
under the bijections of (\ref{eqn:components}),
with $\zeta_{N_2}$ chosen to be $ \zeta_{N_1}^{N_1/N_2}$, 

\subsection{Adelic action on forms of weight $(0,l)$}

Recall that  the map $[U_1 g U_2]$ arises by descent (and reduction mod $\pi$) from maps
\begin{equation}\label{eqn:zaction}
H^0(S_2, \mathcal{L}_{J_2,N_2}^{0,l})
   \to H^0(S_1, \tilde{\rho}_g^*\mathcal{L}_{J_2,N_2}^{0,l})
      \to H^0(S_1, \mathcal{L}_{J_1,N_1}^{0,l})\end{equation}
 where $S_i = Y_{J_i,N_i}$.  Moreover we have isomorphisms
  $\mathcal{L}_{J_i,N_2}^{0,l} \cong \mathcal{O}_{S_i}$ obtained
by tensoring powers of the components of the composite 
$$\bigoplus_{\tau\in\Sigma}\delta_\tau = 
\wedge^2_{O_F\otimes \mathcal{O}_{S_i}}\mathcal{H}^1_{\mathrm{DR}}(A_i/S_i)
\cong J_i\mathfrak{d}^{-1}\otimes \mathcal{O}_{S_i}
 \cong \bigoplus_{\tau\in \Sigma} \mathcal{O}_{S_i},$$
 where the first isomorphism is the canonical one following (\ref{eqn:duality}),
 and the second arises from  the isomorphisms 
$J_i\mathfrak{d}^{-1} \otimes \mathcal{O} \cong 
O_F \otimes \mathcal{O} \cong \oplus_{\tau\in \Sigma} \mathcal{O}$
induced by the inclusions $J_i\mathfrak{d}^{-1} \subset F$ (a choice
permitted by our assumption that $J_1$, $J_2$ and $\mathfrak{d}$
are prime to $p$).
Since the global sections of $\mathcal{O}_{S_i}$ are constant on components,
we may realise (\ref{eqn:zaction}) as a map
$$\{\, Z_{J_2,N_2}(\mathcal{O}) \to \mathcal{O}\, \}
  \longrightarrow \{\, Z_{J_1,N_1}(\mathcal{O}) \to \mathcal{O}\, \}.$$
Under the canonical isomorphisms
$\wedge^2_{O_F\otimes \mathcal{O}_{S_i}}\mathcal{H}^1_{\mathrm{DR}}(A_i/S_i)
\cong J_i\mathfrak{d}^{-1}\otimes \mathcal{O}_{S_i}$, we find that the map
$$\tilde{\rho}_g^*\left(\wedge^2_{O_F\otimes \mathcal{O}_{S_2}}\mathcal{H}^1_{\mathrm{DR}}(A_2/S_2)\right)
\longrightarrow \wedge^2_{O_F\otimes \mathcal{O}_{S_1}}\mathcal{H}^1_{\mathrm{DR}}(A_1/S_1)$$
in the definition of $[U_1 g U_2]$
corresponds to the map $J_2\mathfrak{d}^{-1} \otimes \mathcal{O}_{S_1} \to 
J_1\mathfrak{d}^{-1} \otimes \mathcal{O}_{S_1}$
induced by multiplication by $\beta\alpha^{-2} \in (O_F\otimes\mathcal{O})^\times$.
We therefore realise (\ref{eqn:action}) as the map sending $s: Z_{J_2,N_2}(\mathcal{O}) \to \mathcal{O}$ to the map $Z_{J_1,N_1}(\mathcal{O}) \to \mathcal{O}$ sending
$\zeta$ to $||\det(g)||(\beta\alpha^{-2})^l s(\tilde{\rho}_g(\zeta))$.   Note in particular that
if $\det(g) = 1$ and $U_1 \subset U_2$, then we may choose $\beta = \alpha^2$
and conclude that $[U_1 g U_2]$ coincides with the natural inclusion
$M_{0,l}(U_2;E) \to M_{0,l}(U_1;E)$ defined by $[U_1 1 U_2]$. It follows that
the action of $\{\,g\in \mathrm{GL}_2(\mathbb{A}_F^\infty)\,|\, g_p \in \mathrm{GL}_2(O_{F,p})\,\}$
on $M_{0,l}(E)$ factors  via $\det$ through that of
$\{\,a \in (\mathbb{A}_F^\infty)^\times \,|\, a_p \in O_{F,p}^\times\,\}$, so
we get an action of 
$\{\,g\in \mathrm{GL}_2(\mathbb{A}_F^\infty)\,|\, g_p \in \mathrm{GL}_2(O_{F,p})\,\}$ on
$M_{0,l}(U;E)$ factoring through
$$\{\,a \in (\mathbb{A}_F^\infty)^\times \,|\, a_p \in O_{F,p}^\times\,\} / \det(U).$$

We now determine the corresponding representation of the latter group
on $M_{0,l}(U;E)$.  Note that we have an exact sequence
$$\begin{array}{rcl} 1& \longrightarrow &O_{F,+}^\times \cap \det(U) 
        \longrightarrow  F_+^\times \cap O_{F,p}^\times  \\      &&\\     
        &\longrightarrow &\{\,a \in (\mathbb{A}_F^\infty)^\times \,|\, a_p \in O_{F,p}^\times\,\} / \det(U)
        \longrightarrow (\mathbb{A}_F^\infty)^\times / F_+^\times \det(U)
        \longrightarrow 1,\end{array}$$
where the maps are all induced by the canonical inclusions.  Note that the last quotient
is finite.  If $\overline{\nu}^l=1$ for all $\nu \in \det(U) \cap O_{F,+}^\times$, then 
$\mu \mapsto \overline{\mu}^l$ defines an $E^\times$-valued character of 
$(F_+^\times \cap O_{F,p}^\times)/(O_{F,+}^\times \cap \det(U))$, hence of a finite
index subgroup of  $\{\,a \in (\mathbb{A}_F^\infty)^\times \,|\, a_p \in O_{F,p}^\times\,\} / \det(U)$.

\begin{lemma} \label{lem:components}
If $\overline{\nu}^l = 1$ for all $\nu \in \det(U) \cap O_{F,+}^\times$, then
$M_{0,l}(U;E)$ is isomorphic, as a representation of
$\{\,a \in (\mathbb{A}_F^\infty)^\times \,|\, a_p \in O_{F,p}^\times\,\} / \det(U)$,
to the induction of the character 
$$\psi_l: (F_+^\times \cap O_{F,p}^\times)/(O_{F,+}^\times \cap \det(U)) \to E^\times$$
defined by $\psi_l(\mu) =  \mathrm{Nm}_{F/\mathbb{Q}}(\mu)^{-1}\overline{\mu}^l$;
otherwise $M_{0,l}(U;E) = 0$.
\end{lemma}
\begpf Note that the conclusion of the lemma is equivalent to the assertion that
$M_{0,l}(U;E)$ is isomorphic to
$$I_U = \{\, f: G \to E \,|\, \mbox{$f(\mu x w) = 
\psi_l(\mu) f(x)$ for all $\mu \in G \cap F_+^\times$,
$x\in G$, $w \in \det U$}\,\}$$
as a representation of $G = \{\,a \in (\mathbb{A}_F^\infty)^\times \,|\, a_p \in O_{F,p}^\times\,\}$.
We may therefore replace $L$ by a finite extension and $U$ by an open subgroup $U_2$
for which the hypotheses of  Proposition~\ref{prop:twist} are satisfied.

Next observe that if $\det(g) = \mu \in F_+^\times \cap O_{F,p}^\times$ and
$g^{-1}U_1g \subset U_2$, then we may take
$\beta = \mu \alpha^2$ in the definition of $[U_1 g U_2]$, so that $\tilde{\rho}_g$ induces the
natural projection $Z_{J,N_1}(\mathcal{O}) \to Z_{J,N_2}(\mathcal{O})$ for each $J$,
and the map in (\ref{eqn:zaction}) is the composite of the natural inclusion with
multiplication by $\mathrm{Nm}_{F/\mathbb{Q}}(\mu)^{-1}\overline{\mu}^l$.
Therefore $F_+^\times \cap O_{F,p}^\times$ acts
on $M_{0,l}(U_2;E)$ via the character $\psi_l$.

Let $e$ be a non-zero element of $M_{0,l}(U_2;E)$ supported on a single component
of $ Z_{U_2}(\mathcal{O})$.  Since $F_+^\times \cap O_{F,p}^\times$ acts
via $\psi_l$ on $e$, there is a $G$-equivariant homomorphism $I_{U_2} \to M_{0,l}(U_2;E)$
whose image contains $e$.  Since $G$ acts transitively on $Z_{U_2}(\mathcal{O})$,
the $G$-orbit of $e$ spans $M_{0,l}(U_2;E)$, so the homomorphism is surjective.
Since $I_{U_2}$ and $M_{0,l}(U_2;E)$ both have dimension equal to the cardinality
of $(\mathbb{A}_F^\infty)^\times / F_+^\times \det(U_2)$, it follows that the map is in
fact an isomorphism. \epf

\subsection{Twisting by characters}
\label{subsec:exi}
It follows from Lemma~\ref{lem:components} that for any character 
$\xi: \{\,a \in (\mathbb{A}_F^\infty)^\times \,|\, a_p \in O_{F,p}^\times\,\} / \det(U) \to E^\times$ such that $\xi(\alpha) = \overline{\alpha}^l$ for all $\alpha \in 
F_+^\times \cap O_{F,p}^\times$, the eigenspace consisting of those $e \in M_{0,l}(U;E)$ satisfying
$$ \mbox{$ge = ||\det(g)||\xi(\det(g))e$ for all 
$g \in \mathrm{GL}_2(\mathbb{A}_F^\infty)$ such that  $g_p \in \mathrm{GL}_2(O_{F,p})$}$$
is one-dimensional. We let $e_\xi$ be a basis element.

\begin{lemma}  \label{lem:twist}
If $U$, $l$ and $\xi$ are as above, then for any $k,m \in \mathbb{Z}^\Sigma$, the
map $f \mapsto e_\xi \otimes f$ defines an isomorphism $M_{k,m}(U;E) \to M_{k,l+m}(U;E)$
such that 
$$[UgU] (e_\xi \otimes f) = \xi(\det(g)) e_\xi \otimes [UgU] f$$
 for all $f \in M_{k,m}(U;E)$,
$g \in \mathrm{GL}_2(\mathbb{A}_F^\infty)$ such that  $g_p \in \mathrm{GL}_2(O_{F,p})$;
in particular $T_v(e_\xi \otimes f) = \xi(\varpi_v) e_\xi \otimes T_v f$ and 
$S_v(e_\xi \otimes f) = \xi(\varpi_v)^2 e_\xi \otimes S_v f$ for all $v$ such that
$v\nmid p$ and $\mathrm{GL}_2(O_{F,v}) \subset U$.
\end{lemma}
\begpf We first prove that the map is an isomorphism.  The existence of $\xi$
implies that $\overline{\nu}^l = 1$ for all $\nu \in \det(U) \cap O_{F,+}^\times$, so replacing
$L$ by a finite extension, we may assume that the hypotheses of 
Proposition~\ref{prop:twist} are satisfied and hence view $e_\xi$ as a function
$Z_U(\mathcal{O}) \to E$.  Since $e_\xi$ is non-zero and the action of the group
$\{\,a \in (\mathbb{A}_F^\infty)^\times \,|\, a_p \in O_{F,p}^\times\,\}$ on $Z_U(\mathcal{O})$
is transitive, it follows that $e_\xi$ is everywhere non-zero.
We therefore have a section $e_\xi^{-1} \in M_{0,-l}(U;E)$ such that $f \mapsto
f \otimes e_\xi^{-1}$ defines the inverse of our map.

We now establish the compatibility with the Hecke action.
The definition of $[UgU]$ gives 
$$[UgU] (e_\xi \otimes f) = \sum_i g_i (e_\xi \otimes f) = 
\sum_i || \det(g) ||^{-1} g_i e_\xi \otimes g_i f,$$
where $UgU = \displaystyle\coprod_i  g_i U$.  
Noting that $g_i e_\xi = g e_\xi
 =  ||\det(g)||\xi(\det(g))e_\xi$ since $\det(g_i) \in \det(g) \det U$,
 it follows that 
 $$[UgU] (e_\xi \otimes f) =  \xi(\det(g)) e_\xi \otimes \sum_i g_i f = 
\xi(\det(g)) e_\xi \otimes [UgU] f$$
as required.
\epf

\section{Partial Hasse invariants}

We next adapt the definition of partial Hasse invariants from \cite{AG} to our setting.

\label{sec:Hasse}

\subsection{Definition of partial Hasse invariants}
We write $\mathrm{Ver}_A$ for the Verschiebung isogeny of an abelian scheme $A$ over a base
$S$ of characteristic $p$, i.e., the morphism $A^{(p)} \to A$ defined as the dual of the relative
Frobenius morphism $A^\vee \to (A^\vee)^{(p)} =(A^{(p)})^\vee$, where $A^{(p)}$ denotes the
pull-back $A \times_S S$ with respect to the absolute Frobenius morphism $\mathrm{Fr}_S: S \to S$.
Taking $A$ to be the universal HBAV over $S = \overline{Y}_{J,N}$, the pull-back $\mathrm{Ver}_A^*$
defines an $O_F \otimes \mathcal{O}_S$-linear morphism
$$s_* \Omega^1_{A/S}  \to s_* \Omega^1_{A^{(p)}/S} =  \mathrm{Fr}_S^*s_* \Omega^1_{A/S},$$
where $s:A\to S$ denotes the structure morphism.
 Writing $s_* \Omega^1_{A/S} = \oplus_{\tau} \overline{\omega}_\tau$, we see that the
$\tau$-component of $\mathrm{Fr}_S^* s_*\Omega^1_{A/S}$ is canonically isomorphic to 
$\overline{\omega}_{\mathrm{Fr}^{-1}\circ\tau}^{\otimes p}$, where $\mathrm{Fr}$ denotes the absolute
Frobenius on $\overline{\mathbb{F}}_p$.  The $\tau$-component of $\mathrm{Ver}_A^*$
is therefore a section of $\overline{\mathcal{L}}_{J,N}^{k,0} = \overline{\omega}_{\mathrm{Fr}^{-1}\circ\tau}^{\otimes p}
\overline{\omega}_{\tau}^{\otimes (-1)}$, where:
\begin{itemize}
\item if $\mathrm{Fr}\circ\tau = \tau$, then $k_\tau= p-1$ and $k_{\tau'}= 0$ if $\tau' \neq \tau$;
\item if $\mathrm{Fr}\circ\tau \neq \tau$, then $k_\tau = -1$, $k_{\mathrm{Fr}^{-1}\circ\tau}= p$, and
$k_{\tau'} = 0$ if $\tau' \not\in \{\mathrm{Fr}^{-1}\circ\tau,\tau\}$.
\end{itemize}

For each $\tau$, we denote this weight by $k_{\mathrm{Ha}_\tau}$, and let $\mathrm{Ha}_{J,N,\tau}$ be
the element of $H^0(\overline{Y}_{J,N},\overline{\mathcal{L}}_{J,N}^{k_{\mathrm{Ha}_\tau},0})$ just constructed.
Then $\mathrm{Ha}_{J,N,\tau}$ has non-zero restriction to each component of $\overline{Y}_{J,N}$;
moreover if we let $Z_{\tau}$ denote the associated divisor of zeros, then $Z_{\tau}$ is non-trivial on
each component and $\sum_\tau Z_\tau$ is reduced.
(This follows from the corresponding result proved in \cite[\S8]{AG} for the partial Hasse invariants on the 
variety they denote $\mathfrak{M}(\mathbb{F}_p,\mu_N)$: Choosing $\zeta_N \in \mu_N(E)$ for
sufficiently large $E$ yields an \'etale cover $\overline{Y}_{J,N} \to \mathfrak{M}(E,\mu_N)$
which identifies $\mathfrak{M}(E,\mu_N)$ with the quotient of $\overline{Y}_{J,N}$ by the image of
$U_1(N)$ in $G_{U_1(N),N}$ and our $\mathrm{Ha}_\tau$ with the pull-back of their partial
Hasse invariant $h_{\mathfrak{P},i}$ for the pair $(\mathfrak{P},i)$ corresponding to $\tau$.)

Note that $\mu^{k_{\mathrm{Ha}_\tau}} \equiv 1 \bmod \pi$ for all $\mu \in O_F^\times$, so the line
bundle $\overline{\mathcal{L}}_U^{k_{\mathrm{Ha}_\tau},0}$ is defined for all $U$ under consideration.
By the compatibility of the Verschiebung with base-change and isomorphisms, we see that the
sections $\mathrm{Ha}_{J,N,\tau}$ on $\coprod \overline{Y}_{J,N}$ descend to define a mod $p$
Hilbert modular form of weight $(k_{\mathrm{Ha}_\tau},0)$ and level $U$, which we denote by $\mathrm{Ha}_{U,\tau}$.
Moreover, from the compatibility of Verschiebung with isogenies, in particular with $\pi_\alpha$ as defined in \S\ref{subsec:adelicaction}, we see that $[U_1 g U_2]\mathrm{Ha}_{U_2,\tau} = ||\det g|| \mathrm{Ha}_{U_1,\tau}$
for any $g \in \mathrm{GL}_2(\mathbb{A}_F^\infty)$ such that $g_p \in \mathrm{GL}_2(O_{F,p})$ and $g^{-1}U_1g \subset U_2$.
In particular, the element
$$\mathrm{Ha}_{U,\tau} \in M_{k_{\mathrm{Ha}_\tau},0}(E) :=   \varinjlim M_{k_{\mathrm{Ha}_\tau},0}(U;E)$$
is independent of the choice of $U$, so we henceforth omit the subscript $U$ and write
simply $\mathrm{Ha}_{\tau}$ for this mod $p$ Hilbert modular form, which we call {\em the partial Hasse
invariant (associated to $\tau$}).

We record the following immediate consequence of the assertions above:
\begin{proposition} \label{prop:Ha1} The partial Hasse invariant $\mathrm{Ha}_\tau$ satisfies
$g \mathrm{Ha}_\tau = ||\det(g)|| \mathrm{Ha}_\tau$ for all $g \in \mathrm{GL}_2(\mathbb{A}_F^\infty)$ such that $g_p \in \mathrm{GL}_2(O_{F,p})$.
For any weight $(k,l)$, multiplication by $\mathrm{Ha}_\tau$ defines an injective map:
$$M_{k,l}(E) \to M_{k+k_{\mathrm{Ha}_\tau},l}(E)$$
commuting with the action of $g$ for all such $g$.  In particular, for any open compact subgroup $U$ of $\mathrm{GL}_2(\mathbb{A}_F^\infty)$ containing
$\mathrm{GL}_2(O_{F,p})$, multiplication by $\mathrm{Ha}_\tau$ defines an injective map
$$M_{k,l}(U;E) \to M_{k+k_{\mathrm{Ha}_\tau},l}(U;E)$$
commuting with the operators $T_v$ and $S_v$ for all $v \nmid p$ such that $\mathrm{GL}_2(O_{F,v}) \subset U$.
\end{proposition}

\subsection{Minimal weights}\label{subsec:minimal}
We now recall the definition of the minimal weight of a mod $p$ Hilbert modular form, again adapting notions from \cite{AG} to our setting  (see also \cite{DK}).  This is an analogue of the weight filtration for mod $p$ modular forms in the classical setting 
$F=\mathbb{Q}$.  For $F = \mathbb{Q}$, the vanishing of the spaces of mod $p$ modular forms of negative
weight forces the weight filtration to be non-negative, but in the Hilbert case, the partial negativity of the weights
of partial Hasse invariants already shows the situation is more subtle.  We let
$$\Xi_{\mathrm{AG}} = \left\{\left. \sum_{\tau\in\Sigma} n_\tau k_{\mathrm{Ha}_\tau} \,\right| \, \mbox{$n_\tau \in \mathbb{Z}_{\ge 0}$ for all $\tau \in \Sigma$}\,\right\}$$
be the set of non-negative integer linear combinations of the weights of the partial Hasse invariants.  Note that the weights $k_{\mathrm{Ha}_\tau}$ are linearly
independent, so each $k \in \Xi_{\mathrm{AG}}$ is of the form $\sum_{\tau\in\Sigma} n_\tau k_{\mathrm{Ha}_\tau}$ for a unique $n \in \mathbb{Z}^\Sigma_{\ge 0}$.
We define a partial ordering $\le_{\mathrm{Ha}}$ on $\mathbb{Z}^\Sigma$ by stipulating that $k' \le_{\mathrm{Ha}} k$ if and only if $k - k' \in \Xi_{\mathrm{AG}}$.

For any non-zero $f \in M_{k,l}(U;E)$, consider the set $W(f)$ defined as
$$
\left\{\left.\,  k' = k - \sum_\tau n_\tau k_{\mathrm{Ha}_\tau} \, \right|\, \mbox{$n \in \mathbb{Z}_{\ge 0}^\Sigma$, $f = f'\displaystyle\prod_\tau \mathrm{Ha}_\tau^{n_\tau}$ for some $f' \in M_{k',l}(U;E)$}\,\right \}.$$
Since the divisor $\sum_\tau Z_\tau$ is reduced, the set  $W(f)$ contains a unique minimal element under the partial ordering $\le_{\mathrm{Ha}}$ (cf. \cite[8.19, 8.20]{AG}),
which we call the {\em minimal weight} of $f$, and denote $\nu(f)$.  Note that replacing $U$ by an open compact subgroup $U' \subset U$ does not alter $\nu(f)$, since
any $f' \in M_{k',l}(U';E)$ satisfying $f = f' \prod_\tau \mathrm{Ha}_\tau^{n_\tau}$ will be invariant under $U$, hence in $M_{k',l}(U;E)$.  We may therefore define $\nu(f)$ for
$f \in M_{k,l}(E)$ without reference to $U$.  Note also that $\nu(f)$ is not affected by replacing $E$ by an extension $E'$.

We note also that the minimal weight of a form is independent of $l$ in the following sense:  Recall from Lemma~\ref{lem:twist} that we have isomorphisms
$M_{k,l}(E) \to M_{k,l'}(E)$ defined by multiplication by eigenvectors $e_\xi \in M_{0,l'-l}(E)$ associated to suitable characters $\xi$ of $(\mathbb{A}_F^\infty)^\times$.
Since these isomorphisms commute with multiplication by the partial Hasse invariants, it follows that $\nu(e_\xi \otimes f) = \nu(f)$ for all $f \in M_{k,l}(E)$.

Finally we define the {\em minimal cone} in $\mathbb{Z}^\Sigma$ to be
$$\Xi_{\mathrm{min}} = \{\, k \in \mathbb{Z}^\Sigma \, |\, \mbox{$pk_\tau \ge k_{\mathrm{Fr}^{-1}\circ\tau}$ for all $\tau \in \Sigma$}\,\}.$$
(Note that $\Xi_{\mathrm{min}} \subset \mathbb{Z}_{\ge 0}^\Sigma$.)
A recent result of the first author and Kassaei~\cite{DK} shows that in fact $\nu(f) \in \Xi_{\mathrm{min}}$ for all non-zero mod $p$ Hilbert
modular forms $f$.

\section{Associated Galois representations}
\label{sec:galois}

The aim of this section is to prove the existence of Galois representations associated to Hecke eigenforms
of arbitrary weight.  We first state the theorem and review some ingredients needed for the proof.

\subsection{Statement of the theorem}

\begin{theorem} \label{thm:galois}  Suppose that $U$ is an open compact subgroup of $\mathrm{GL}_2(\widehat{O}_F)$
containing $\mathrm{GL}_2(O_{F,p})$, and $Q$ is a finite set of primes containing all $v|p$ and
all $v$ such that $\mathrm{GL}_2(O_{F,v}) \not\subset U$.  Suppose that $k,l\in \mathbb{Z}^{\Sigma}$ and
that $f \in M_{k,l}(U;E)$ is an eigenform for
$T_v$ and $S_v$ (defined in (\ref{eqn:TvSv})) for all $v \not\in Q$.   Then there is a Galois representation
$$\rho_f : G_F \to \mathrm{GL}_2(E)$$
such that if $v\not\in Q$, then $\rho_f$ is unramified at $v$ and the characteristic
polynomial of $\rho_f(\mathrm{Frob}_v)$ is 
$$X^2 - a_v  X + d_v \mathrm{Nm}_{F/\mathbb{Q}}(v),$$
where $T_v f = a_v f$ and $S_v f = d_v f$.
\end{theorem}

This has been proved for paritious weights $(k,l)$, independently by Emerton--Reduzzi--Xiao~\cite{ERX} and
Goldring--Koskivirta~\cite{GoK}; in fact their methods yield the result under a weaker parity condition.
The contribution here is to remove the parity hypothesis altogether,  and the
new ingredient is to use congruences to forms of level divisible by $p$.   For this we will need to 
work with the integral models for Hilbert modular varieties with level structure $U_1(p)$ at $p$
studied by Pappas in~\cite{P}.

\subsection{Hilbert modular varieties of level $U'=U\cap U_1(p)$}

Suppose that $J$ is a fractional ideal of $F$ and $N \ge 3$ is an integer, with $J$ and $N$ both prime to $p$. 
We let  $\mathcal{M}^{0}_{J,N}$ denote the functor which associates to an $\mathcal{O}$-scheme $S$
the set of isomorphism classes of pairs $(\underline{A},H)$, where 
\begin{itemize}
\item $\underline{A} = (A,i,\lambda,\eta)$ is a $J$-polarised HBAV with level $N$-structure over $S$, and
\item $H$ is a free rank one $(O_F/p)$-submodule scheme of $A[p]$ over $S$ such that the quotient
isogeny $A \to A' = A/H$ induces an isomorphism $\mathrm{Sym}(A'/S) \to p\mathrm{Sym}(A/S)$.

\end{itemize}
Then $\mathcal{M}^0_{J,N}$ is represented by an $\mathcal{O}$-scheme which we denote $Y_{J,N}^0$,
the forgetful morphism $Y_{J,N}^0 \to Y_{J,N}$ is projective and $Y_{J,N}^0$ is a flat local complete intersection over $\mathcal{O}$
of relative dimension $[F:\mathbb{Q}]$ (\cite[Thm.~2.2.2]{P}).
We let $\mathcal{M}^{1}_{J,N}$ denote the functor which associates to an $\mathcal{O}$-scheme $S$
the set of isomorphism classes of triples $(\underline{A},H,P)$ where $\underline{A}$ and $H$ are as above and
\begin{itemize}
\item $P \in H(S)$ is an $(O_F/p)$-generator of $H$ in the sense of Drinfeld--Katz--Mazur~\cite[1.10]{KM}.
\end{itemize}
Then $\mathcal{M}^1_{J,N}$ is represented by $\mathcal{O}$-scheme which we denote $Y_{J,N}^1$, and
the forgetful morphism $Y_{J,N}^1 \to Y_{J,N}^0$ is finite flat, so $Y_{J,N}^1$ is flat and
Cohen--Macaulay over $\mathcal{O}$ (\cite[Thm.~2.3.3]{P}).

Suppose $U$ is an open compact subgroup of $\mathrm{GL}_2(\widehat{O}_F)$ containing $\mathrm{GL}_2(O_{F,p})$,
and let $U' = U \cap U_1(p)$.  We suppose that $U$ is sufficiently small, and in particular that $U$ is $p$-neat (see Definition~\ref{def:pneat}).
The action of the group $G_{U,N}$ on $Y_{J,N}$ then lifts to one on $Y^1_{J,N}$, corresponding to
the action on $\mathcal{M}^1_{J,N}$ defined by $(\nu,u)\cdot(\underline{A},H,P) = ((\nu,u)\cdot\underline{A},H,P)$.
It follows from the corresponding assertions for $Y_{J,N}$ that $G_{U,N}$
acts freely on $\coprod_{t\in T} Y_{J_t,N}^1$, the quotient is representable by a scheme $Y_{U'}$, and
the quotient map is \'etale and Galois with group $G_{U,N}$.   Since the $Y_{J,N}^1$ are flat and Cohen--Macaulay
over $\mathcal{O}$, so is $Y_{U'}$, and let $\pi_U:Y_{U'} \to Y_U$ denote the natural projection
(writing just $\pi$ when $U$ is clear from the context).   We let $\mathcal{K}_{U'}$ denote the dualising
sheaf on $Y_{U'}$ over $\mathcal{O}$ (see \cite[\S3.5]{Co2}), and similarly let $\mathcal{K}_{U}$ denote
the dualising sheaf on $Y_U$ over $\mathcal{O}$.  Since $Y_U$ is smooth over $\mathcal{O}$, its
dualising sheaf $\mathcal{K}_U$ is canonically identified with $\Omega^{[F:\mathbb{Q}]}_{Y_U/\mathcal{O}}
= \wedge^{[F:\mathbb{Q}]}_{\mathcal{O}_{Y_U}} \Omega^1_{Y_U/\mathcal{O}}$.

Suppose now that $g$, $U_1$ and $U_2$ are as in \S\ref{section:Hecke}, so in particular
$g^{-1}U_1g \subset U_2$, and assume further that $g_p \in U_1(p)$.   We then obtain exactly
as before finite \'etale $\rho_g':  Y_{U'_1} \to Y_{U'_2}$, by descent from morphisms 
$\tilde{\rho}_g^1: \coprod Y_{J_t,N_1}^1 \to \coprod Y_{J_t,N_2}^1$, and 
compatible with $\rho_g: Y_{U_1} \to Y_{U_2}$ via the projections $\pi_{U_i}: Y_{U_i'}\rightarrow Y_{U_i}$.
Since $\rho_g'$ is \'etale, we have a canonical isomorphism $(\rho_{g}')^* \mathcal{K}_{U_2'}
\stackrel{\sim}{\to} \mathcal{K}_{U_1'}$.

\subsection{Hilbert modular forms of level $U'=U\cap U_1(p)$}

For $(m, n) \in \mathbb{Z}^2$ (viewed also as an element of $(\mathbb{Z}^\Sigma)^2$) and $p$-neat $U$, we let 
$\mathcal{L}^{m,n}_{U'} = \pi_U^*\mathcal{L}_U^{m,n}$,
and we similarly define $\mathcal{L}^{m,n}_{U',R}$ for $\mathcal{O}$-algebras $R$,
writing also $\overline{\mathcal{L}}^{m,n}_{U'}$ in the case $R=E$.  For $k,l \in \mathbb{Z}$,
we define the space of Hilbert modular forms over $R$ of weight $k$ and level $U'$ to be
$$M_{k,l}(U';R) :=  H^0(Y_{U',R},  \mathcal{K}_{U',R} \otimes_{\mathcal{O}_{Y_{U',R}}} \mathcal{L}^{k-2,l+1}_{U',R}).$$
Note that we could have made this definition for more general weights $(k,l)$, but we will in fact
only need the case of parallel weight.  Recall also from \cite{Co2} that formation of the dualising sheaf is compatible
with base change, so  $\mathcal{K}_{U',R}$ can be identified with the dualising sheaf of $Y_{U',R}$ over $R$.

For $g$, $U_1$, $U_2$ as above, we define an $R$-linear map
$$[U'_1 g U'_2] : M_{k,l}(U_2'; R) \to M_{k,l}(U_1';R)$$
as $||\det g||$ times the composition of the pull-back from $Y_{U_2'}$ to $Y_{U_1'}$ with the map
on sections induced by the tensor product of the canonical isomorphism $(\rho_{g}')^* \mathcal{K}_{U_2'}
\stackrel{\sim}{\to} \mathcal{K}_{U_1'}$ with the map 
$$(\rho_g')^*\mathcal{L}^{k-2,l+1}_{U_2'} = \pi_{U_1}^*\rho_g^* \mathcal{L}^{k-2,l+1}_{U_2}  
  \to \pi_{U_1}^* \mathcal{L}^{k-2,l+1}_{U_1}  = \mathcal{L}^{k-2,l+1}_{U_1'}$$
given by $\pi_{U_1}^*$ of (\ref{eqn:pullback}).  
We again have the compatibility $[U'_1 g_1 U'_2] \circ [U'_2 g_2 U'_3] = [U'_1 g_1g_2 U'_3]$,
giving rise to an $R$-linear action of the group $\{\, g \in \mathrm{GL}_2(\mathbb{A}_F^\infty)\,|\, g_p \in U_1(p)\,\}$
on $M'_{k,l}(R) := \varinjlim M_{k,l}(U';R)$.  As before we may identify $M_{k,l}(U';R)$ with $(M'_{k,l}(R))^{U'}$,
and define commuting $R$-linear Hecke operators $T_v$ and $S_v$ on $M_{k,l}(U';R)$ for all 
$v$ such that $\mathrm{GL}_2(O_{F,v}) \subset U'$.

Let $S = Y_{J,N}$, and $A = A_{J,N}$ the universal HBAV over $S$.
Since $A$ is smooth over $S$ and $S$ is smooth over $\mathcal{O}$, we have an exact sequence
$$0 \to  s^*\Omega^1_{S/\mathcal{O}}  \to \Omega^1_{A/\mathcal{O}}   \to \Omega^1_{A/S} \to 0$$
of locally free sheaves on $A$.   Applying $R^i s_*$, we obtain the
connecting homomorphism:
\begin{equation}\label{KS0}
s_* \Omega^1_{A/S} \longrightarrow  R^1s_*s^*\Omega^1_{S/\mathcal{O}}  \cong   \Omega^1_{S/\mathcal{O}}  \otimes_{\mathcal{O}_S}  R^1s_*\mathcal{O}_A.\end{equation}
Combined with the canonical isomorphisms
$$\mathcal{H}om_{O_F\otimes\mathcal{O}_S}(\wedge^2_{O_F\otimes\mathcal{O}_S}\mathcal{H}^1_{\mathrm{DR}}(A/S),\,
s_*\Omega^1_{A/S}) 
 \stackrel{\sim}{\longrightarrow}
\mathcal{H}om_{\mathcal{O}_S}(R^1s_*\mathcal{O}_A, \,\mathcal{O}_S )$$
induced by the inclusion $O_F \subset \mathfrak{d}^{-1} \stackrel{\sim}{\longrightarrow}\mathrm{Hom}(O_F,\mathbb{Z})$
and the isomorphism (\ref{wedge de Rham}),
we obtain an $ \mathcal{O}_S$-linear homomorphism
\begin{equation}\label{eqn:KS1}
\mathcal{H}om_{O_F\otimes\mathcal{O}_S} \left( \wedge^2_{O_F\otimes\mathcal{O}_S}\mathcal{H}^1_{\mathrm{DR}}(A/S),\,
{\otimes_{O_F\otimes \mathcal{O}_S}^2}(s_*\Omega^1_{A/S} )\right)  \to \Omega^1_{S/\mathcal{O}},
\end{equation}
which is in fact an isomorphism (see \cite[1.0.21]{K}), called the {\em Kodaira--Spencer isomorphism}.   
Taking $\wedge^{[F:\mathbb{Q}]}_{\mathcal{O}_S}$, we obtain an isomorphism:
$$\xi_{J, N}: \mathcal{L}^{2,-1}_{J,N}  = \otimes_\tau ( \omega_\tau^2 \otimes_{\mathcal{O}_S} \delta_\tau^{-1} )
 \longrightarrow  \wedge_{\mathcal{O}_S}^{[F:\mathbb{Q}]} \Omega^1_{S/\mathcal{O}} = \Omega^d_{S/\mathcal{O}}.$$
The functoriality of the morphisms in the construction ensures that the isomorphism is compatible with the action of
$G_{U,N}$, and therefore descends to an 
 isomorphism $$\xi_U: \mathcal{L}^{2,-1}_{U}  \cong \mathcal{K}_U.$$ 
 Moreover for $g$, $U_1$ and $U_2$ such that $g^{-1}U_1g \subset U_2$, $g_p \in \mathrm{GL}_2(O_{F,p})$,
 one finds similarly that the canonical isomorphism $\rho_g^*\mathcal{K}_{U_2} \to \mathcal{K}_{U_1}$
 is compatible with the morphism of~(\ref{eqn:pullback}).  It follows that the
 isomorphisms
 $$M_{k,l}(U;R) \cong  H^0(Y_{U,R},  \mathcal{K}_{U,R} \otimes_{\mathcal{O}_{Y_{U,R}}} \mathcal{L}^{k-2,l+1}_{U,R})$$
 induced by $\xi_U$ are compatible with the operators $[U_1 g U_2]$.
 Moreover the generic fibre of $Y_{J,N}^1$ is smooth
 over $L$, so that if $p$ is invertible in $R$, the same constructions apply to give isomorphisms
 $$H^0(Y_{U',R}, \mathcal{L}_{U',R}^{k,l})  \cong M_{k,l}(U';R)$$
 such that the operators $[U_1' g U_2']$ are compatible by extension of scalars with those on the
 spaces $A_{k,l}(U')$ of automorphic forms of weight $(k,l)$ and level $U'$.
 
 \subsection{Minimal compactifications}
 \label{subsec:XU}
  
We will also make use of minimal compactifications of Hilbert modular varieties, whose properties
we now recall.  The minimal compactification $X_{J,N}$ of $Y_{J,N}$ is constructed by Chai
in \cite{C2} (see also \cite{Di} and \cite{DT}), and we define $X_U$ to be the quotient of $\coprod X_{J,N}$ under
the natural action of $G_{U,N}$.  Then $X_U$ is a flat, projective scheme over $\mathcal{O}$
with $j:Y_U \to X_U$ as an open subscheme whose complement is finite over $\mathcal{O}$,
and the line bundle $\mathcal{L}_U^{1,0}$ extends to an ample line bundle on $X_U$ which we denote
by $\mathscr{L}_U$.  The Koecher Principle in this setting means that the natural map
$\mathcal{O}_{X_U} \to j_*\mathcal{O}_{Y_U}$ is an isomorphism.

\begin{definition} \label{def:cusp}
Assuming as usual that $\mathcal{O}$ is sufficiently large (i.e., containing the $N^\mathrm{th}$ roots
of unity), then each (reduced) connected component $C$ of $X_U - Y_U$ is isomorphic to
$\mathrm{Spec}\,\mathcal{O}$.  We call $C$ a {\em cusp} of $X_U$. \end{definition}

 If $U$ is of the form $U(\mathfrak{n}):=\mathrm{ker}(\mathrm{GL}_2(\widehat{O}_F)\rightarrow \mathrm{GL}_2(O_F/\mathfrak{n}))$
for a sufficiently small, prime-to-$p$ ideal $\mathfrak{n}$ of $O_F$, then the completion of $X_U$ along $C$
is canonically isomorphic to $\mathrm{Spf}\, \widehat{S}_C$, where 
\begin{equation}\label{eqn:RZ}
\widehat{S}_C :=  \mathcal{O}[[q^\alpha]]_{\alpha\in I_+ \cup \{0\}}^{U_{\mathfrak{n},+}^\times}
\end{equation}
for a fractional ideal $I$ depending on $C$, and
$\mu \in U_{\mathfrak{n},+} = \ker(O_{F,+}^\times \to (O_F/\mathfrak{n})^\times$) 
acts via $q^\alpha \mapsto q^{\mu\alpha}$.
(The $\mathcal{O}$-algebra $\widehat{S}_C$ is obtained from the corresponding one in \cite{C2}
by working over $\mathcal{O}$ instead of $\mathbb{Z}[\mu_N,1/N]$ and taking invariants
under the stabiliser in $G_{U,N}$ of a cusp $\tilde{C}$ of $X_{J,N}$ mapping to $C$.
In particular, the class of the ideal $I$ in (\ref{eqn:RZ}) is given by $\mathfrak{abn}^{-1}$
where $\mathfrak{a}$ and $\mathfrak{b}$ are as in \cite{C2}; a more detailed discussion
in the case of arbitrary $U$ is provided below in \S\ref{sec:qexpansions}, where
Proposition~\ref{prop:QCkl} gives (\ref{eqn:RZ}) as a special case.)

The minimal compactification of $Y_{U'}$ is then obtained as follows.  First one constructs a
toroidal compactification $X_{U'}^{\mathrm{tor}}$ of $Y_{U'}$ as the quotient of a toroidal compactification
of $\coprod Y_{J,N}^1$ defined exactly as for $\coprod Y_{J,N}$, but using the functors $\mathcal{M}_{J,N}^1$
and $\Gamma(Np)$-admissible polyhedral cone decompositions (in the terminology of \cite{C2}).
Then $\pi:Y_{U'} \to Y_U$ extends to
a projective morphism $X_{U'}^{\mathrm{tor}} \to X_U$ such that the connected components
of the pre-image of a cusp $C$ correspond to pairs $(\mathfrak{f},P)$ where $pO_F \subset \mathfrak{f} \subset O_F$
and $P$ is an $(O_F/p)$-generator of $O_F/\mathfrak{f}$ (or more canonically, $\mathfrak{b}/\mathfrak{bf}$).
Moreover a similar calculation to the case of level $U$ shows that if $U = U(\mathfrak{n})$, then
the ring of global sections of
the completion of $X_{U'}^{\mathrm{tor}}$ along the component over $C$ corresponding to $(\mathfrak{f},P)$
is isomorphic to the $\widehat{S}_C$-algebra
\begin{equation}\label{eqn:RZ1}
\mathcal{O}'_{\mathfrak{f}}[[q^\alpha]]_{\alpha\in (\mathfrak{f}^{-1}I)_+ \cup \{0\}}^{U_{\mathfrak{n},+}},
\end{equation}
where $\mathrm{Spec}\,\mathcal{O}'_{\mathfrak{f}}$ is the finite flat 
$\mathcal{O}$-scheme representing $(O_F/p)$-generators of
$\mu_p\otimes \mathfrak{f}/pO_F$ (or more canonically, $\mu_p\otimes\mathfrak{fa}^{-1}\mathfrak{d}^{-1}/
p\mathfrak{a}^{-1}\mathfrak{d}^{-1}$).

Now let $X_U^{\mathrm{ord}}$ denote the ordinary locus of $X_U$, so $X_U^{\mathrm{ord}}$ is an open
subscheme of $X_U$ containing the cusps, and let  $Y_U^{\mathrm{ord}} = Y_U \cap X_U^{\mathrm{ord}}$.
Let $X_{U'}^{\mathrm{tord}}$ (resp.~$Y_{U'}^{\mathrm{ord}}$) denote the pre-image of $X_U^{\mathrm{ord}}$
(resp.~$Y_U^{\mathrm{ord}}$) in $X_{U'}^{\mathrm{tor}}$ (resp.~$Y_{U'}$), and define
$$X_{U'}^{\mathrm{ord}} = \mathbf{Spec}\,f_*\mathcal{O}_{X_{U'}^{\mathrm{tord}}}$$
where $f:X_{U'}^{\mathrm{tord}} \to X_U^{\mathrm{ord}}$ is the restriction of $X_{U'}^{\mathrm{tor}} \to X_U$.
Since $f$ is proper, $X_{U'}^{\mathrm{ord}}$ is finite over $X_U^{\mathrm{ord}}$, and since 
$Y_{U'}^{\mathrm{ord}} \to Y_U^{\mathrm{ord}}$ is finite, we can identify $Y_{U'}^{\mathrm{ord}}$
with an open subscheme of $X_{U'}^{\mathrm{ord}}$.  We then define the minimal compactification
$j':Y_{U'} \to X_{U'}$ by gluing $Y_{U'}$ and $X_{U'}^{\mathrm{ord}}$ along $Y_{U'}^{\mathrm{ord}}$.

Thus $X_{U'}$ is flat over $\mathcal{O}$, and the morphism $\pi$ extends to a projective morphism
$\tilde{\pi}: X_{U'} \to X_U$, so in particular $X_{U'}$ is projective over $\mathcal{O}$.
Furthermore the restriction $\tilde{\pi}^{\mathrm{ord}}:X_{U'}^{\mathrm{ord}} \to X_U^{\mathrm{ord}}$
is finite, and $\pi^{\mathrm{ord}}:Y_{U'}^{\mathrm{ord}} \to Y_U^{\mathrm{ord}}$ is finite flat.
The cusps $C'$ of $X_{U'}$ (i.e., the reduced connected components of $X_{U'} - Y_{U'}$) lying over
a cusp $C$ of $X_U$ correspond to pairs $(\mathfrak{f},P)$ as above, and in the case $U= U(\mathfrak{n})$,
the completion of
$X_{U'}$ along $C'$ is isomorphic to $\mathrm{Spf}\,\widehat{S}_{C'}$ where $\widehat{S}_{C'}$
is the $\widehat{S}_{C}$-algebra defined by (\ref{eqn:RZ1}) above.  Note in particular that
if $\mathfrak{f} = O_F$, then $\widehat{S}_{C'} = \mathcal{O}_{O_F}'\otimes_{\mathcal{O}}\widehat{S}_C$
is flat over $\widehat{S}_C$.  The Koecher Principle carries over to show that
$j'_*\mathcal{O}_{Y_{U'}} = \mathcal{O}_{X_{U'}}$, and we let
$\mathscr{L}_{U'}$ denote the pull-back $\tilde{\pi}^*\mathscr{L}_U$ of the ample line bundle
$\mathscr{L}_U$.

\subsection{Proof of Theorem~\ref{thm:galois}}\label{subsec:proof}

We begin the proof with some preliminary reductions.  

First we claim we can replace the field $E$
by a finite extension $E'$.  Indeed if $\rho: G_F \to \mathrm{GL}_2(E')$ satisfies the
conclusion of the theorem with $E$ replaced by $E'$, then in fact $\rho$ is defined
over $E$. For $p>2$, this follows by an elementary argument using an
element of $g \in G_F$ (a complex conjugation, for example) such that $\rho(g)$
has distinct eigenvalues in $E$.  For $p=2$, one can twist by the character $\xi:G_F \to E^\times$
such that $\xi^2 = \det\rho$ so as to assume $\det\rho = 1$, and then use the classification
of subgroups of $\mathrm{SL}_2(E') \cong \mathrm{PGL}_2(E')$ to arrive at the desired
conclusion.

Next we claim that we can assume $U = U(\mathfrak{n})$ for a sufficiently small ideal $\mathfrak{n}$
prime to $p$.  Indeed by the proof of Chevalley's Theorem
on congruence subgroups of $O_F^\times$,  we can choose ideals $\mathfrak{n}_1$ and
$\mathfrak{n}_2$ relatively prime to each other and to $p$ so that the kernels of
reduction mod $\mathfrak{n}_i$ for $i=1,2$ are contained in that of reduction mod $p$.
We may then apply the theorem with $U$ replaced by $U(\mathfrak{mn}_i)$, where 
$U(\mathfrak{m}) \subset U$ and $\mathfrak{m}$ is divisible only by primes such
that $\mathrm{GL}_2(O_{F,v}) \not\subset U$.  This produces
representations $\rho_i$ satisfying the conclusions with $Q$ augmented by the set
of primes dividing $\mathfrak{n}_i$.  Moreover we can replace the $\rho_i$ by their
semi-simplifications, which are isomorphic to each other by the Brauer--Nesbitt and Cebotarev Density
theorems.  We therefore obtain the desired conclusion for all $v \not\in Q$.
 
Next we show that we can assume $l = -1$,\footnote{In fact any parallel $l$ will do; the choice
of $l=-1$ is made for later convenience.} i.e., that $l_\tau = -1$ for all $\tau \in \Sigma$.
Given any $l$, define $l' \in \mathbb{Z}^\Sigma$ by $l'_\tau = l_\tau+1$.
Recall from the discussion before Lemma~\ref{lem:components} that our hypothesis on $U$
ensures that $\mu \mapsto \overline{\mu}^{l'}$ is a well-defined $E^\times$-valued character
on the finite index subgroup $(F_+^\times \cap O_{F,p}^\times)/(O_{F,+}^\times \cap \det(U))$
of $\{\,a \in (\mathbb{A}_F^\infty)^\times \,|\, a_p \in O_{F,p}^\times\,\} / \det(U)$, for which
we may choose an extension $\xi$ as in Lemma~\ref{lem:twist} (enlarging $E$ if necessary).
The case $l=-1$ of the theorem then furnishes a Galois representation $\rho_{f\otimes e_{\xi}^{-1}}$
unramified at all $v\not\in Q$ with $\mathrm{Frob}_v$ having  characteristic polynomial
$$X^2 - \xi(\varpi_v)^{-1} a_v X + \xi(\varpi_v)^{-2} d_v\mathrm{Nm}_{F/\mathbb{Q}}(v).$$
Let $V = \{\,b\in \det(U)\,|\,b_p \equiv 1 \bmod p\,\}$, and define
$$\xi': \mathbb{A}_F^\times/F^\times F_{\infty,+}^\times V \to E^\times$$
by $\xi'(\alpha za) = \xi(a)\overline{a}_p^{-l'}$ for $\alpha \in F^\times$, $z\in F_{\infty,+}^\times$
and $a \in (\mathbb{A}_F^\infty)^\times$ with $a_p \in O_{F,p}^\times$.  
Letting $\rho_{\xi'}: G_F \to E^\times$ be the character corresponding to $\xi'$ by class field
theory, we have $\rho_{\xi'}(\mathrm{Frob}_v) = \xi(\varpi_v)$ for all $v\not\in Q$, so the
representation $\rho_{\xi'} \otimes \rho_{f\otimes e_{\xi}^{-1}}$ satisfies the conclusion
of the theorem.

Now we reduce to the case where $f$ is of arbitrarily large, ``nearly parallel'' weight.
More precisely, we claim that, given any $M \in \mathbb{Z}$, we can assume that
$k = (k_{\tau})_{\tau\in\Sigma}$ has the form $k = m + 2 - \kappa = (m+2-\kappa_\tau)_\tau$, where
\begin{itemize}
\item $m\in \mathbb{Z}$, $m \ge M$;
\item $0 \le \kappa_\tau \le p-1$ for all $\tau \in \Sigma$;
\item for each $v|p$, $\kappa_\tau < p-1$ for some $\tau \in \Sigma_v$.
\end{itemize}
Here we have identified $\Sigma$ with the set of embeddings $O_F \to \mathcal{O}$
and written 
\begin{equation}  \label{eqn:Sigmav} \Sigma = \coprod_{v|p} \Sigma_v, \mbox{where 
$\Sigma_v = \{\,\tau\in \Sigma\,|\, v = \tau^{-1}(\pi\mathcal{O})\,\}$.}\end{equation}

To prove the claim, suppose $f \in M_{k,-1}(U;E)$, and choose any $m \in \mathbb{Z}$
such that $m \ge M$ and $m \ge k_\tau+p-3$ for all $\tau \in \Sigma$.
For each $v|p$, choose some $\tau_{v,0} \in \Sigma_v$ and let
$\tau_{v,i} = \mathrm{Fr}^i\circ \tau_{v,0}$, so 
$\Sigma_v = \{\,\tau_{v,i}\,|\,i=0,\ldots,f_v-1\,\}$
where $f_v = [O_F/v : \mathbb{F}_p]$.
Now let $r \in \mathbb{Z}$ be such that $0 \le r < p^{f_v} -1$ and
$$r\,  \equiv  \,\sum_{i=0}^{f_v-1} (m+2-k_{\tau_{v,i}})p^i \,\, \bmod (p^{f_v} - 1).$$
We then define $\kappa_\tau$ for $\tau\in \Sigma_v$ by requiring
that $0 \le \kappa_{\tau_{v,i}} \le p-1$ for $i=0,\ldots,f_v-1$ and
$r = \sum  \kappa_{\tau_{v,i}}p^i$.  Note that the resulting $\kappa_\tau$
is independent of the choice of $\tau_{v,0}$ and that $\kappa_\tau < p-1$
for some $\tau \in \Sigma_v$.  Now define $k' = (k'_\tau)_\tau \in \mathbb{Z}^{\Sigma}$ by
setting $k'_\tau = m +2 - \kappa_\tau$.   We then have $k' - k = \sum n_\tau k_{\mathrm{Ha}_\tau}$
where
$$n_\tau = (p^{f_v} - 1)^{-1} \sum_{i=0}^{f_v-1} (k'_{\mathrm{Fr}^i\circ\tau} - k_{\mathrm{Fr}^i\circ\tau}) p^i$$
for $\tau \in \Sigma_v$.  Note that $n_\tau \in \mathbb{Z}_{\ge 0}$ for all $\tau \in \Sigma$, so
$k' - k \in \Xi_{\mathrm{AG}}$.  By Proposition~\ref{prop:Ha1} there is a Hecke-equivariant
injection $M_{k,-1}(U;E) \to M_{k',-1}(U;E)$, so the theorem for forms of weight $(k,-1)$ follows
from the case of weight $(k',-1)$.

The heart of the proof is to construct, for $k = m+2- \kappa$ as above,
a Hecke-equivariant injective homomorphism
$$M_{k,-1}(U;E) \to M_{m+2,-1}(U';E).$$
Letting $\underline{A}$ denote the universal HBAV over $S = \overline{Y}_{J,N}$,
$\mathrm{Frob}_A:A \to A^{(p)}$ the relative Frobenius morphism and $H^\mu = \ker\mathrm{Frob}_A$, the pair $(\underline{A},H^\mu)$
defines a section $\overline{Y}_{J,N} \to \overline{Y}^0_{J,N}$, where as usual we use
$\overline{Y}$ to denote the special fibre of an $\mathcal{O}$-scheme $Y$.
Moreover the section identifies $\overline{Y}_{J,N}$ with a union of irreducible
components of $\overline{Y}^0_{J,N}$, whose pre-image in $\overline{Y}^1_{J,N}$
we denote by $Y^\mu_{J,N}$.  The action of $G_{U,N}$ on $\overline{Y}^1_{J,N}$
restricts to one on $Y^\mu_{J,N}$, and we let $Y^\mu_U$ denote the 
corresponding quotient of $\coprod Y^\mu_{J_t,N}$.  Thus 
$i:Y^\mu_U \to \overline{Y}_{U'}$ is a closed immersion identifying $Y^\mu_U$ with a union of irreducible
components of $\overline{Y}_{U'}$, and $\overline{\pi}\circ i:  Y^\mu_U \to \overline{Y}_U$ is finite flat.
In particular $Y^\mu_U$ is Cohen--Macaulay (over $E$), and we let $\mathcal{K}^\mu_U$
denote its dualising sheaf.   By Grothendieck--Serre duality (\cite[Thm.~3.4.4]{Co2}, and the 
compatibility~\cite[(3.3.14)]{Co2}) applied to the finite morphisms $i$ and $\overline{\pi}\circ i$, we
have canonical isomorphisms:
\begin{equation}\label{eqn:dualising} \begin{array}{rcl}
i_*\mathcal{K}^\mu_{U}& \cong& \mathcal{H}om_{\mathcal{O}_{\overline{Y}_{U'}}}(i_*\mathcal{O}_{Y^\mu_U}, \overline{\mathcal{K}}_{U'})\\
\mbox{and}\quad \overline{\pi}_*i_* \mathcal{K}^\mu_{U} & \cong  &\mathcal{H}om_{\mathcal{O}_{\overline{Y}_{U}}}(\overline{\pi}_*i_*\mathcal{O}_{Y^\mu_U}, \overline{\mathcal{K}}_{U}).
\end{array}
\end{equation}
Since $i$ is a closed immersion, the first of these isomorphisms identifies $i_*\mathcal{K}^\mu_{U}$ with a subsheaf of $\overline{\mathcal{K}}_{U'}$.
To exploit the second isomorphism, we recall that \cite[Prop.~5.1.5]{P} identifies $Y^1_{J,N}$ with a closed subscheme of
the universal submodule scheme $H$ over $Y^0_{J,N}$.  In particular, if $\underline{A}$ is the universal HBAV on $S = \overline{Y}_{J,N}$,
then %
$$H^\mu = \mathbf{Spec}\,\left(\mathrm{Sym}_{\mathcal{O}_S} \left(\oplus_{\tau \in \Sigma} \mathcal{L}_\tau\right)/ 
  \langle \mathcal{L}_\tau^{\otimes p}\,\mbox{for}\,\tau\in \Sigma \rangle\right)$$
as a Raynaud $(O_F/p)$-module scheme (i.e., the morphisms $\Delta_\tau:\mathcal{L}_\tau^{\otimes p} \to \mathcal{L}_{\mathrm{Fr}\circ\tau}$ of \cite[4.4.1]{P}
are zero), so that
$$Y^\mu_{J,N} = \mathbf{Spec}\,\left(\mathrm{Sym}_{\mathcal{O}_S} \left( \oplus_{\tau \in \Sigma}  \mathcal{L}_\tau\right)/ 
\langle \mathcal{L}_\tau^{\otimes p}\,\mbox{for}\,\tau\in \Sigma,\, \otimes_{\tau\in\Sigma_v} \mathcal{L}_{\tau}^{\otimes(p-1)}\,\mbox{for}\, v|p
  \rangle\right),$$
where the $\mathcal{L}_\tau$ are line bundles on $S$.  Moreover the inclusion $H^\mu \to A$ induces
a canonical $\mathcal{O}_S\otimes O_F$-linear isomorphism:
$$s_*\Omega^1_{A/S}  \cong e^*\Omega^1_{A/S}\, \cong \,e^*\Omega^1_{H^\mu/S} \, \cong\, \oplus_{\tau\in \Sigma}  \mathcal{L}_\tau,$$
and hence isomorphisms $\overline{\omega}_\tau \cong \mathcal{L}_\tau$ of line bundles on $S$ for $\tau \in \Sigma$.
These isomorphisms are compatible with the action of $G_{U,N}$, and so give rise to an isomorphism
$$Y^\mu_{U} \cong \mathbf{Spec}\,\left(\mathrm{Sym}_{\mathcal{O}_{\overline{Y}_U}} \left( \oplus_{\tau \in \Sigma}  \overline{\omega}_\tau\right)/ 
\langle \overline{\omega}_\tau^{\otimes p}\,\mbox{for}\,\tau\in \Sigma,\, \otimes_{\tau\in\Sigma_v} \overline{\omega}_{\tau}^{\otimes(p-1)}\,\mbox{for}\, v|p
  \rangle\right),$$
which in turn gives an isomorphism $\overline{\pi}_*i_*\mathcal{O}_{Y^\mu_U} \cong  \oplus_{\kappa}   \overline{\mathcal{L}}_{U}^{\kappa,0}$
where the direct sum is over $\kappa = (\kappa_\tau)_{\tau\in \Sigma}$ such that $0 \le \kappa_\tau \le p-1$ for each $\tau$,
and $\kappa_\tau < p-1$ for some $\tau$ in each $\Sigma_v$.  Combined with the Kodaira--Spencer isomorphism on $\overline{Y}_U$,
we deduce from (\ref{eqn:dualising}) that
$$\overline{\pi}_*i_* \mathcal{K}^\mu_{U} \, \cong    \,
\mathcal{H}om_{\mathcal{O}_{\overline{Y}_{U}}}(\oplus_{\kappa}  \overline{\mathcal{L}}_U^{\kappa,0} ,  \overline{\mathcal{K}}_U )
   \,  \cong \, \oplus_{\kappa}  \overline{\mathcal{L}}_U^{2-\kappa,-1}.$$
Tensoring with $\overline{\mathcal{L}}_U^{m,0}$, we get injective morphisms $\overline{\mathcal{L}}_U^{k,-1}
\to \overline{\pi}_*i_* (\mathcal{K}^\mu_{U} \otimes_{\mathcal{O}_Y^\mu} i^*\overline{\mathcal{L}}^{m,0}_{U'})$
for $k = m+2 - \kappa$ as above.  Composing the homomorphism on sections with the one induced
by the inclusion $i_*\mathcal{K}_U^\mu \to \overline{\mathcal{K}}_{U'}$ obtained from (\ref{eqn:dualising}),
we obtain the desired injective homomorphism
$$H^0(\overline{Y}_U,\overline{\mathcal{L}}_U^{k,-1})  \to H^0(Y_U^\mu, \mathcal{K}^\mu_{U} \otimes_{\mathcal{O}_Y^\mu} i^*\overline{\mathcal{L}}^{m,0}_{U'})
 \to H^0(\overline{Y}_{U'}, \overline{\mathcal{K}}_{U'} \otimes_{\mathcal{O}_{\overline{Y}_{U'}}} \overline{\mathcal{L}}_{U'}^{m,0}).$$
Moreover one finds that for $g \in \mathrm{GL}_2(\mathbb{A}_F^\infty)$ with $g_p \in U_1(p)$, the isomorphisms
$\overline{\pi}_*i_*\mathcal{O}_{Y^\mu_U} \cong  \oplus_{\kappa}   \overline{\mathcal{L}}_{U}^{\kappa,0}$
are compatible with (\ref{eqn:pullback}) under the restriction of $\rho_g'$ to the subschemes $Y_U^\mu$,
and deduces that the maps $M_{k,-1}(U;E) \to M_{m+2,-1}(U';E)$ are compatible with the Hecke action;
in particular they commute with the operators $T_v$ and $S_v$ for $v\not\in Q$.

Next we show that if $m$ is sufficiently large, then the image of $M_{k,-1}(U;E)$ in  $M_{m+2,-1}(U';E)$
is contained in that of the reduction map from $M_{m+2,-1}(U';\mathcal{O})$ to $M_{m+2,-1}(U';E)$.
For this we will make use of the minimal compactifications $j:Y_U\to X_U$ and $j':Y_{U'}\to X_{U'}$ and
their properties recalled above.

We first compute the completion of $j'_*\mathcal{K}_{U'}$ along the cusps of $X_{U'}$.  We let 
 $j^{\mathrm{ord}}:Y_U^{\mathrm{ord}} \to X_U^{\mathrm{ord}}$ denote the restriction of $j$.
 Recall also the notation $\tilde{\pi}:X_{U'} \to X_U$, 
$\tilde{\pi}^{\mathrm{ord}}:X_{U'}^\mathrm{ord} \to X_U^\mathrm{ord} $ and $\pi^{\mathrm{ord}}:Y_{U'}^\mathrm{ord}
\to Y_U^\mathrm{ord}$ for the morphisms extending and restricting $\pi$.
Since $\pi^{\mathrm{ord}}$ is finite flat, we have
$$\begin{array}{rcl} \tilde{\pi}^{\mathrm{ord}}_*(j'_*\mathcal{K}_{U'})|_{X_{U'}^{\mathrm{ord}}}
 &=& j^{\mathrm{ord}}_*\pi^{\mathrm{ord}}_*(\mathcal{K}_{U'}|_{Y_{U'}^{\mathrm{ord}}}) \\
& \cong& j^{\mathrm{ord}}_*((\mathcal{H}om_{\mathcal{O}_{Y_U^{\mathrm{ord}}}}(\pi^{\mathrm{ord}}_*\mathcal{O}_{Y_{U'}^{\mathrm{ord}}},\mathcal{O}_{Y_U}^{\mathrm{ord}})\otimes_{\mathcal{O}_{Y_U^{\mathrm{ord}}}}( \mathcal{L}^{2,-1}_U|_{Y_U^{\mathrm{ord}}}))\\
&\cong& \mathcal{H}om_{j^{\mathrm{ord}}_*\mathcal{O}_{Y_U^{\mathrm{ord}}}}
    (j^{\mathrm{ord}}_*\pi^{\mathrm{ord}}_*\mathcal{O}_{Y_{U'}^{\mathrm{ord}}},j^{\mathrm{ord}}_*\mathcal{O}_{Y_U^{\mathrm{ord}}})\otimes_{\mathcal{O}_{X_U^{\mathrm{ord}}}}  (\mathscr{L}_U^2|_{X_U^{\mathrm{ord}}})\\
 &=& \mathcal{H}om_{\mathcal{O}_{X_U^{\mathrm{ord}}}}(\tilde{\pi}^{\mathrm{ord}}_*\mathcal{O}_{X_{U'}^{\mathrm{ord}}},
  \mathcal{O}_{X_U^{\mathrm{ord}}}) \otimes_{\mathcal{O}_{X_U^{\mathrm{ord}}}}  (\mathscr{L}_U^2|_{X_U^{\mathrm{ord}}}),
  \end{array}$$
where we made use of the canonical trivialisation of $\mathcal{L}_U^{0,1}$ and 
the Koecher Principle (for the last equality).  
Moreover the isomorphism is of $\tilde{\pi}^{\mathrm{ord}}_*\mathcal{O}_{X_{U'}^{\mathrm{ord}}}$-modules.

Since $\tilde{\pi}^{\mathrm{ord}}$ is finite, it follows that the completion of 
$\tilde{\pi}_*j'_*\mathcal{K}_{U'}\otimes_{\mathcal{O}_{X_{U}}}\mathscr{L}_{U}^{-2}$ along a cusp
 $C \subset X_U$ is canonically isomorphic to
the coherent sheaf on $\mathrm{Spf}\,\widehat{S}_C$ associated to the $\oplus \widehat{S}_{C'}$-module
$$\mathrm{Hom}_{\widehat{S}_C}  (\oplus \widehat{S}_{C'}, \widehat{S}_C),$$
where the direct sums are over the cusps
$C'$ of $X_{U'}$ in the pre-image of $C$ and the rings $\widehat{S}_C$ and $\widehat{S}_{C'}$
are defined by (\ref{eqn:RZ}) and (\ref{eqn:RZ1}) above.  Therefore the completion of 
$j'_*\mathcal{K}_{U'} \otimes_{\mathcal{O}_{X_{U'}}} \mathscr{L}_{U'}^{-2}$
along a cusp $C'$ of $X_{U'}$ 
is canonically isomorphic to $\mathrm{Hom}_{\widehat{S}_C}(\widehat{S}_{C'},\widehat{S}_C)$ as an 
$\widehat{S}_{C'}$-module if $C' \subset \tilde{\pi}^{-1}(C)$.  
  
Now consider the natural inclusion $\overline{j'_*\mathcal{K}}_{U'}  \to \overline{j}'_* \overline{\mathcal{K}}_{U'}$
of coherent sheaves on $\overline{X}_{U'}$,  where as usual we write
 $\overline{\,\cdot\,}$ for the special fibres of (quasi-coherent sheaves on and morphisms of) 
 schemes over $\mathcal{O}$.   This inclusion is an isomorphism on $\overline{Y}_{U'}$, so its cokernel
is supported on the cusps of $\overline{X}_{U'}$.  The same computation as above shows the completion of 
 $\overline{j}'_*\overline{\mathcal{K}}_{U'} \otimes_{\mathcal{O}_{\overline{X}_{U'}}}\overline{\mathscr{L}}_{U'}^{-2}$
 along $\overline{C}' \subset \overline{X}_{U'}$ is canonically isomorphic to the sheaf associated to
 $\mathrm{Hom}_{\widehat{S}_{\overline{C}}}(\widehat{S}_{\overline{C}'},\widehat{S}_{\overline{C}})$,
 where $\widehat{S}_{\overline{C}} = \widehat{S}_{C} \otimes_{\mathcal{O}}E$
 and $\widehat{S}_{\overline{C}'} = \widehat{S}_{C'} \otimes_{\mathcal{O}}E$.
Let  $X_U^\mu$ denote the
 closure of $Y_U^\mu$ in $\overline{X}_{U'}$,  so $X_U^\mu$ is a union of irreducible components of
$\overline{X}_{U'}$.  If $C'\subset \tilde{\pi}^{-1}(C)$ is a cusp of $X_{U'}$ such that  $\overline{C}' \subset X_U^\mu$, then
$\mathfrak{f} = O_F$, so $\widehat{S}_{C'}$ is flat over $\widehat{S}_C$ and the natural inclusion
$$\mathrm{Hom}_{\widehat{S}_C}(\widehat{S}_{C'},\widehat{S}_C) \otimes_{\mathcal{O}} E  \to
 \mathrm{Hom}_{\widehat{S}_{\overline{C}}}(\widehat{S}_{\overline{C}'},\widehat{S}_{\overline{C}})$$
is an isomorphism.  It follows that $\overline{j'_*\mathcal{K}}_{U'}  \to \overline{j}'_* \overline{\mathcal{K}}_{U'}$
is an isomorphism after completing along $C'$, and so an isomorphism on stalks at the (closed points of) cusps of $X_U^\mu$.
Therefore the cokernel of $\overline{j'_*\mathcal{K}}_{U'}  \to \overline{j}'_* \overline{\mathcal{K}}_{U'}$
is supported on the complement of $X_U^\mu$.  
It follows that $\overline{j}'_*$ of the inclusion $i_*\mathcal{K}_U^\mu \to \overline{\mathcal{K}}_{U'}$
factors through $\overline{j'_*\mathcal{K}}_{U'}$, and hence that the image of $M_{k,-1}(U;E)$ is contained in
the subspace
$$\begin{array}{rl}H^0(\overline{X}_{U'}, \overline{j'_*\mathcal{K}}_{U'}\otimes_{\mathcal{O}_{\overline{X}_{U'}}}
\overline{\mathscr{L}}_{U'}^m)
\subset 
&H^0(\overline{X}_{U'}, \overline{j}'_*\overline{\mathcal{K}}_{U'}\otimes_{\mathcal{O}_{\overline{X}_{U'}}}\overline{\mathscr{L}}_{U'}^m)
\\
=& H^0(\overline{X}_{U'}, \overline{j}'_*(\overline{\mathcal{K}}_{U'}\otimes_{\mathcal{O}_{\overline{Y}_{U'}}}
\overline{\mathcal{L}}_{U'}^m))\\
=& H^0(\overline{Y}_{U'}, \overline{\mathcal{K}}_{U'}\otimes_{\mathcal{O}_{\overline{Y}_{U'}}}\overline{\mathcal{L}}_{U'}^m))
=  M_{m+2,-1}(U';E).\end{array}$$

A key ingredient we need at this point is Theorem~E of~\cite{DKS}, which states that $R^i\pi_*\mathcal{K}_{U'} = 0$ for $i>0$,
so that in particular  $R^1\pi_*\mathcal{K}_{U'} = 0$.   Since $\tilde{\pi}^\mathrm{ord}$ is finite, it follows that
$R^1\tilde{\pi}_*(j'_*\mathcal{K}_{U'}) = 0$,
and hence the morphism $\tilde{\pi}_*j'_*\mathcal{K}_{U'} \to \tilde{\pi}_*(\overline{j'_*\mathcal{K}}_{U'})$
is surjective.  Since $\mathscr{L}_U$ is ample, we have
$H^1(X_U, \tilde{\pi}_*j'_*\mathcal{K}_{U'} \otimes_{\mathcal{O}_{X_U}}\mathscr{L}_U^m) =0$
for sufficiently large $m$, and it follows that the homomorphism
$$\begin{array}{ccc} H^0(X_U, \tilde{\pi}_* j'_*\mathcal{K}_{U'}  \otimes_{\mathcal{O}_{X_U}}\mathscr{L}_U^m)
 &\longrightarrow & H^0(X_U, \tilde{\pi}_* (\overline{j'_*\mathcal{K}}_{U'}) \otimes_{\mathcal{O}_{X_U}}\mathscr{L}_U^m))\\
&&\\ \parallel &  & \parallel \\ &&\\
 M_{m+2,-1}(U';\mathcal{O})  & \longrightarrow  &  H^0(\overline{X}_{U'}, \overline{j'_*\mathcal{K}}_{U'}
  \otimes_{\mathcal{O}_{\overline{X}_{U'}}}  \overline{\mathscr{L}}_{U'}^m))\end{array}
 $$
is surjective.   This completes the proof of the claim that the image of $M_{k,-1}(U;E)$ in  $M_{m+2,-1}(U';E)$
is contained in that of  $M_{m+2,-1}(U';\mathcal{O})$.

The theorem now follows from a standard argument.   Let $\mathbb{T}$ denote the ring of endomorphisms of
$M_{m+2,-1}(U';\mathcal{O})$ generated over $\mathcal{O}$ by the operators $T_v$ and $S_v$ for $v\not\in Q$.
Then $\mathbb{T}$ is a finite flat $\mathcal{O}$-algebra, and $M_{m+2,-1}(U';\mathcal{O})$ is a faithful $\mathbb{T}$-module
with $M_{k,-1}(U;E)$ as a subquotient.  The formula $Tf = \theta_f(T)f$ defines an $E$-algebra homomorphism
$\mathbb{T} \to E$ whose kernel is a maximal ideal $\mathfrak{m}$ generated by the operators $T_v - a_v$
and $S_v - d_v$ for $v\not\in Q$.  By the Going Down Theorem, there is a prime ideal $\mathfrak{p} \subset \mathfrak{m}$
such that $\mathfrak{p}\cap \mathcal{O} = 0$, and hence (enlarging $L$, $\mathcal{O}$ and $E$ if necessary), an $\mathcal{O}$-algebra
homomorphism $\tilde{\theta}:\mathbb{T} \to \mathcal{O}$ whose kernel is $\mathfrak{p} \subset \mathfrak{m}$.
Since $\mathfrak{p}$ is in the support of $M_{m+2,-1}(U';L) = M_{m+2,-1}(U';\mathcal{O})\otimes_{\mathcal{O}}L$, there
is an eigenform $\tilde{f} \in M_{m+2,-1}(U';L)$ such that $T\tilde{f}  = \tilde{\theta}(T)\tilde{f}$ for all $T \in \mathbb{T}$.  
By the existence of Galois representations associated to characteristic zero eigenforms~\cite{Ca} and~\cite{T} (together
with the usual association of reducible representations to Eisenstein series), we have a representation:
$$\rho_{\tilde{f}}:   G_F \longrightarrow \mathrm{GL}_2(L)$$
such that if $v\not\in Q$, then $\rho_{\tilde{f}}$ is unramified at $v$ and the characteristic
polynomial of $\rho_{\tilde{f}}(\mathrm{Frob}_v)$ is 
$$X^2 - a_v  X + d_v \mathrm{Nm}_{F/\mathbb{Q}}(v).$$
where $a_v = \tilde{\theta}(T_v)$ and $d_v = \tilde{\theta}(S_v)$.
Choosing a stable lattice and reducing modulo $\pi$ gives the desired representation $\rho_f$. This concludes the proof of Theorem~\ref{thm:galois}.  \epf

\begin{remark}  \label{rmk:det}

Note that by construction, if $\alpha \in F^\times \cap O_{F,p}^\times$, then 
$\left(\begin{array}{cc}  \alpha & 0 \\ 0 & \alpha \end{array}\right)$ acts on $M_{k,l}(E)$
as $\overline{\alpha}^{k+2l-2}$.  Therefore if $f \in M_{k,l}(U;E)$ is an eigenform
for $S_v$ with eigenvalue $d_v$ for all $v\not\in Q$, then there is a character 
$$\psi: (\mathbb{A}_F^\infty)^\times /(U\cap (\mathbb{A}_F^\infty)^\times) \to E^\times$$
such that $\psi(\alpha) = \overline{\alpha}^{k+2l-2}$ for all $\alpha \in F^\times \cap O_{F,p}^\times$
and $\psi(\varpi_v) = d_v$ for all $v\not\in Q$.
It follows from the description of $\rho_f$ in Theorem~\ref{thm:galois} that 
$\det(\rho_f)\chi_{\mathrm{cyc}}$ (where $\chi_{\mathrm{cyc}}$ is the cyclotomic character)
corresponds via class field theory to the character
$$\psi': \mathbb{A}_F^\times /F^\times F_\infty^\times V \to E^\times$$
defined by $\psi'(\alpha z a) = \psi(a)\overline{a}_p^{2-k-2l}$ for $\alpha \in F^\times, z \in F_{\infty}^\times$
and $a \in (\mathbb{A}_F^\infty)^\times$ with $a_p \in O_{F,p}^\times$, where
$V = \{\,a\in \mathbb{A}_F^\infty\,|\, a \in U, a_p \equiv 1 \bmod p\,\}$.\end{remark}

\section{Geometric weight conjectures}
In this section we formulate our geometric Serre weight conjectures and discuss
the relation with \cite{BDJ}.

\label{sec:conj}

\subsection{Geometric modularity}
Let 
$$\rho: G_F = \mathrm{Gal}(\overline F/F)\rightarrow \mathrm{GL}_2(\overline{\mathbb{F}}_p)$$
be an irreducible, continuous, totally odd representation of the absolute Galois group of $F$.

\begin{definition} We say that $\rho$ is {\em geometrically modular of weight $(k,l)$} if
$\rho$ is equivalent to the extension of scalars of $\rho_f$
for some open compact subgroup $U\subset \mathrm{GL}_2(\widehat{O}_F)$
and eigenform $f \in M_{k,l}(U;E)$ as in the statement of Theorem~\ref{thm:galois}.\end{definition}

Note that the level $U$ is unspecified, but required to contain $\mathrm{GL}_2(O_{F,p})$.  
Also unspecified are the field $E$ (and thus implicitly the field $L \subset \overline{\mathbb{Q}}_p$,
by which we view $E = \mathcal{O}/\pi \subset \overline{\mathbb{F}}_p$) and the finite
set of primes of $Q$.  Thus $\rho$ is geometrically modular of weight $(k,l)$ if there
is a non-zero element $f \in M_{k,l}(U;E)$ for some $U\supset \mathrm{GL}_2(O_{F,p})$
and $E\subset \overline{\mathbb{F}}_p$ such that
$$T_v f = \mathrm{tr}(\rho(\mathrm{Frob}_v))f\quad\mbox{and} \quad\mathrm{Nm}_{F/\mathbb{Q}}(v) S_v f = \det(\rho(\mathrm{Frob}_v))f$$
for all but finitely many primes $v$.  (Note that both sides of both equations are defined whenever
$v\nmid p$, $\mathrm{GL}_2(O_{F,v}) \subset U$ and $\rho$ is unramified at $v$.)

\begin{remark}  \label{rmk:folklore}
Folklore conjectures predict that every $\rho$ as above is indeed geometrically modular of {\em some}
weight $(k,l)$.  The focus of this paper is to give a conjectural recipe for {\em all} such weights $(k,l)$ in terms
of the local behaviour of $\rho$ at primes  over $p$.\end{remark}

\subsection{Crystalline lifts}  \label{subsec:cryslift}
In order to formulate our conjectures, we recall the notion of labelled Hodge-Tate weights.
Let $K$ be a finite extension of $\mathbb{Q}_p$, and let
$$\sigma:G_K \to \mathrm{GL}_d(L)  = \mathrm{Aut}_L(V)$$
be a continuous representation on a $d$-dimensional $L$-vector space $V$.
Recall that $V$ is  crystalline
if $D_{\mathrm{crys}}(V) = (V \otimes_{\mathbb{Q}_p} B_{\mathrm{crys}})^{G_K}$
is free of rank $d$ over 
$$(L\otimes_{\mathbb{Q}_p} B_{\mathrm{crys}})^{G_K} = L\otimes_{\mathbb{Q}_p} K_0$$
where $B_{\mathrm{crys}}$ is Fontaine's ring of crystalline periods \cite{F} and $K_0$ is the maximal unramified
subfield of $K$.  One similarly defines the notion of a {\em de Rham} (resp.~{\em Hodge--Tate}) representation
and an associated filtered (resp.~graded) free
module $D_{\mathrm{dR}}(V)$ (resp.~$D_{\mathrm{HT}}(V)$) of rank $d$ over $L\otimes_{\mathbb{Q}_p} K$
in terms of the rings $B_{\mathrm{dR}}$ (resp.~$B_{\mathrm{HT}}$).
Moreover if
$V$ is crystalline, then it is de Rham, and if $V$ is de Rham then it is also Hodge--Tate.
Thus if $V$ is crystalline, then $D_{\mathrm{HT}}(V)$ is a graded free module of rank $d$
over $L\otimes_{\mathbb{Q}_p} K$.  If $L$ is sufficiently large that it contains the image
of each embedding of $K$ into $\overline{\mathbb{Q}}_p$, then $L\otimes_{\mathbb{Q}_p} K
\cong \prod_{\tau \in \Sigma_K} L$ where $\Sigma_K = \{\,\tau:K\to L\,\}$, and for each 
$\tau \in \Sigma_K$, the corresponding component of $D_{\mathrm{HT}}(V)$
is a graded $d$-dimensional vector space over $L$.  
\begin{definition}  \label{def:HTtype}
If $V$ is crystalline, then the $\tau$-{\em labelled weights}
of $V$ are defined as the $d$-tuple of integers $(w_1,w_2,\ldots,w_d) \in \mathbb{Z}^d$ such that
$w_1 \ge w_2 \ge \cdots \ge w_d$
and the $\tau$-component of $D_{\mathrm{HT}}(V)$ is isomorphic to $\oplus_{i=1}^d L[w_i]$,
where $L[w_i]$ has degree $w_i$. 
We define the {\em Hodge--Tate type} of $V$ to be the element of $(\mathbb{Z}^d)^{\Sigma_K}$
whose $\tau$-component is given by the $\tau$-labelled weights of $V$; thus to give the Hodge--Tate
type of $V$ is equivalent to giving the isomorphism class of $D_{\mathrm{HT}}(V)$ as a
graded $L\otimes_{\mathbb{Q}_p} K$-module.
\end{definition}

We now specialise to the case $d=2$ and $K=F_v$ where $v$ is a prime of $F$ dividing $p$,
so $\Sigma_K$ is identified with the subset $\Sigma_v  \subset \Sigma = \{\,\tau:F \to L \,\}$
defined by (\ref{eqn:Sigmav}), and consider the representation 
$$\sigma: G_K \to \mathrm{GL}_2(\overline{\mathbb{F}}_p).$$

\begin{definition}  \label{def:weightlift} For a pair $(k,l) \in \mathbb{Z}_{\ge 1}^{\Sigma_v} \times \mathbb{Z}^{\Sigma_v}$,
we say that $\sigma$ has a {\em crystalline lift} of weight $(k,l)$ if for some
sufficiently large extension $L\subset\overline{\mathbb{Q}}_p$ of $\mathbb{Q}_p$ with ring $\mathcal{O}$ of integers and residue field 
$E\subset\overline{\mathbb{F}}_p$, there exists a continuous representation:
$$\tilde{\sigma}:  G_K \to  \mathrm{GL}_2(\mathcal{O})$$
such that $\tilde{\sigma}\otimes_{\mathcal{O}} \overline{\mathbb{F}}_p$ is isomorphic 
to $\sigma$, and $\tilde{\sigma}\otimes_{\mathcal{O}}L$ is crystalline with
Hodge--Tate type $(k+l-1,l)$.\end{definition}

\subsection{Statement of the conjectures}
First recall from \S\ref{subsec:minimal} the definition of the minimal cone:
$$\Xi_{\mathrm{min}} = \{\, k \in \mathbb{Z}^\Sigma \, |\, \mbox{$pk_\tau \ge k_{\mathrm{Fr}^{-1}\circ\tau}$ for all $\tau \in \Sigma$}\,\},$$
and let $\Xi_{\mathrm{min}}^+ = \Xi_{\mathrm{min}} \cap \mathbb{Z}_{\ge 1}^\Sigma$.

\begin{conjecture} \label{conj:geomweights}
Let $\rho: G_F \rightarrow \mathrm{GL}_2(\overline{\mathbb{F}}_p)$
be an irreducible, continuous, totally odd representation, and let $l \in \mathbb{Z}^\Sigma$.
There exists $k_{\mathrm{min}} = k_{\mathrm{min}}(\rho,l) \in \Xi_{\mathrm{min}}^+$ such that the following hold:
\begin{enumerate}
\item $\rho$ is geometrically modular of weight $(k,l)$ if and only if $k \ge_{\mathrm{Ha}} k_{\mathrm{min}}$;
\item if $k \in \Xi_{\mathrm{min}}^+$, then $k \ge_{\mathrm{Ha}}  k_{\mathrm{min}}$ if and only if $\rho|_{G_{F_v}}$
has a crystalline lift of weight of $(k_\tau,l_\tau)_{\tau\in \Sigma_v}$ for all $v|p$.
\end{enumerate}
\end{conjecture}

Note that the conjecture, in particular the existence of $k_{\mathrm{min}}$ as in 1), incorporates the ``folklore conjecture''
(see Remark~\ref{rmk:folklore}) that $\rho$ is geometrically modular of {\em some} weight $(k,l)$.  Moreover, for any $l \in \mathbb{Z}^\Sigma$
there should be weights $(k,l)$ for which $\rho$ is geometrically modular.  In fact one can show using partial
$\Theta$-operators (defined below in \S\ref{sec:Theta})
 that for any given $l$ and $l'$, if $\rho$ is irreducible and geometrically modular of some weight $(k,l)$,
then $\rho$ is geometrically modular of some weight $(k',l')$.  
We explain this, and the dependence of $k_{\mathrm{min}}$ on $l$ (for fixed $\rho$), in \S\ref{subsec:consequences}. 
A simpler observation is that the conjecture is compatible with twists by arbitrary
characters $\xi:G_F \to \overline{\mathbb{F}}_p^\times$.    More precisely, by Lemma~\ref{lem:twist} and the
well-known computation of reductions of crystalline characters (see for example \cite[Prop.~B4]{Co}), we see that the conjecture holds
for the pair $(\rho,l)$ if and only if it holds for the pair $(\rho\otimes\xi,l-m)$ for any $m\in \mathbb{Z}^\Sigma$
such that $\xi|_{I_{F_v}} = \prod_{\tau \in \Sigma_v} \epsilon_\tau^{m_\tau}$ for all $v|p$, where $I_{F_v}$ is the inertia
subgroup of $G_{F_v}$ and 
$$\epsilon_\tau : I_{F_v}  \longrightarrow O_{F,v}^\times \longrightarrow \overline{\mathbb{F}}_p^\times$$
is the fundamental character defined as the composite of the maps induced by $\tau$ and local class field theory.
Thus Conjecture~\ref{conj:geomweights} for all pairs $(\rho,l)$ reduces to the case $l=0$, with the resulting minimal weights
related by $k_{\mathrm{min}}(\rho,l) = k_{\mathrm{min}}(\rho\otimes\xi,0)$ for any character $\xi$
chosen so that $\xi|_{I_{F_v}} = \prod_{\tau \in \Sigma_v} \epsilon_\tau^{l_\tau}$ for all $v|p$.
 We remark also that $\rho|_{G_{F_v}}$ always has a crystalline lift of some weight 
$(k_\tau,l_\tau)_{\tau\in \Sigma_v}$ with $2 \le k_\tau \le p+1$ for all $\tau \in \Sigma_v$, %
from which it follows that $\rho$ has a twist for which $k_{\mathrm{min}}$ as in 2) would satisfy
$k_\tau \le p + 1$ all $\tau \in \Sigma$.

Assuming that $\rho$ is geometrically modular of some weight, then the existence of a weight
$k_{\mathrm{min}}$ satisfying 1) in Conjecture~\ref{conj:geomweights}
is strongly suggested by Corollary~1.2 of \cite{DK}, which implies that
the minimal weight $\nu(f)$  of the eigenform $f$ satisfies $\nu(f) \in \Xi_{\mathrm{min}}$, but it is
not an immediate consequence.   Indeed there are two issues: firstly, we would need
$\nu(f) \in \Xi_{\mathrm{min}}^+$ (which we expect to hold if $\rho_f$ is irreducible), and secondly,
the eigenform $f$ giving rise to $\rho$ is not unique.  %
However if we grant the existence of $k_{\mathrm{min}}$ as in 1), then Conjecture~\ref{conj:geomweights}
reduces to the following:
\begin{conjecture} \label{conj:geomweights2}
Suppose that $\rho: G_F \rightarrow \mathrm{GL}_2(\overline{\mathbb{F}}_p)$ is irreducible and
geometrically modular some weight, and that $k \in \Xi_{\mathrm{min}}^+$.  Then
$\rho$ is geometrically modular of weight $(k,l)$ if and only if $\rho|_{G_{F_v}}$
has a crystalline lift of weight of $(k_\tau,l_\tau)_{\tau\in \Sigma_v}$ for all $v|p$.
\end{conjecture}

\begin{remark} The existence of $k_{\mathrm{min}}$ satisfying part 2) of Conjecture~\ref{conj:geomweights}
is a purely $p$-adic Hodge-theoretic statement,  and it is strongly suggested by the Breuil--M\'ezard
Conjecture (of~\cite{BM} as generalised by~\cite{GKi})
and the modular representation theory of $\mathrm{GL}_2(O_F/p)$, but again not an immediate
consequence.  We remark also that the condition $k \in \Xi_{\mathrm{min}}^+$ is needed; indeed R.~Bartlett
has constructed local Galois representations with crystalline lifts of weight\footnote{For notational consistency,
assume here that $p$ is inert in $F$.} $(k,l)$, but none of weight 
$(k',l)$, where $k' = k+ k_{\mathrm{Ha}_\tau}$ is in $\mathbb{Z}_{\ge 1}^{\Sigma}$ but not in $\Xi_{\mathrm{min}}^+$.
Granting the existence of a weight $k_{\mathrm{min}}$ as in 2), then Conjecture~\ref{conj:geomweights}
follows from Conjecture~\ref{conj:geomweights2} under the assumptions that $\rho$ is geometrically
modular of some weight and that $\nu(f) \in \Xi_{\mathrm{min}}^+$ if $\rho_f \sim \rho$ is irreducible.
\end{remark}

\subsection{The case $k=1$}
We now consider a special case of Conjecture~\ref{conj:geomweights2}.  Since a representation
$G_K \to \mathrm{GL}_d(L)$ is crystalline of Hodge--Tate type $0 \in (\mathbb{Z}^d)^{\Sigma_K}$ if and
only if it is unramified, it follows that $\sigma:G_K \to \mathrm{GL}_2(\overline{\mathbb{F}}_p)$ has
a crystalline lift of weight $(1,0)$ if and only if it is unramified.  Thus Conjecture~\ref{conj:geomweights2}
incorporates the prediction that $\rho$, assumed to be geometrically modular, is of weight $(1,0)$ if and
only if it is unramified at all primes $v|p$.  One direction of this, that if $\rho$ is geometrically modular
of weight $(1,0)$ then it is unramified at all $v|p$, is a theorem of Dimitrov and Wiese~\cite{DW}
(also proved independently by Emerton, Reduzzi and Xiao~\cite{ERX2} under additional hypotheses),
and the other direction is proved under technical hypotheses by Gee and Kassaei~\cite{GK}.
By twisting, these results extend to the case of weight $(1,l)$ for aribtrary $l$.

\subsection{Relation to algebraic modularity}
We now explain how our conjecture is consistent with results on the weight part of Serre's
Conjecture as formulated by Buzzard, Jarvis and one of the authors in \cite{BDJ}.   These results
provide information about {\em algebraic weights}, meaning weights $(k,l)$ such that $k_\tau \ge 2$
for all $\tau$, but with a different notion of modularity, which we call {\em algebraic modularity}.  
We will next explain this notion and its relation with the conjectures above.  The remainder of the
paper will then focus on developing methods applicable to the case of partial weight one, which lies
outside both settings just mentioned, namely weights that are algebraic or of the form $(1,l)$.

Recall that in \cite{BDJ}, a {\em Serre weight} is an irreducible representation of $\mathrm{GL}_2(O_F/p)$
over $\overline{\mathbb{F}}_p$.   For an algebraic weight $(k,l) \in \mathbb{Z}^\Sigma_{\ge 2} \times \mathbb{Z}^\Sigma$,
we let $V_{k,l}$ denote the representation 
$$\bigotimes_{\tau\in \Sigma}  \left(\det{}^{l_\tau}  \otimes \mathrm{Sym}^{k_\tau-2}\, \overline{\mathbb{F}}_p^2 \right),$$
where $\mathrm{GL}_2(O_F/p)$ acts on the factor indexed by $\tau$ via the homomorphism to
$\mathrm{GL}_2(\overline{\mathbb{F}}_p)$ induced  by $\tau$.   The irreducible representations of $\mathrm{GL}_2(O_F/p)$
(i.e., Serre weights) are precisely the $V_{k,l}$ such that $2 \le k_\tau \le p+1$ for all $\tau \in \Sigma$;
moreover for such $(k,l)$, we have that $V_{k,l}$ is isomorphic to
$V_{k',l'}$ if and only if $k = k'$ and $l-l' \in \bigoplus_{\tau\in \Sigma} \mathbb{Z}\cdot\mathrm{Ha}_\tau$.
(More concretely, the latter condition means that $\sum_{i=0}^{f_\tau-1} l_{\mathrm{Fr}^i\circ\tau} p^i 
\equiv \sum_{i=0}^{f_\tau-1} l'_{\mathrm{Fr}^i\circ\tau} p^i \bmod (p^{f_\tau}-1)$ for all $\tau\in \Sigma$, where
$f_\tau = [\tau(O_F):\mathbb{F}_p]$.)

 For an irreducible representation $\rho:G_F\to \mathrm{GL}_2(\overline{\mathbb{F}}_p)$ and an aribtrary finite-dimensional
 representation $V$ of $\mathrm{GL}_2(O_F/p)$ over $\overline{\mathbb{F}}_p$, we say $\rho$ is
 {\em modular of weight} $V$ if it arises in the \'etale cohomology of a suitable quaternionic Shimura curve
 over $F$ with coefficients in a lisse sheaf associated to $V$; we refer the reader to Section~2 of~\cite{BDJ}
 for the precise definition.\footnote{Alternatively, but not a priori equivalently, one can define the notion of
 modularity of weight $V$ in terms of the presence of the corresponding system of Hecke eigenvalues
 on spaces of mod $p$ automorphic forms on totally definite quaternion algebras over $F$.}
 It is also proved in {\em loc.~cit.} that $\rho$ is modular of weight $V$ if and only
 if it is modular of weight $W$ for some Jordan--H\"older consitituent $W$ of $V$, so the determination
 of the weights $V$ for which $\rho$ is modular reduces to the consideration of Serre weights.
 \begin{definition} \label{def:algmod}
 For an algebraic weight $(k,l) \in \mathbb{Z}^\Sigma_{\ge 2} \times \mathbb{Z}^\Sigma$, we will
 say that $\rho$ is {\em algebraically modular} of weight $(k,l)$ if it is modular of weight $V_{k,1-k-l}$
 in the sense of \cite{BDJ} (the presence of the twist being to reconcile the conventions
 of this paper with the ones of \cite{BDJ}).  
 \end{definition}
 
 A conjecture is formulated in \cite{BDJ} for the set of Serre weights for which $\rho$ is modular.
 Under the assumption that $\rho$ is algebraically modular of some weight and mild technical hypotheses,
 the conjecture is proved in a series of papers by Gee and coauthors, culminating in~\cite{GLS} and~\cite{GKi},
 with an independent alternative to the latter (deducing the conjecture from its analogue in the context of
certain unitary groups) provided by Newton~\cite{N}.  They also prove variants of
 the conjecture (under the same hypotheses), including that if $2 \le k_\tau \le p+1$ for all $\tau$, then
 $\rho$ is algebraically modular of weight $(k,l)$ if and only if $\rho|_{G_{F_v}}$
has a crystalline lift of weight of $(k_\tau,l_\tau)_{\tau\in \Sigma_v}$ for all $v|p$.
The generalised Breuil--M\'ezard Conjecture (as in \cite{GKi}) would imply
that this result extends to arbitrary algebraic weights.   We are therefore led to conjecture:

\begin{conjecture} \label{conj:algvsgeomweights}
Let $\rho: G_F \rightarrow \mathrm{GL}_2(\overline{\mathbb{F}}_p)$
be an irreducible, continuous, totally odd representation, and let 
$(k,l) \in \mathbb{Z}^\Sigma_{\ge 2} \times \mathbb{Z}^\Sigma$. 
If $\rho$ is algebraically modular of weight $(k,l)$, then $\rho$ is geometrically modular of weight $(k,l)$.
Moreover, if in addition $k \in \Xi_{\mathrm{min}}^+$, then the converse holds.
\end{conjecture}

\begin{remark} \label{rmk:TobyQ} The assumption $k \in \Xi_{\mathrm{min}}^+$ appears for reasons related to its presence
in part 2) of Conjecture~\ref{conj:geomweights}; here however one can see its necessity more readily from
modular representation-theoretic considerations.  Indeed if $k,k' \in \mathbb{Z}^\Sigma_{\ge 2}$ with
$k' = k + k_{{\mathrm{Ha}}_\tau} \not\in \Xi_{\mathrm{min}}^+$, then $V_{k,l}$ may have Jordan--H\"older constituents
not present in $V_{k',l}$, so there are representations which are algebraically modular of weight $(k,l)$,
but not of weight $(k',l)$.    Note though that if $2 \le k_\tau \le  p+1$ for all $\tau$,
then $k \in \Xi_{\mathrm{min}}^+$, so we conjecture that algebraic and geometric modularity are equivalent for
weights associated to Serre weights.
\end{remark}

From our construction of the Galois representation associated to an eigenform $f$,
we see that $\rho_f$ is the reduction of some representation associated to a characteristic zero eigenform,
from which it follows (e.g. from \cite[Prop.~2.10]{BDJ}) that $\rho_f$ is modular of some weight $V$.
Thus if $\rho$ is geometrically modular of some weight, then it is algebraically modular of
some weight.  

Conversely, suppose that $\rho$ is algebraically modular of some {\em paritious} weight 
$(k,l) \in \mathbb{Z}_{\ge 2}^\Sigma \times \mathbb{Z}^\Sigma$ (see Definition~\ref{def:paritious}).  Then \cite[Prop.~2.5]{BDJ}
implies that $\rho$ is the reduction of some representation associated to a characteristic zero
eigenform of weight $(k,l)$ and level prime to $p$, and hence that $\rho$ is geometrically
modular of weight $(k,l)$.   More generally, if $(k,l)$ is any algebraic weight such $k_\tau \equiv k_{\tau'}\bmod 2$
for all $\tau,\tau'\in \Sigma$, then we can choose $l'$ so that $(k,l')$ is paritious and a character
$\xi$ so that $\xi|_{I_{F_v}} = \prod_{\tau \in \Sigma_v} \epsilon_\tau^{l_\tau-l'_\tau}$.
If $\rho$ is algebraically modular of weight $(k,l)$, then $\rho\otimes\xi$ is algebraically
modular of weight $(k,l')$ (by \cite[Prop.~2.11]{BDJ}), so the above argument shows that
$\rho\otimes\xi$ is geometrically modular of weight $(k,l')$ and hence that $\rho$ is geometrically
modular of weight $(k,l)$.  We have thus proved the following:

\begin{proposition} \label{prop:algvsgeomweights} If $\rho$ is geometrically modular of some weight, then it is algebraically
modular of some (algebraic) weight.  Conversely, if $\rho$ is algebraically modular of some
algebraic weight $(k,l)$ such that $k_\tau \equiv k_{\tau'}\bmod 2$ for all $\tau,\tau'\in \Sigma$, then
$\rho$ is geometrically modular of the same weight $(k,l)$.
\end{proposition}

\section{$\Theta$ operators}
\label{sec:Theta}

In this section we recall the definition due to Andreatta and Goren of partial $\Theta$-operators (see \cite[\S12]{AG}),
with some simplifications and adaptations to our setting.  See also~\cite{PRtheta}, where the approach taken here
is applied in the context of Pappas--Rapoport models defined in \cite{PR}, leading to a substantial refinement of the
results in~\cite{AG} when $p$ is ramified in $F$.

\subsection{Igusa level structure}
We assume that $U$ is $p$-neat, as in Definition~\ref{def:pneat}.  For $\tau \in \Sigma$, we will write $\overline{\omega}_\tau$ 
(resp.~$\overline{\delta}_\tau$)
for the line bundle  $\overline{\mathcal{L}}^{k,0}_U$ (resp.~$\overline{\mathcal{L}}^{0,k}_U$) on $\overline{Y}_U$,
where $k$ is such that $k_\tau = 1$ and $k_{\tau'} = 0$ for $\tau' \neq \tau$.
We view the partial Hasse invariant
$$\mathrm{Ha}_\tau \in 
H^0(\overline{Y}_U,\overline{\mathcal{L}}_U^{k_{\mathrm{Ha}_\tau},0})
= H^0(\overline{Y}_U,\overline{\omega}_\tau^{-1}\otimes \overline{\omega}_{\mathrm{Fr}^{-1}\circ\tau}^p)$$
as a morphism $\overline{\omega}_{\mathrm{Fr}^{-1}\circ\tau}^{-p}\to  \overline{\omega}_\tau^{-1}$.
For each $v|p$, we let $\mathrm{Ha}_v = \prod_{\tau\in \Sigma_v} \mathrm{Ha}_\tau$, which we view
as a morphism $(\otimes_{\tau\in \Sigma_v} \overline{\omega}_{\mathrm{Fr}^{-1}\circ\tau} ^{-p})\otimes(\otimes_{\tau\in \Sigma_v} \overline{\omega}_\tau)\rightarrow \mathcal{O}_{\overline{Y}_U}$, i.e., $\otimes_{\tau\in \Sigma_v} \overline{\omega}_{\tau}^{1-p} \to  \mathcal{O}_{\overline{Y}_U}$
(where $\Sigma_v$ is defined in (\ref{eqn:Sigmav})).

We define the scheme
$$Y_{U}^{\mathrm{Ig}}  =  \mathbf{Spec}  \left(\mathrm{Sym}_{\mathcal{O}_{\overline{Y}_U}}
(\oplus_{\tau\in\Sigma} \overline{\omega}_\tau^{-1})/\mathcal{I}\right),$$
where $\mathcal{I}$ is the sheaf of ideals of 
$\mathrm{Sym}_{\mathcal{O}_{\overline{Y}_U}}(\oplus_{\tau\in\Sigma} \overline{\omega}_\tau^{-1})$
generated by the sheaves of $\mathcal{O}_{\overline{Y}_U}$-submodules
$$(\mathrm{Ha}_\tau-1)\omega_{\mathrm{Fr}^{-1}\circ\tau}^{-p}\,\mbox{for $\tau \in \Sigma$}, \quad
  (\mathrm{Ha}_v-1)\left(\otimes_{\tau\in \Sigma_v} \omega_{\tau}^{1-p}\right)\,\mbox{for $v|p$.}$$
We define an action of $(O_F/pO_F)^\times$ on $Y_{U}^{\mathrm{Ig}}$ over $\overline{Y}_U$
by having $\alpha \in (O_F/pO_F)^\times$ act on the structure sheaf as the $\mathcal{O}_{\overline{Y}_U}$-algebra
automorphism defined by multiplication by $\tau(\alpha)^{-1}$ on the summand $\overline{\omega}_\tau^{-1}$.
(Note that the action is well-defined since the $\mathrm{Ha}_\tau$ are invariant under this action and hence 
$\mathcal{I}$ is preserved.)
\begin{proposition}\label{prop:igusa}  
Let $\pi_U:Y_{U}^{\mathrm{Ig}} \to \overline{Y}_U$ denote the natural projection.  Then
\begin{enumerate}
\item The morphism $\pi_U$ is finite and flat, and identifies $\overline{Y}_U$ with the quotient
of $Y_{U}^{\mathrm{Ig}}$ by the action of $(O_F/pO_F)^\times$.
\item The restriction of $\pi_U$ to the preimage of $Y_U^{\mathrm{ord}}$ is \'etale.
\item The scheme $Y_{U}^{\mathrm{Ig}}$ is normal.
\end{enumerate}
\end{proposition}
\begpf Each assertion can be checked over affine open subschemes 
$V \subset \overline{Y}_U$ on which the line bundles $\overline{\omega}_\tau^{-1}$
are trivial.  For each $\tau\in \Sigma$, let $x_\tau$ be a generator of
$M_\tau = \Gamma(V,\overline{\omega}_\tau^{-1})$
over $R= \Gamma(V,\mathcal{O}_{\overline{Y}_U})$.
Then $\mathrm{Ha}_\tau(x^p_{\mathrm{Fr}^{-1}\circ\tau}) = r_\tau x_\tau$ for some
$r_\tau \in R$, and $\pi_U^{-1}(R) = \mathrm{Spec}\,T$ where
$$T = R[x_\tau]_{\tau \in \Sigma} / \langle\,
\mbox{$x_{\mathrm{Fr}^{-1}\circ\tau}^p - r_\tau x_\tau$ for $\tau\in \Sigma$, 
 $\displaystyle\prod_{\tau\in\Sigma_v} x_\tau^{p-1} - \displaystyle\prod_{\tau\in\Sigma_v} r_\tau$ for $v|p$}\,\rangle.$$
Thus $T$ is free over $R$ with basis $\left\{ \, \prod_{\tau\in \Sigma}  t_\tau^{\kappa_\tau} \, \right\}$,
where $t_\tau$ denotes the image of $x_\tau$ in $T$ and the tuples $\kappa = (\kappa_\tau)_{\tau\in \Sigma}$ are those
satisfying
\begin{itemize}
\item $0 \le \kappa_\tau \le p-1$ for each $\tau \in \Sigma$, 
\item and $\kappa_\tau < p-1$ for some $\tau$ in each $\Sigma_v$.
\end{itemize}
Note that $(O_F/pO_F)^\times$ acts on $ \prod_{\tau\in \Sigma}  t_\tau^{\kappa_\tau}$ by the character
$\prod_{\tau\in\Sigma} \tau^{-\kappa_\tau}$, and these are precisely the distinct characters of the 
$(O_F/pO_F)^\times$.  Therefore $T^{(O_F/pO_F)^\times} = R$, and 1) follows.

To prove 2), recall that $Y_U^{\mathrm{ord}}$ is the complement of $\cup_{\tau\in \Sigma} Z_{U,\tau}$
where $Z_{U,\tau}$ is vanishing locus of $\mathrm{Ha}_\tau$ on $\overline{Y}_U$.  We must therefore show that if
all $r_\tau$ are invertible in $R$, then $T$ is \'etale over $R$.  From the above description
of $T$, we see that 
$$r_\tau dt_\tau = d(t_{\mathrm{Fr}^{-1}\circ\tau}^p - r_\tau t_\tau) = 0$$
in $\Omega^1_{T/R}$.  It follows that $dt_\tau = 0$ for all $\tau$, and hence
$\Omega^1_{T/R} = 0$, so $T$ is \'etale over $R$.

To prove 3), we use Serre's Criterion.  Since $R$ is regular and $T$ is finite and flat over $R$,
$T$ is Cohen--Macaulay, so it suffices to prove that $T$ is regular in codimension 1.
Thus it suffices to prove that the semi-local ring $T_{\mathfrak{p}}$  is regular for every
height one prime $\mathfrak{p}$ of $R$.  If $\mathfrak{p} \in Y_U^{\mathrm{ord}}$,
then $T_{\mathfrak{p}}$ is \'etale over $R_{\mathfrak{p}}$, so $T_{\mathfrak{p}}$
is regular.  Otherwise $\mathfrak{p}$ defines an irreducible component of 
$V \cap Z_{U,\tau}$ for some $\tau = \tau_0 \in \Sigma$.  Since $V \cap Z_{U,\tau}$ is
defined by $r_\tau$ and $\sum Z_{U,\tau}$ is reduced, the DVR
$R_{\mathfrak{p}}$ has uniformiser $r_{\tau_0}$, and $r_\tau \in R_{\mathfrak{p}}^\times$
for $\tau \neq \tau_0$.   Letting $\tau_0 \in \Sigma_{v_0}$, $f = \#\Sigma_{v_0}$ and
$$S = R[x_\tau]_{\tau\in \Sigma_{v_0}}/  \langle\,
\mbox{$x_{\mathrm{Fr}^{-1}\circ\tau}^p - r_\tau x_\tau$ for $\tau\in \Sigma_{v_0}$, 
 $\displaystyle\prod_{\tau\in\Sigma_{v_0}} x_\tau^{p-1} - \displaystyle\prod_{\tau\in\Sigma_{v_0}} r_\tau$}\,\rangle,$$
the formulas $t_{\mathrm{Fr}^{i}\circ \tau} = r_{\mathrm{Fr}^{i}\circ \tau}^{-1}t_{\mathrm{Fr}^{i-1}\circ \tau}^p$
for $i=1,\ldots,f-1$ show that 
$$S_{\mathfrak{p}} = R_{\mathfrak{p}}[x_{\tau_0}]/\langle\, x_{\tau_0}^{p^f-1} - r_{\tau_0} \prod_{i=1}^{f-1}  r_{\mathrm{Fr}^{f-i}\circ\tau_0}^{p^i}\rangle,$$
which is a DVR with uniformiser $t_{\tau_0}$.  Since $r_\tau$ is invertible in $S_{\mathfrak{p}}$ for $\tau \not\in \Sigma_{v_0}$,
we see as above that $T_{\mathfrak{p}}$ is \'etale over $S_{\mathfrak{p}}$, and is therefore also regular.\epf

\begin{remark}

We could similarly have defined schemes $Y_{J,N}^{\mathrm{Ig}}$ as above by
replacing $\overline{Y}_U$ with $\overline{Y}_{J,N}$.  Then $Y_{J,N}^{\mathrm{Ig}}$ is isomorphic
to the closed subscheme (in fact a union of irreducible components) of $\overline{Y}^1_{J,N}$ for
which the subgroup scheme $H \subset A[p]$ is generically \'etale.
However under this isomorphism, the natural projection $Y_{J,N}^{\mathrm{Ig}}\to \overline{Y}_{J,N}$
corresponds to the restriction of the morphism $\overline{Y}_{J,N}^1\to \overline{Y}_{J,N}$
defined by $(\underline{A},H,P) \mapsto \underline{A/H}$.  Furthermore, we can realise $Y^{\mathrm{Ig}}_U$
as the quotient of $\coprod Y_{J,N}^{\mathrm{Ig}}$ by the action of $G_{U,N}$ obtained from
the one on $\overline{Y}_{J,N}^1$ defined by $(\nu,u)\cdot(\underline{A},H,P) = ((\nu,u)\cdot\underline{A},H,\nu P)$;
as this differs from the one already defined, it does not yield an identification of $Y_U^{\mathrm{Ig}}$ with
a union of irreducible components of $\overline{Y}_{U'}$.\end{remark}

\begin{remark} \label{rmk:logdiffs}
We note also that the ordinary locus of $Y^{\mathrm{Ig}}_{J,N}$ can instead be viewed as parametrising
pairs $(\underline{A},\iota)$ where $\iota:\mu_p\otimes O_F \stackrel{\sim}{\to} \ker \mathrm{Frob}_A$.
Since $Y^{\mathrm{Ig}}_{J,N}$ is normal, it is essentially the scheme defined as 
$\overline{\mathfrak{M}}(E,\mu_{pN})^{\mathrm{Kum}}$ in \cite[\S9]{AG}; the differences are that
we are working with full level $N$ structure and not including the cusps.  We will not however make
any direct use of the fact that $Y^{\mathrm{Ig}}_U$ or $Y^{\mathrm{Ig}}_{J,N}$ is normal; in particular
we will not compute divisors on them as in \cite[\S12]{AG}, appealing instead in the proof of Theorem~\ref{thm:theta}
below to general properties of logarithmic differentiation in order to descend the problem to $\overline{Y}_U$.\end{remark}

\subsection{Construction of $\Theta$-operators}
For each $\tau\in \Sigma$, we consider the inclusion 
$\overline{\omega}_\tau^{-1}\subset \mathrm{Sym}_{\mathcal{O}_{\overline{Y}_U}}
(\oplus_{\tau\in\Sigma} \overline{\omega}_\tau^{-1})$, which induces an injective morphism
$$\overline{\omega}_\tau^{-1} \to 
 \pi_{U,*}\mathcal{O}_{Y^{\mathrm{Ig}}_U} 
 = \mathrm{Sym}_{\mathcal{O}_{\overline{Y}_U}}
(\oplus_{\tau\in\Sigma} \overline{\omega}_\tau^{-1})/\mathcal{I},$$
hence an injective morphism $\pi_U^*\overline{\omega}_\tau^{-1} \to \mathcal{O}_{Y^{\mathrm{Ig}}_U}$,
which we view as a section of $\pi_U^*\overline{\omega}_\tau$.
We denote this section by $h_\tau$, and call it a {\em fundamental Hasse invariant}.
The definition of $Y_U^{\mathrm{Ig}}$ implies that these satisfy the relation
$$h_{\mathrm{Fr}^{-1}\circ\tau}^p  = h_\tau  \pi_U^*(\mathrm{Ha}_\tau).$$

Recall now the Kodaira--Spencer isomorphism (\ref{eqn:KS1}).  Taking 
$A$ to be the universal HBAV over $S = \overline{Y}_{J,N}$ and decomposing
over embeddings $\tau$ yields a $G_{U,N}$-equivariant isomorphism
$$\bigoplus_{\tau\in \Sigma} \left( \overline{\omega}_\tau^2 \otimes_{\mathcal{O}_{\overline{Y}_{J,N}}} \overline{\delta}_\tau^{-1}\right)
   \simeq  \Omega^1_{\overline{Y}_{J,N}/E}$$
whose union over $J$ descends to an isomorphism
\begin{equation}\label{eqn:KS2}
\bigoplus_{\tau\in \Sigma} \left( \overline{\omega}_\tau^2 \otimes_{\mathcal{O}_{\overline{Y}_U}} \overline{\delta}_\tau^{-1}\right)
   \simeq  \Omega^1_{\overline{Y}_U/E}\end{equation}
of vector bundles on $\overline{Y}_U$.
We let $\mathrm{KS}_\tau: \Omega^1_{\overline{Y}_U/E}  \to  
\overline{\omega}_\tau^2 \otimes_{\mathcal{O}_{\overline{Y}_U}} \overline{\delta}_\tau^{-1}$
denote the composite of its inverse with the projection to the $\tau$-component.

Let $\mathcal{F}_U$ denote the sheaf of total fractions on $\overline{Y}_U$, and
$F_U = H^0(\overline{Y}_U,\mathcal{F}_U)$ the ring of meromorphic functions on
$\overline{Y}_U$, so $F_U$ is the product of the function fields of the components
of $\overline{Y}_U$.  Similarly let $\mathcal{F}_U^{\mathrm{Ig}}$ be the sheaf of
total fractions on $Y_U^{\mathrm{Ig}}$ (so $\mathcal{F}_U^{\mathrm{Ig}} = \pi_U^*\mathcal{F}_U$)
and let $F_U^{\mathrm{Ig}}$ be the ring of meromorphic functions on  $Y_U^{\mathrm{Ig}}$,
so $F_U^{\mathrm{Ig}}$ is Galois over $F_U$ with Galois group $(O_F/pO_F)^\times$.
Since the natural map $\pi_U^*\Omega^1_{\overline{Y}_U/E} \to \Omega^1_{Y_U^{\mathrm{Ig}}/E}$
is generically an isomorphism, i.e., 
$$\pi_U^*(\Omega^1_{\overline{Y}_U/E}\otimes_{\mathcal{O}_{\overline{Y}_U}} 
\mathcal{F}_U) \, \simeq \, 
\pi_U^*\Omega^1_{\overline{Y}_U/E}\otimes_{\mathcal{O}_{Y^{\mathrm{Ig}}_U}} 
\mathcal{F}_U^{\mathrm{Ig}} \, \simeq \, \Omega^1_{Y^{\mathrm{Ig}}_U/E}\otimes_{\mathcal{O}_{Y^{\mathrm{Ig}}_U}} 
\mathcal{F}_U^{\mathrm{Ig}},$$
the pull-back of $\mathrm{KS}_\tau$ induces a morphism
$$ \Omega^1_{Y^{\mathrm{Ig}}_U/E}\otimes_{\mathcal{O}_{Y^{\mathrm{Ig}}_U}}
\mathcal{F}_U^{\mathrm{Ig}}  \to \pi_U^*\left(  
\overline{\omega}_\tau^2 \otimes_{\mathcal{O}_{\overline{Y}_U}} \overline{\delta}_\tau^{-1}\right)
\otimes_{\mathcal{O}_{Y^{\mathrm{Ig}}_U}} 
\mathcal{F}_U^{\mathrm{Ig}},$$
which we will denote by $\mathrm{KS}_\tau^{\mathrm{Ig}}$.

Suppose now that $f \in M_{k,l}(U;E)$.  Let $h^k = \prod_\tau h_\tau^{k_\tau}$ and
$g^l = \prod_\tau g_\tau^{l_\tau}$, where $h_\tau$ is the fundamental Hasse invariant
and $g_\tau$ is any trivialisation of $\overline{\delta}_\tau$.  Then 
$h^{-k} \pi_U^*(g^{-l}f) \in F_{U}^{\mathrm{Ig}}$, so we may apply $\mathrm{KS}_\tau^{\mathrm{Ig}}$ to
$$d(h^{-k}\pi_U^*(g^{-l}f)) \in \Omega^1_{F_U^{\mathrm{Ig}}/E} =   H^0(Y_U^{\mathrm{Ig}}, \Omega^1_{Y^{\mathrm{Ig}}_U/E}  \otimes_{\mathcal{O}_{Y^{\mathrm{Ig}}_U}}
\mathcal{F}_U^{\mathrm{Ig}}).$$
\begin{definition}\label{def:igusatheta} We define
$$\Theta_\tau^{\mathrm{Ig}}(f) =  h^k \pi_U^*(g^l\mathrm{Ha}_\tau) \mathrm{KS}_\tau^{\mathrm{Ig}}
(d(h^{-k}\pi_U^*(g^{-l}f)))  \in H^0(Y_{U}^{\mathrm{Ig}}, \pi_U^* \overline{\mathcal{L}}^{k',l'} \otimes
_{\mathcal{O}_{Y^{\mathrm{Ig}}_U}}
\mathcal{F}_U^{\mathrm{Ig}}),$$
where
\begin{itemize}
\item if $\mathrm{Fr}\circ\tau = \tau$, then $k'_\tau= k_\tau + p +1$ and $k'_{\tau'}= k_{\tau'}$ if $\tau' \neq \tau$;
\item if $\mathrm{Fr}\circ\tau \neq \tau$, then $k'_\tau = k_\tau +1$, $k'_{\mathrm{Fr}^{-1}\circ\tau}= k_{\mathrm{Fr}^{-1}\circ\tau}+p$, and
$k_{\tau'} = k_{\tau'}$ if $\tau' \not\in \{\mathrm{Fr}^{-1}\circ\tau,\tau\}$;
\item $l_{\tau}' = l_\tau - 1$, and $l'_{\tau'}= l_{\tau'}$ if $\tau' \neq \tau$.
\end{itemize}
\end{definition}
Since the ratio of any two trivialisations of $\overline{\delta}_\tau$ is locally constant,
we see that $\Theta_\tau^{\mathrm{Ig}}(f)$ is independent of the choice of $g_\tau$.  Moreover
it is straightforward to check that $\Theta_\tau^{\mathrm{Ig}}(f)$ is invariant under the action of
$(O_F/p)^\times$, hence descends to a section of 
$\overline{\omega}^{k',l'} \otimes_{\mathcal{O}_{\overline{Y}_U}} \mathcal{F}_U$,
which we denote by $\Theta_\tau(f)$.
What is more difficult is that $\Theta_\tau(f)$ is in fact a section of $\overline{\omega}^{k',l'}$.
In fact we have the following result, essentially due to Andreatta and Goren~\cite{AG}:

\begin{theorem} \label{thm:theta}  If $f \in M_{k,l}(U;E)$ and $\tau \in \Sigma$, then
$\Theta_\tau(f)  \in M_{k',l'}(U;E)$.  Moreover
$\Theta_\tau(f)$ is divisible by $\mathrm{Ha}_\tau$ if and only if either
$f$ is divisible by $\mathrm{Ha}_\tau$ or $k_\tau$ is divisible by $p$.
\end{theorem}
\begpf  First note that the formula
$$\prod_{\tau'\in \Sigma} h_{\tau'}^{p-1} = \prod_{\tau'\in \Sigma} \pi_U^*(\mathrm{Ha}_{\tau'})$$
implies that $h^k$ is non-vanishing on $\pi_U^{-1}(Y_U^{\mathrm{ord}})$, so
$h^{-k} \pi_U^*(g^{-l}f)$ and hence $d(h^{-k} \pi_U^*(g^{-l}f))$ are regular on
$\pi_U^{-1}(Y_U^{\mathrm{ord}})$.  Since $\pi_U$ is \'etale on $\pi_U^{-1}(Y_U^{\mathrm{ord}})$,
it follows that $d(h^{-k} \pi_U^*(g^{-l}f))$ restricts to a section of
$$\pi_U^*\Omega^1_{Y_U^{\mathrm{ord}}/E}
    = \Omega^1_{\pi_U^{-1}(Y_U^{\mathrm{ord}})/E}.$$
Therefore $\Theta_\tau^{\mathrm{Ig}}(f)$ is regular on $\pi_U^{-1}(Y_U^{\mathrm{ord}})$,
and $\Theta_\tau(f)$ is regular on $Y_U^{\mathrm{ord}}$.

To complete the proof of the first assertion, we must show that if $z$ is the
generic point of an irreducible component $Z \subset Z_{U,\tau_0}$, then (the
germ at $z$) of $\Theta_\tau(f)$ lies in $\overline{\omega}^{k',l'}_z$, i.e., that 
$\mathrm{ord}_z(\Theta_\tau(f)) \ge 0$.  For the second assertion, it suffices
to show further that if $\tau_0 = \tau$, then $\mathrm{ord}_z(\Theta_\tau(f)) > 0$
if and only if either $p|k_\tau$ or $\mathrm{ord}_z(f) > 0$.

Let us revert to the notation of the proof of Proposition~\ref{prop:igusa}, so 
$\mathcal{O}_{\overline{Y}_U,z} = R_{\mathfrak{p}}$ and $x_{\tau'}$ is
a basis for the stalk of $\overline{\omega}_{\tau'}^{-1}$ at $z$.
Let $y_{\tau'}$ be the dual basis for $\overline{\omega}_{\tau'}$, so
that $y^kg^l$ is a basis for $\overline{\omega}^{k,l}_z$ over $R_{\mathfrak{p}}$
(where as usual $y^k$ denotes $\prod_{\tau'} y_{\tau'}^{k_{\tau'}}$),
and we have $f = \phi_f y^k g^l$ for some $\phi_f \in R_{\mathfrak{p}}$.
In terms of the basis $\pi_U^*(y_{\tau'})$ for $\pi_U^*\overline{\omega}_{\tau',z}$,
the fundamental Hasse invariant $h_{\tau'}$ is given by $t_{\tau'} \pi_U^*(y_{\tau'})$,
so that $h^{-k}\pi_U^*(g^{-l}f) = t^{-k}\phi_f$ in the total fraction ring of $T_{\mathfrak{p}}$,
over which we deduce that
$$\Theta_\tau^{\mathrm{Ig}}(f)
   = \mathrm{KS}_\tau^{\mathrm{Ig}}(t^k\,d(t^{-k}\phi_f)) \pi_U^*(\mathrm{Ha}_\tau y^k g^l).$$
The formulas $t_{\mathrm{Fr}^{-1}\circ\tau}^p = r_\tau t_\tau$ in $T_{\mathfrak{p}}$ imply that 
$r_{\tau'}\, dt_{\tau'} = - t_{\tau'}\, dr_{\tau'}$ in $\Omega^1_{T_{\mathfrak{p}}/E}$, and it follows that
\begin{equation} \label{eqn:Theta} \begin{array}{rcl}  \Theta_\tau(f)
  &  = &\mathrm{KS}_\tau(d\phi_f -  \phi_f r^k\,d(r^{-k}))  \mathrm{Ha}_\tau y^k g^l \\
   & = &  \mathrm{KS}_\tau\left( d\phi_f  + 
    \phi_f \displaystyle\sum_{\tau'\in \Sigma}  k_{\tau'} \frac{dr_{\tau'}}{r_{\tau'}} \right) \mathrm{Ha}_\tau y^kg^l
   \end{array}\end{equation}
(locally at $z$).  Since $\mathrm{Ha}_\tau = r_\tau y^{k_{\mathrm{Ha}_\tau}}$, we conclude that
$$\mathrm{ord}_z(\Theta_\tau(f))
 = \mathrm{ord}_z\left( r_\tau \mathrm{KS}_\tau(d\phi_f) +  k_\tau \phi_f \mathrm{KS}_\tau (dr_\tau)
          + \displaystyle\sum_{\tau'\neq \tau} k_{\tau'} r_\tau \phi_f \frac{\mathrm{KS}_\tau(dr_{\tau'})}{r_{\tau'}}
          \right).$$

In particular, if $\tau = \tau_0$, then we see immediately that $\mathrm{ord}_z(\Theta_\tau(f)) \ge 0$,
with equality if and only if  $\mathrm{ord}_z( k_\tau \phi_f \mathrm{KS}_\tau (dr_\tau) ) = 0$,
so in this case we are reduced to proving that $\mathrm{ord}_z (\mathrm{KS}_\tau (dr_\tau) ) = 0$.
On the other hand if $\tau \neq \tau_0$, then we are reduced to proving that 
$\mathrm{ord}_z( \mathrm{KS}_\tau (dr_{\tau_0}) ) > 0$.  Both cases are treated 
by the following lemma. \epf

The following is essentially the unramified case of \cite[Prop.~12.34]{AG}, which we prove
using a computation of Koblitz~\cite{Ko} (as presented in \cite{I}) instead of the theory of displays.

\begin{lemma}  \label{lemma:HW}  Let $z$ be a generic point of $Z_{U,\tau_0}$ and
let $r$ be a generator for the maximal ideal of $\mathcal{O}_{\overline{Y}_U,z}$.  Then
$\mathrm{ord}_z(  \mathrm{KS}_\tau(dr) )= 0$ if and only if $\tau = \tau_0$.
\end{lemma}
\begpf First note that since the projection $\coprod \overline{Y}_{J,N} \to \overline{Y}_U$
is \'etale, we may replace $\overline{Y}_U$ by $\overline{Y}_{J,N}$ and $Z_{U,\tau_0}$ by $Z_{\tau_0}$
in the statement of the lemma.  Note also that the conclusion of the lemma is independent of the choice
of uniformising parameter $r$, since if $u\in \mathcal{O}_{\overline{Y}_{J,N},z}^\times$, then
$$\mathrm{KS}_\tau(d(ur)) = u \mathrm{KS}_\tau(dr) + r \mathrm{KS}_\tau(du).$$
We will prove that for every closed point $x$ of $Z_{\tau_0}$, there is a choice
of parameter $r$, regular at $x$, such that
the fibre of $\mathrm{KS}_\tau(dr)$ at $x$ vanishes if and only if $\tau \neq \tau_0$.
By the formula above, the equivalence then holds for all $r$ regular at $x$, hence
for all $x$ at which any given $r$ is regular, and this implies the lemma.

Let $R = \widehat{\mathcal{O}}_{\overline{Y}_{J,N},x}$, $A$ the pull-back to $\mathrm{Spec}\,R$
of the universal HBAV over $Y_{J,N}$ and $M = H^1_{\mathrm{DR}}(A/R)$.  Letting
$$\nabla: M \to \Omega^1_{R/E} \otimes_R M$$
denote the Gauss--Manin connection and $\phi:M \to M$ the morphism induced by
the absolute Frobenius morphism on $A$, we are in the situation of 
\cite[A2.1]{I}.\footnote{The notations $A$ and $F$ in \cite{I}, being in other use here, have been
replaced by $R$ and $\phi$.}   We then have $L = H^0(A,\Omega^1_{A/R})$ and
$N = H^1(A,\mathcal{O}_A)$, and the $O_F$ action on $A$ yields decompositions
$L = \oplus L_\tau$, $M = \oplus M_\tau$ and $N = \oplus N_\tau$ into free
$R$-modules indexed by $\tau\in \Sigma$.  The Hasse--Witt endomorphism
of \cite[(A2.1.1)]{I} decomposes as the sum of Frobenius semi-linear morphisms
$N_{\mathrm{Fr}^{-1}\circ\tau} \to N_\tau$ such that, after choosing a generator
for $J/pJ$ and applying the isomorphisms $\mathrm{Lie}(A)\otimes_{O_F}J\simeq N$, the induced
morphism $N_{\mathrm{Fr}^{-1}\circ\tau}^{\otimes p} \to N_\tau$ corresponds to
the completion at $x$ of the partial Hasse invariant $\mathrm{Ha}_{N,J,\tau}$.
The general properties of the construction of the Gauss--Manin connection ensure
its compatibility with the $O_F$-action, so that it decomposes as a direct sum
of connections $\nabla_\tau$ on $M_\tau$.  Furthermore the morphism \cite[(A2.1.2)]{I}
induced by $\nabla$ is the completion at $x$ of the reduction mod $\pi$ of (\ref{KS0}),
and hence the induced morphism
$$\bigoplus_\tau \mathrm{Hom}_R (\wedge^2_R M_\tau, L_\tau^{\otimes_R 2}) = 
\bigoplus_\tau \mathrm{Hom}_R (N_\tau,L_\tau) \to \Omega^1_{R/E}$$
is the completion at $x$ of the Kodaira--Spencer isomorphism on $\overline{Y}_{J,N}$.

Following \cite{I}, we let $R_i = R/\mathfrak{m}_R^{i+1}$ (the cases of interest being
$i=0,1$), and similarly use subscript $i$ for reductions mod $\mathfrak{m}^{i+1}$ of $R$-modules,
morphisms and matrices.  We now choose a basis for $M_1$ as in \cite[(A2.1.6)]{I} as follows:
First choose a basis for $M_0$ consisting of vectors $e_{\tau,0} \in L_{\tau,0}$, $f_{\tau,0} \in M_{\tau,0}$
for $\tau \in \Sigma$.  Then 
$$\phi_0(e_{\mathrm{Fr}^{-1}\circ \tau,0}) = 0 \qquad \mbox{and}\qquad  \phi_0(f_{\mathrm{Fr}^{-1}\circ\tau,0}) =  
  c_\tau e_{\tau,0} + d_\tau f_{\tau,0}$$
for some $c_\tau,d_\tau \in R_0$ not both zero.   Replacing $f_{\tau,0}$ by $e_{\tau,0}+f_{\tau,0}$ whenever
$c_\tau = 0$, we may assume $c_\tau \neq 0$, and then replacing $e_{\tau,0}$ by $c_\tau^{-1}e_{\tau,0}$,
we may assume $c_\tau = 1$ for all $\tau$.   Now lift each pair $(e_{\tau,0},f_{\tau,0})$ to a basis
$(e_\tau,f_\tau')$ of $M_{\tau,1}$
with $e_\tau \in L_{\tau,1}$ and let $f_\tau = P(f_\tau')$ where $P$ is defined in \cite[(A2.1.3)]{I}. 
Since $\nabla$ respects the decomposition $M = \oplus M_\tau$, so does $P$, and hence
$f_\tau \in M_{\tau,1}$.  Moreover $f_\tau \equiv f_\tau' \bmod \mathfrak{m}_R$, so in the matrix
$\left( \begin{array}{cc} 0 & B_1 \\ 0 & H_1 \end{array}\right)$ 
of \cite[(A2.1.7)]{I} representing $\phi$ on $M_1$ with respect to the basis $((e_\tau),(f_\tau))$:
\begin{itemize}
\item the reduction $B_0$ of $B_1$ is defined by $b_{\tau,\tau'} = \delta_{\tau,\mathrm{Fr}\circ\tau'}$;
\item the matrix $H_1$ represents the Hasse--Witt  endomorphism $N_1 \to N_1$ with respect to the basis induced by
$(f_\tau)$, so $h_{\tau,\tau'} = 0$ if $\tau' \neq \mathrm{Fr}^{-1}\circ \tau$ and $h_{\tau,\mathrm{Fr}^{-1}\circ\tau} = r_\tau$,
where $r_\tau$ represents the pull-back of $\mathrm{Ha}_{J,N,\tau}$ to $R_1$ with respect to the basis induced by the
map sending $f_{\mathrm{Fr}^{-1}\circ\tau}$ to $f_\tau$.
\end{itemize}
In particular, $B_0$ is invertible, and Proposition~A2.1.8 of \cite{I} gives that the matrix
$$K_0 = (H_0 - H_1)B_0^{-1}$$
with entries in $\mathfrak{m}_R/\mathfrak{m}_R^2 \cong \Omega^1_{R/E} \otimes_R R_0$
is diagonal with $(\tau,\tau)$-entry $-dr_\tau$.  Note that
the map $L_0 \to (\mathfrak{m}_R / \mathfrak{m}_R^2) \otimes_{R_0}  N_0$ is
the fibre of (\ref{KS0}) at $x$, and is represented by $K_0$ with respect to the bases
$(e_{\tau,0})$ of $L_0$ and $(f_{\tau,0}\bmod L_0)$ of $N_0$.    It follows that the fibre at $x$
of the Kodaira--Spencer isomorphism is the map
$$ \bigoplus_\tau \mathrm{Hom}_{R_0} ( N_{\tau,0}, L_{\tau,0})
  \simeq \mathrm{Hom}_{R_0\otimes O_F} (N_0, L_0)
   \longrightarrow \Omega^1_{R/E} \otimes_R R_0$$
 under which the basis vector in the $\tau$-component induced by $f_{\tau,0} \mapsto e_{\tau,0}$
 corresponds to $-dr_{\tau}$.   Note that $r_{\tau_0}$ is the image in $\mathfrak{m}_R/\mathfrak{m}_R^2$
 of a uniformising parameter for $Z_{\tau_0}$ in a neighbourhood of $x$, and the fibre at $x$
 of $\mathrm{KS}_\tau$ sends $dr_{\tau_0}$ to $0$ if and only if $\tau \neq \tau_0$, so this
 completes the proof of the lemma.
\epf

It is straightforward to check that the maps $\Theta_\tau$ are compatible with the maps
$[U_1 g U_2]$ for sufficiently small $U_1,U_2$ and $g \in \mathrm{GL}_2(\mathbb{A}_F^\infty)$ 
such that $g^{-1}U_1 g \subset U_2$ and $g_p \in \mathrm{GL}_2(O_{F,p})$.
Taking limits over open compact subgroups $U$ therefore gives:
\begin{corollary}  \label{cor:theta} For any weight $(k,l)$, $\Theta_\tau$ defines a map:
$$M_{k,l}(E) \to M_{k',l'}(E)$$
commuting with the action of all
$g \in \mathrm{GL}_2(\mathbb{A}_F^\infty)$ such that $g_p \in \mathrm{GL}_2(O_{F,p})$.
 In particular, for any open compact subgroup $U$ of $\mathrm{GL}_2(\mathbb{A}_F^\infty)$ containing
$\mathrm{GL}_2(O_{F,p})$, $\Theta_\tau$ defines a map
$$M_{k,l}(U;E) \to M_{k',l'}(U;E)$$
commuting with the operators $T_v$ and $S_v$ for all $v \nmid p$ such that $\mathrm{GL}_2(O_{F,v}) \subset U$.
\end{corollary}

\section{$q$-expansions} \label{sec:qexpansions}

We review the definition and properties of $q$-expansions, including the effect
on them of Hecke and partial $\Theta$-operators, and we generalise a result of Katz
on the kernel of $\Theta$.  See~\cite{PRtheta} for a further generalisation to the case
where $p$ is ramified in $F$.

\subsection{Definition and explicit descriptions}\label{subsec:definition}
Suppose as usual that $U$ is a sufficiently small open compact subgroup of
$\mathrm{GL}_2(\widehat{O}_F)$ containing $\mathrm{GL}_2(\widehat{O}_{F,p})$, 
with $k,l \in \mathbb{Z}^\Sigma$ and $R$ a Noetherian $\mathcal{O}$-algebra such that
$\nu^{k+2l} = 1$ in $R$ for all $\nu \in O_F^\times \cap U$.  Recall 
from \S\ref{subsec:XU} that $X_U$ is the minimal compactification of $Y_U$,
and a cusp of $X_U$ is a connected component of $X_U - Y_U$.
When it is clear and convenient, we will suppress the subscript $R$ from base-changes of morphisms
to which schemes have been extended, writing just $j$ for example for the inclusion 
$Y_{U,R} \hookrightarrow X_{U,R}$.

\begin{definition}  \label{def:qexpansion} For each cusp $C$
of $X_U$, we let $Q_{C,R}^{k,l}$ denote the completion of $j_*\mathcal{L}_{U,R}^{k,l}$
at $C_R$, and for $f \in M_{k,l}(U;R)$, we define the {\em $q$-expansion} of $f$
at $C$ to be its image in $Q_{C,R}^{k,l}$.  
\end{definition}

We now proceed to describe $Q_{C,R}^{k,l}$ more explicitly.   (See Sections~7.1--2 of~\cite{PRtheta}
for a more general and detailed analysis.)   We first recall (e.g. from \cite{C2}) the description
in the context of $X_{J,N}$, supposing that $N \ge 3$ and $\mu_N(\overline{\mathbb{Q}}) \subset \mathcal{O}$.
The cusps $\tilde{C}$ of $X_{J,N}$ are in bijection with equivalence classes of data:
\begin{itemize}
\item fractional ideals $\mathfrak{a}$, $\mathfrak{b}$ of $O_F$;
\item an exact sequence of $O_F$-modules $0 \to (\mathfrak{ad})^{-1} \to H \to \mathfrak{b} \to 0$;
\item an isomorphism $J \stackrel{\sim}{\longrightarrow}\mathfrak{ab}^{-1}$;
\item an isomorphism $(O_F/NO_F)^2 \stackrel{\sim}{\longrightarrow} H/NH$.
\end{itemize}
By the Koecher Principle (as in the discussion preceding~\cite[Prop.~6.11]{R}; see~\cite[Thm.~2.5]{L} for
a statement in the required generality), we have $\pi_*\mathcal{L}_{J,N,R}^{k,l,\mathrm{tor}} = \tilde{j}_*\mathcal{L}_{J,N,R}^{k,l}$,
where $\pi: X_{J,N}^{\mathrm{tor}} \to  X_{J,N}$ is the projection from a toroidal compactification of $Y_{J,N}$, and
$\mathcal{L}_{J,N}^{k,l,\mathrm{tor}}$ is the canonical extension of $\mathcal{L}_{J,N}^{k,l}$ to a line bundle over 
$X_{J,N}^{\mathrm{tor}}$ associated to the semi-abelian scheme $A^{\mathrm{tor}}$ extending the universal abelian scheme $A$ over $Y_{J,N}$.
Since $\pi$ is projective,  the Formal Functions Theorem identifies the completion $Q_{\tilde{C},R}^{k,l}$ of the sheaf
$\tilde{j}_*\mathcal{L}_{J,N,R}^{k,l}$ at $\tilde{C}_R$ with the set of global sections of the completion of the coherent sheaf
$\mathcal{L}_{J,N,R}^{k,l,\mathrm{tor}}$ at the fibre of $\pi$ over $\tilde{C}_R$.  

The construction of the toroidal compactification identifies the completion of $X_{J,N}^{\mathrm{tor}}$ at the
fibre over $\tilde{C}$ with the quotient $\widehat{S}/U_N$, where $\widehat{S}$ is the formal scheme denoted
$S_N(\{\sigma_\alpha^{\tilde{C}}\})^\wedge$ in \cite[3.4.2]{C2} (for a choice of polyhedral cone decomposition
$\{\sigma_\alpha^{\tilde{C}}\}$), but we define the action of  $\alpha \in U_N = \ker(O_F^\times \to (O_F/N)^\times)$
as being induced by multiplication by $\alpha^2$.  We therefore obtain an identification of the $R$-algebra 
$Q_{\tilde{C},R}^{0,0}$ with
$$\widehat{S}_{\tilde{C},R} = (R[[q^m]]_{m \in (N^{-1}\mathfrak{ab})_+\cup \{0\}})^{U_N}$$
where $\alpha \in U_N$ acts via $q^m \mapsto q^{\alpha^2m}$ on power series, under which
$\tilde{C}$ corresponds to the closed subscheme of $\mathrm{Spec}\,\widehat{S}_{\tilde{C}}$
defined by $\sum r_mq^m \mapsto r_0$.  
Furthermore, the pull-back of $A^{\mathrm{tor}}$ to $\widehat{S}$ is the Tate semi-HBAV $T_{\mathfrak{a},\mathfrak{b}}$
associated to the  quotient $(\mathbb{G}_m  \otimes (\mathfrak{ad})^{-1})/q^{\mathfrak{b}}$
(by Mumford's construction).  The canonical trivialisations
$$\begin{array}{cccccc}
&\mathcal{L}ie(T_{\mathfrak{a},\mathfrak{b}}/\widehat{S}) &\cong& (\mathfrak{ad})^{-1} \otimes \mathcal{O}_{\widehat{S}}
&\cong& \mathrm{Hom}(\mathfrak{a},\mathcal{O}_{\widehat{S}})\\
 \mbox{and}&
 \mathcal{L}ie(T_{\mathfrak{a},\mathfrak{b}}^\vee/\widehat{S})& \cong& (\mathfrak{bd})^{-1} \otimes \mathcal{O}_{\widehat{S}}
&\cong& \mathrm{Hom}(\mathfrak{b},\mathcal{O}_{\widehat{S}})\end{array}$$
then give a trivialisation of the pull-back of $\mathcal{L}_{J,N,R}^{k,l,\mathrm{tor}}$ to $\widehat{S}$, and hence an
identification
$$Q_{\tilde{C},R}^{k,l} = (D^{k,l} \otimes_\mathcal{O} R[[q^m]]_{m \in (N^{-1}\mathfrak{ab})_+\cup \{0\}})^{U_N},$$
where 
$$D^{k,l} = \bigotimes_\tau \left((\mathfrak{a}\otimes \mathcal{O})_\tau^{\otimes k_\tau}
           \otimes_{\mathcal{O}}(J\mathfrak{d}^{-1} \otimes \mathcal{O})_\tau^{\otimes l_\tau}\right)$$
and $\alpha \in U_N$ acts as $\alpha^k$ on $D^{k,l}$.   Note that $D^{k,l}$ is a free of rank one over $\mathcal{O}$,
and letting $b$ be a basis, we have
$$Q_{\tilde{C},R}^{k,l} = \left\{ \left.\, \sum_{m\in (N^{-1}\mathfrak{ab})_+\cup \{0\}}(b \otimes  r_m)q^m \,\right|\,
   r_{\alpha^2 m} = \alpha^k r_m \ \mbox{for all $\alpha \in U_N$}\,\right\}.$$
We will also write $\overline{D}^{k,l}$ for $D^{k,l}\otimes_{\mathcal{O}}E$.

If we fix the data of $\mathfrak{a}$, $\mathfrak{b}$ and 
$\mathfrak{ab}^{-1} \simeq J$, then the corresponding cusps of $X_{J,N}$ are in
bijection with $P_N \backslash \mathrm{GL}_2(O_F/N)$, where
$$P_N = \left\{\, \left. \left(\begin{array}{cc} \alpha^{-1} & * \\ 0 & \alpha \end{array}\right) \bmod N
 \,\right|\, \alpha\in O_F^\times \,\right\}.$$
Here we have chosen isomorphisms $s:O_F/NO_F \simeq N^{-1}\mathfrak{b}/\mathfrak{b}$
and $t:O_F/NO_F \simeq\mu_N  \otimes (\mathfrak{ad})^{-1} $ 
to define a level $N$-structure $\eta$
on $T_{\mathfrak{a},\mathfrak{b}}$ by $\eta(x,y) = t(y)q^{s(x)}$, and then associated the coset $P_N g$ to the cusp 
of $X_{J,N}$ defined by $T_{\mathfrak{a},\mathfrak{b}}$ with level 
$N$-structure $\eta \circ r_{g^{-1}}$.
Under this bijection, the (right) action of $G_{U,N}$ is defined by
$$P_N g \cdot (\nu,u)  =  P_N \left(\begin{array}{cc} \nu^{-1} & 0 \\ 0 & 1 \end{array}\right) gu.$$
The stabiliser in $G_{U,N}$ of (the cusp corresponding to) $P_Ng$ is therefore the image
of the group
$$\left\{\,(\nu,u)  \in  O_{F,+}^\times \times U\,\left|\, gug^{-1} \equiv
\left(\begin{array}{cc} \nu\alpha^{-1} & * \\ 0 & \alpha \end{array}\right)  \bmod N
\ \mbox{for some $\alpha \in O_F^\times$}\,\right.\right\}.
$$
We find that the (left) action on $Q_{\tilde{C},R}^{k,l}$ of such an element $(\nu,u)$, with
$$gug^{-1} = \left(\begin{array}{cc} \nu\alpha^{-1} & -\nu\alpha^{-1}x \\ 0 & \alpha \end{array}\right)  \bmod N,$$
is defined by $\alpha^k\nu^l$ on $D^{k,l}$ and
$$\sum r_mq^m \mapsto \sum \zeta(xm) r_m q^{\alpha^2\nu^{-1}m}$$
on $R[[q^m]]_{m \in (N^{-1}\mathfrak{ab})_+\cup \{0\}}$, where $\zeta: N^{-1}\mathfrak{ab}/\mathfrak{ab} \to \mu_N$
is the composite of the $O_F$-linear isomorphism $N^{-1}\mathfrak{ab}/\mathfrak{ab}  \to \mathfrak{d}^{-1} \otimes \mu_N$
induced by $t\circ s^{-1}$ with $\mathrm{tr}_{F/\mathbb{Q}} \otimes 1$.
The module $Q_{C,R}^{k,l}$ (over $\widehat{S}_{C,R} = Q_{C,R}^{0,0}$) is then given by the invariants 
in $Q_{\tilde{C},R}^{k,l}$ under the action of the above stabiliser.  
In particular, we note the following two special cases:
\begin{proposition}  \label{prop:QCkl}  Suppose that $\mu^{k+2l} = 1$ in $R$ for all $\mu \in O_F^\times$ such that
$\mu \equiv 1 \bmod \mathfrak{n}$,
and let $b$ be a generator of $D^{k,l}$.
\begin{itemize}
\item If $U = U(\mathfrak{n})$, then
$$Q_{C,R}^{k,l} \simeq \left\{ \left.\, \sum_{m\in (\mathfrak{n}^{-1}\mathfrak{ab})_+\cup \{0\}}(b \otimes  r_m)q^m \,\right|\,
   r_{\nu m} = \nu^{-l} r_m \ \mbox{for all $\nu \in U_{\mathfrak{n},+}$}\,\right\}.$$
\item If $U = U_1(\mathfrak{n})$ and $g=1$, then
$$Q_{C,R}^{k,l} = \left\{ \left.\, \sum_{m\in (\mathfrak{ab})_+\cup \{0\}}(b \otimes  r_m)q^m \,\right|\,
   r_{\nu m} = \nu^{-l} r_m \ \mbox{for all $\nu \in O_{F,+}^\times$}\,\right\}.$$
\end{itemize}
\end{proposition}
Note that the isomorphism (in the case of $U = U(\mathfrak{n})$) depends on the choice of
representative $g$. Note also that the description of $Q_{C,R}^{k,l}$ is compatible in the obvious senses
with the morphisms induced by base-changes $R \to R'$, and inclusions $U' \subset U$ (for 
cusps $C'$ of $X_{U'}$ mapping to $C$).

\subsection{The $q$-expansion Principle}

The $q$-expansion at $C$ of a form $f \in M_{k,l}(U;R)$ vanishes if and only (the extension to $X_{U,R}$ of)
$f$ vanishes on a neighbourhood of $C_R$, which is equivalent to the vanishing of $f$ on all connected
components of $X_{U,R}$ intersecting $C_R$.   (Note that if $\mathrm{Spec}\,R$ is connected, then there is a unique
component containing $C_R$.)
Recall from \S\ref{subsec:components} that $Z_U$ is the scheme representing the set of components of $Y_U$ and hence $X_U$,
so we have the following:
\begin{lemma} \label{lem:qexp}
If $\mathcal{S}$ is any set of cusps of $X_U$ such that $\coprod_{C\in \mathcal{S}} C \to Z_U$ is surjective, then the
$q$-expansion map:
$$ M_{k,l}(U;R)  \longrightarrow \bigoplus_{C \in \mathcal{S}} Q_{C,R}^{k,l} $$
is injective.  
\end{lemma}

If $U = U_1(\mathfrak{n})$, then $\det U = \widehat{O}_F^\times$, so $Z_U$ is in bijection
with the strict class group of $F$.  For each representative $J$, we choose $\mathfrak{b} = J^{-1}$, $\mathfrak{a} = O_F$,
and $\mathcal{S}$ consisting of  a single $C$ {\em at infinity} (i.e. as in the second part of Proposition~\ref{prop:QCkl}) on each
component associated to a fixed $t:O_F/NO_F \simeq  \mu_N \otimes\mathfrak{d}^{-1}$ (independent of $J$)
and $s:O_F/NO_F \simeq (NJ)^{-1}/J^{-1}$ (of which $C$ is independent).
Using the isomorphism $D^{k,l} \cong \mathcal{O}$ obtained from the inclusion
$J\mathfrak{d}^{-1} \subset F$ for each $J$, we obtain an injective $q$-expansion map
(defined for arbitrary $\mathfrak{n}$):
\begin{equation}\label{eqn:Qexp}  M_{k,l}(U_1(\mathfrak{n});R)  \to \bigoplus_J  
\left\{ \left.\, \sum_{m\in J^{-1}_+\cup \{0\}}r_mq^m \,\right|\,
   r_{\nu m} = \nu^{-l} r_m \ \mbox{for all $\nu \in O_{F,+}^\times$}\,\right\}.\end{equation}

\subsection{$q$-expansions of partial Hasse invariants}
Let us now return to the case of arbitrary (sufficiently small) $U$, take $R=E$ and consider
the consider the $q$-expansions of the partial Hasse invariants $\mathrm{Ha}_\tau$.
Since the pull-back
$\mathrm{Ver}^*_{T_{\mathfrak{a},\mathfrak{b}}}$ of the relative Verschiebuung on the Tate semi-HBAV $T_{\mathfrak{a}, \mathfrak{b}}$
(see \S\ref{subsec:definition}) is induced by the canonical isomorphism
$$\bigoplus_{\tau} \,(\mathfrak{a}\otimes E)_\tau  \to 
 \bigoplus_{\tau}\, (\mathfrak{a}\otimes E)^{\otimes p}_{\mathrm{Fr}^{-1}\circ\tau},$$
we see that the $q$-expansion of $\mathrm{Ha}_\tau$ at any cusp is the constant $1$, or more precisely
$\iota_\tau \otimes 1$  where $\iota_\tau \in \overline{D}^{k_{\mathrm{Ha}_\tau},0}$ is defined by
the canonical isomorphism
$(\mathfrak{a}\otimes E)_\tau  \to (\mathfrak{a}\otimes E)^{\otimes p}_{\mathrm{Fr}^{-1}\circ\tau}$.

Furthermore the same is true for the $q$-expansions of $G_\tau \in M_{0,k_{\mathrm{Ha}_\tau}}(U;E)$
(with $\overline{D}^{k_{\mathrm{Ha}_\tau},0}$ replaced by  $\overline{D}^{0,k_{\mathrm{Ha}_\tau}}$
and $\mathfrak{a}$ by $J\mathfrak{d}^{-1}$), where $G_\tau$ is obtained by descent from
the canonical trivialisations defined in \S\ref{subsec:cantriv}.

\subsection{$\Theta$-operators on $q$-expansions} \label{subsec:qH}
We continue to assume $R=E$.  We now describe the effect of $\Theta$-operators on $q$-expansions.
We first assume $U = U(\mathfrak{n})$ for some $\mathfrak{n}$ sufficiently small that $\nu^l = 1 \bmod p$
for all $\nu \in U_{\mathfrak{n},+}$.   By Proposition~\ref{prop:QCkl}, we can identify $\overline{Q}_C^{k,l}
= Q_{C,E}^{k,l}$ with
$$ \overline{D}^{k,l} \otimes_E \widehat{S}_{C,E} = \left\{ \left.\, \sum_{m\in (\mathfrak{n}^{-1}\mathfrak{ab})_+\cup \{0\}}
(\overline{b} \otimes  r_m)q^m \,\right|\,
   r_{\nu m} = r_m \ \mbox{for all $\nu \in U_{\mathfrak{n},+}$}\,\right\},$$
where $\overline{b}$ is any basis for $\overline{D}^{k,l}$.
In particular note that $\overline{Q}_{C}^{k,l}$ is free 
over $\widehat{S}_{C,E}$ for all $k,l$.

We now appeal to the formula (\ref{eqn:Theta}), and observe that it is compatible
with the analogous formula defining a map $\overline{Q}_{C}^{k,l} \to \overline{Q}_{C}^{k',l'}$,
where $\mathrm{KS}_\tau$ is replaced by the completion of $j_*\mathrm{KS}_\tau$ at $\overline{C} =C_E$.
Moreover the formula is valid for any choices of bases $y_{\tau'}$
for the completions of $j_*\overline{\omega}_{\tau'}$ (which are invertible thanks to
our choice of $U$).   In particular we can choose the $y_{\tau'}$ of the form
$a_{\tau'} \otimes 1$ where the $a_{\tau'}$ are bases for $(\mathfrak{a}\otimes E)_{\tau'}$
such that $\iota_{\tau'} \otimes a_{\tau'} = a_{\mathrm{Fr}^{-1}\circ \tau'}^{\otimes p}$ for all $\tau'$.
This gives $y^kg^l = \overline{b}\otimes 1$ for some basis $\overline{b}$ of $\overline{D}^{k,l}$, and
in view of the $q$-expansions of the partial Hasse invariants, $r_{\tau'} = 1$ for all $\tau'$.
Thus if $f$ has $q$-expansion $\sum (\overline{b}\otimes r_m)q^m$ at $C$, then we are reduced to computing
the image of $\phi_f = \sum r_mq^m$ under the composite 
\begin{equation}
\label{eqn:KShat} \widehat{S}_{C,E}  \stackrel{d}{\to} (\Omega^1_{\overline{X}_U/E})_{\overline{C}}^{\wedge}  \to 
        (j_*\Omega^1_{\overline{Y}_U/E})_{\overline{C}}^{\wedge}   \to  \overline{Q}_C^{(2,-1)_\tau},\end{equation}
where $(2,-1)_\tau = (k',l') - (k,l) - (k_{\mathrm{Ha}_\tau},0)$ and the last map
is induced by $j_*\mathrm{KS}_\tau$.

A computation on the toroidal compactification identifies $(j_*\Omega^1_{\overline{Y}_U/E})_{\overline{C}}^{\wedge}$
with
$$\mathfrak{n}^{-1}\mathfrak{ab} \otimes \widehat{S}_{C,E}$$
(in view of our assumption that $U = U(\mathfrak{n})$ for sufficiently small $\mathfrak{n}$) and
the composite of the first two maps of (\ref{eqn:KShat}) with $q^m \mapsto m \otimes q^m$.
Moreover identifying $\overline{Q}_C^{(2,-1)_\tau}$ with
$$\overline{D}^{(2,-1)_\tau} \otimes_{E}  \widehat{S}_{C,E}  =
 (\mathfrak{dab} \otimes_{\mathcal{O}}E )_\tau \otimes_{E}  \widehat{S}_{C,E},$$
\cite[(1.1.20)]{K} gives that the last map of (\ref{eqn:KShat}) is
the inverse of the isomorphism induced by the inclusion
$\mathfrak{dab} \to \mathfrak{n}^{-1}\mathfrak{ab}$, followed by projection to the
$\tau$-component.  Therefore the $q$-expansion of $\Theta_\tau(f)$ at $C$ is given by
$$\sum  (\iota_\tau \tau(m)\overline{b}\otimes r_m)q^m.$$
In view of the compatibility of $q$-expansions with the morphisms
induced by inclusions $U' \subset U$, this formula is in fact valid for all sufficiently
small $U$.  We have thus proved:

\begin{proposition} \label{prop:Thetaq} If the $q$-expansion at $C$ of $f \in M_{k,l}(U;E)$ is
$\sum (\overline{b}\otimes r_m)q^m$, then the $q$-expansion at $C$ of $\Theta_\tau(f)$ is
$\sum  (\iota_\tau \tau(m)\overline{b}  \otimes r_m)q^m$.
\end{proposition}

Recall from Lemma~\ref{lem:qexp} that a form is determined by its $q$-expansions.
Using also that $\iota_\tau\tau(m) = (\mathrm{Fr}^{-1}\circ\tau(m))^p$, we deduce
the following (see \S\ref{subsec:qH} for the definition of $G_\tau$):

\begin{corollary}  \label{cor:Thetarelns}  For all $\tau,\tau' \in \Sigma$ and $f \in M_{k,l}(U;E)$, we have the relations
\begin{itemize}
\item $\Theta_\tau\Theta_{\tau'}(f) = \Theta_{\tau'}\Theta_\tau(f)$, and
\item  $G_\tau\Theta_{\mathrm{Fr}^{-1}\circ\tau}^p(f) = \mathrm{Ha}_{\mathrm{Fr}^{-1}\circ\tau}^p \mathrm{Ha}_\tau \Theta_\tau(f)$.
\end{itemize}
\end{corollary}

\subsection{Hecke operators on $q$-expansions}
We now describe the effect of the Hecke operators $T_v$ on $q$-expansions in the case
of $U = U_1(\mathfrak{n})$.  For $f \in M_{k,l}(U_1(\mathfrak{n});R)$ and $m \in J_+^{-1} \cup \{0\}$,
we write $r_m^J(f)$ for the coefficient of $q^m$ in the $J$-component of its $q$-expansion as
in (\ref{eqn:Qexp}).

\begin{proposition}  \label{prop:Tvonqexps}
If $f \in M_{k,l}(U_1(\mathfrak{n});R)$, $v\nmid \mathfrak{n}p$ and $m \in J_+^{-1} \cup \{0\}$, then
$$r_m^J(T_vf) = \beta_1^l  r_{\beta_1 m}^{J_1}(f)  + \mathrm{Nm}_{F/\mathbb{Q}}(v) \beta_2^l r_{\beta_2 m}^{J_2}S_vf),$$
where the $J_i$ and $\beta_i \in F_+$ are such that $vJ = \beta_1J_1$ and $v^{-1}J = \beta_2J_2$
(and we interpret $r_{\beta_2 m}$ as $0$ if $\beta_2m \not\in J_2^{-1}$).
\end{proposition}
\begpf This is a standard computation which we briefly indicate how to 
carry out in our context.  Let $U = U_1(\mathfrak{n})$,
$g = \left(\begin{array}{cc}1&0\\0&\varpi_v\end{array}\right)$, denote the rational prime in $\varpi_v$ by $r$, and
choose a sufficiently large $N$ prime to $pr$.  We may extend scalars so as to assume
$\mu_{Nr}(\overline{\mathbb{Q}}) \subset \mathcal{O}$.

Note that we have 
$$UgU =U \left(\begin{array}{cc}\varpi_v&0\\0&1\end{array}\right) U= \coprod_{i \in \mathbb{P}^1(O_F/\varpi_v)} g_iU,$$
with $g_i \in \mathrm{GL}_2(O_{F,\varpi_v})$ defined by
$\left(\begin{array}{cc}\varpi_v & [i] \\0&1\end{array}\right)$ if $i \in O_F/\varpi_v$, where $[i]$ is the Teichmuller
(or indeed any) lift of $i$, and $g_\infty = \left(\begin{array}{cc}0&1\\ \varpi_v&0\end{array}\right)$.
To define the maps $[U'g_iU]:M_{k,l}(U;R) \to M_{k,l}(U';R)$ (where $U' = U(rN)$ for example), we may
take $N_1 = rN$, $N_2 = N$ and $\alpha = 1$ in the notation of \S\ref{section:Hecke}.

Recall that the $J$-component of the $q$-expansion of $T_vf$ is given by its image in $Q_{C,R}^{k,l}$
where the cusp $C$ of $X_U$ is the image of a cusp $\tilde{C}$ of $X_{J,rN}$ associated to the Tate semi-HBAV
$T_{\mathfrak{a},\mathfrak{b}}$ with $\mathfrak{a} = O_F$, $\mathfrak{b} = J^{-1}$, canonical polarisation data
(i.e., associated to the identity $\mathfrak{ab}^{-1} = J$), and level structure defined by $\eta(x,y) = t(y)q^{s(x)}$
for some choice of isomorphisms $s:O_F/rNO_F \simeq (rNJ)^{-1}/J^{-1}$ and $t:O_F/rNO_F \simeq
\mu_{rN} \otimes \mathfrak{d}^{-1}$. 

Suppose first that $i \in O_F/\varpi_v$.  Choosing $\beta = \beta_1$ in the definition of 
$\tilde{\rho}_{g_i}: Y_{J,rN}  \to Y_{J_1,N}$ and extending to minimal compactifications,
we find that $\tilde{\rho}_{g_i}(\tilde{C}) = \tilde{C}_1$ where $\tilde{C}_1$ is the cusp of
$X_{J_1,N}$ associated to $T_{O_F,J_1^{-1}}$ with canonical polarisation data
and level $N$ structure defined by $(x,y)\mapsto t(ry)q^{\beta_1rs(x)}$.  
Moreover, the induced morphism $\widehat{S}_{\tilde{C}_1} \to \widehat{S}_{\tilde{C}}$
on completions is defined by $q^m \mapsto \zeta_i(\beta_1^{-1}m)q^{\beta_1^{-1}m}$,
with $\zeta_i$ running through the distinct homomorphisms $(vNJ)^{-1}/(NJ)^{-1} \to \mu_r$ as
$i$ runs through the distinct elements of $O_F/\varpi_v$, and the pull-back to $\widehat{S}$
of the isogeny denoted $\pi$ in \S\ref{section:Hecke} is just the natural projection
$T_{O_F,J^{-1}}  \to \tilde{\rho}_{g_i}^* T_{O_F,J_1^{-1}}$ induced by the identity
on $\mathbb{G}_m \otimes \mathfrak{d}^{-1}$.
Taking into account  the normalisation by $|| \det g_i || = (\mathrm{Nm}_{F/\mathbb{Q}} v)^{-1}$,
we conclude that $[U'g_i U]$ is compatible with the morphism $Q_{\tilde{C}_1,R}^{k,l} \to Q_{\tilde{C},R}^{k,l}$
on $q$-expansions defined by
$$ \sum_{m\in (NJ_1)_+^{-1}\cup\{0\}}   (b \otimes r_m) q^m\,\,\, \mapsto \,\,\, (\mathrm{Nm}_{F/\mathbb{Q}} v)^{-1}  
\!\!\!\!\!\!\!\!\!\!\sum_{m\in (vNJ)_+^{-1}\cup\{0\}}
( \beta_1^l  b \otimes  \zeta_i(m)  r_{\beta_1m} ) q^m.$$

As for  $i=\infty$, note that $[U' g_\infty U] = [U'hU] S_v$, where 
$h = \left(\begin{array}{cc} 0 & \varpi_v^{-1} \\ 1 &0\end{array}\right)$.
Choosing $\alpha = r$ and $\beta = r^2\beta_2$ in the definition of 
$\tilde{\rho}_{h}: Y_{J,rN}  \to Y_{J_2,N}$ and extending to minimal compactifications,
we find that $\tilde{\rho}_{h}(\tilde{C}) = \tilde{C}_2$ where $\tilde{C}_2$ is the cusp of
$X_{J_2,N}$ associated to $T_{O_F,J_2^{-1}}$ with canonical polarisation data
and level $N$ structure defined by $(x,y)\mapsto t(ry)q^{\beta_2rs(x)}$.  
Moreover, the induced morphism $\widehat{S}_{\tilde{C}_2} \to \widehat{S}_{\tilde{C}}$
on completions is defined by $q^m \mapsto q^{\beta_2^{-1}m}$
and the pull-back of $\pi$ to $\widehat{S}$ is the map
$T_{O_F,J^{-1}}  \to \tilde{\rho}_{h}^* T_{O_F,J_2^{-1}}$ induced by multiplication by $r$
on $\mathbb{G}_m \otimes \mathfrak{d}^{-1}$.
Taking into account  the normalisation by $|| \det h || = \mathrm{Nm}_{F/\mathbb{Q}} v$,
we conclude that $[U'hU]$ is compatible with the morphism $Q_{\tilde{C}_2,R}^{k,l} \to Q_{\tilde{C},R}^{k,l}$
on $q$-expansions defined by
$$ \sum_{m\in (NJ_2)_+^{-1}\cup\{0\}}   (b \otimes r_m) q^m\,\,\, \mapsto\,\,\,  \mathrm{Nm}_{F/\mathbb{Q}} v 
\!\!\!\!\!\!\!\!\!\!\sum_{m\in v(NJ)_+^{-1}\cup\{0\}}
( \beta_2^l  b \otimes   r_{\beta_2m} ) q^m.$$
Summing over $i$ then gives the desired formula.\epf

\subsection{Hecke operators at primes dividing the level} 
We shall also make use of the operator $T_v = [U \left(\begin{array}{cc}\varpi_v&0\\0&1 \end{array}\right)U]$
on $M_{k,l}(U;R)$ for $U =  U_1(\mathfrak{n})$ and $v|\mathfrak{n}$.  Note that the operators $T_v$
on $M_{k,l}(U_1(\mathfrak{n});R)$ for all $v\nmid p$ commute with each other, as well as the $S_v$
for $v\nmid p\mathfrak{n}$.   The effect of $T_v$ on $q$-expansions for $v|\mathfrak{n}$ is
computed exactly as in the proof of Proposition~\ref{prop:Tvonqexps} except for the absence of the
coset indexed by $i=\infty$:

\begin{proposition} \label{prop:Uvonqexps}
  If $f \in M_{k,l}(U_1(\mathfrak{n});R)$, $v| \mathfrak{n}$ and $m \in J_+^{-1} \cup \{0\}$, then
$$r_m^J(T_vf) = \beta_1^l  r_{\beta_1 m}^{J_1}(f),$$
where the $J_1$ and $\beta_1 \in F_+$ are such that $vJ = \beta_1J_1$.
\end{proposition}

\subsection{Hecke operators at primes dividing $p$}  \label{subsec:Uv}
We also need the operators $T_v$ for $v|p$ in the case $R = E$, $l_\tau = 0$, $k_\tau \ge 2$ for all $\tau$;
we recall the definition.  Again let $J,J_1,\beta_1$ (in $F_+$) be such that $vJ = \beta_1J_1$.
Let $A_1 = A_{J_1,N}$ denote the universal HBAV over $\overline{Y}_{J_1,N}$.   Letting $H$ denote the
kernel of $\mathrm{Ver}_{A_1}:A_1^{(p)} \to A_1$,  we may decompose $H = \prod_{w|p}  H_w$ where
each $H_w$ is a free rank one $(O_F/w)$-module scheme over $\overline{Y}_{J_1,N}$ and set
$A_1' = A_1^{(p)}/H_v'$ where $H_v' = \prod_{w\neq v} H_w$.    The $O_F$-action on $A_1^{(p)}$, 
polarisation $p\beta_1^{-1}\lambda^{(p)}$ and level $N$-structure $p^{-1}\eta^{(p)}$ induce ones
on $A_1'$ making it a $J$-polarised HBAV over $\overline{Y}_{J_1,N}$, corresponding to a finite,
flat morphism
$$\tilde{\rho}:  \overline{Y}_{J_1,N} \to \overline{Y}_{J,N}$$
of degree $\mathrm{Nm}_{F/\mathbb{Q}}v$.   Taking the union over $J_1$,  the resulting
morphism descends, for sufficiently small $U \supset U(N)$, to a finite, flat  endomorphism of $\overline{Y}_U$
which we denote by $\rho$.  

To define $T_v$, recall that the Kodaira-Spencer isomorphism (\ref{eqn:KS1}) induces
$\overline{\mathcal{L}}_U^{2,-1} \cong \overline{\mathcal{K}}_U$, where $\overline{\mathcal{K}}_U$
is the dualising sheaf on $\overline{Y}_U$, and hence an isomorphism
$$\overline{\mathcal{L}}_U^{k,0} \cong \overline{\mathcal{K}}_U \otimes_{\mathcal{O}_{\overline{Y}_U}} 
\overline{\mathcal{L}}_U^{k-2,1}.$$
Letting $s:A \to \overline{Y}_{J,N}$ and $s_1:A_1 \to \overline{Y}_{J_1,N}$ denote the structure
morphisms for the universal HBAV's, the isogenies $\pi:A_1' \to A_1$ induced by $\mathrm{Ver}_{A_1}$
yield morphisms $s_{1,*}\Omega^1_{A_1/\overline{Y}_{J_1,N}} \to \tilde{\rho}^*s_*\Omega^1_{A/\overline{Y}_{J,N}}$,
which in turn induce morphisms $\overline{\mathcal{L}}_{J_1,N}^{k-2,0} \to \tilde{\rho}^*\overline{\mathcal{L}}_{J,N}^{k-2,0}$
(using here that $k_\tau \ge 2$ for all $\tau$) whose union over $J_1$ descends to 
\begin{equation} \label{eqn:verstar} \overline{\mathcal{L}}_U^{k-2,0} \to \rho^*\overline{\mathcal{L}}_U^{k-2,0}.\end{equation}
Making use of the canonical trivialisations 
$\overline{\mathcal{L}}_{J,N}^{0,1} \simeq \mathrm{Nm}_{F/\mathbb{Q}} (J \mathfrak{d}^{-1})\otimes 
\mathcal{O}_{\overline{Y}_{J,N}}$ and 
$\overline{\mathcal{L}}_{J_1,N}^{0,1} \simeq \mathrm{Nm}_{F/\mathbb{Q}} (J_1 \mathfrak{d}^{-1})\otimes 
\mathcal{O}_{\overline{Y}_{J_1,N}}$,
we define $\overline{\mathcal{L}}_{J_1,N}^{0,1} \stackrel{\sim}{\to} \tilde{\rho}^*\overline{\mathcal{L}}_{J,N}^{0,1}$
by multiplication by $\mathrm{Nm}_{F/\mathbb{Q}}(JJ_1^{-1})$.  (Note that this is {\em not} the
morphism induced by $\pi$, which is in fact $0$.)  The union over $J_1$ then descends to an isomorphism
$\overline{\mathcal{L}}_U^{0,1} \stackrel{\sim}{\to} {\rho}^*\overline{\mathcal{L}}_U^{0,1}$, which we tensor
with (\ref{eqn:verstar}) to obtain a morphism
$\overline{\mathcal{L}}_U^{k-2,1}  \to \rho^*\overline{\mathcal{L}}_U^{k-2,1}$.
We then define $T_v$ as the composite 
$$\begin{array}{l}
H^0(\overline{Y}_U,\overline{\mathcal{K}}_U\otimes_{\mathcal{O}_{\overline{Y}_U}} \overline{\mathcal{L}}_U^{k-2,1})
 \to H^0(\overline{Y}_U,\overline{\mathcal{K}}_U\otimes_{\mathcal{O}_{\overline{Y}_U}} \rho^*\overline{\mathcal{L}}_U^{k-2,1})
 \\ \ \qquad \stackrel{\sim}{\to}
 H^0(\overline{Y}_U,\rho_*\overline{\mathcal{K}}_U\otimes_{\mathcal{O}_{\overline{Y}_U}} \overline{\mathcal{L}}_U^{k-2,1}) \to
H^0(\overline{Y}_U,\overline{\mathcal{K}}_U\otimes_{\mathcal{O}_{\overline{Y}_U}}\overline{\mathcal{L}}_U^{k-2,1}),\end{array}$$
where the first map is given by the one just defined, the second is the canonical isomorphism, and the third is induced
by the trace map $\rho_*\overline{\mathcal{K}}_U \to \overline{\mathcal{K}}_U$.

\begin{proposition}  \label{prop:Uponqexps} Suppose that $v|p$, $l=0$ and $k_\tau \ge 2$ for all $\tau \in \Sigma$.
If $f \in M_{k,l}(U_1(\mathfrak{n});E)$ and $m \in J_+^{-1} \cup \{0\}$, then
$$r_m^J(T_vf) = r_{\beta_1 m}^{J_1}(f),$$
where the $J_1$ and $\beta_1 \in F_+$ are such that $vJ = \beta_1J_1$.
\end{proposition}
\begpf Let $C$ be a cusp at infinity on $X_U$ where $U = U_1(\mathfrak{n})$, so that
$C$ is the image of a cusp $\tilde{C}$ of $X_{J,N}$ associated to the Tate semi-HBAV
$T_{\mathfrak{a},\mathfrak{b}}$ with $\mathfrak{a} = O_F$, $\mathfrak{b} = J^{-1}$, canonical polarisation data
and level structure $\eta(x,y) = t(y)q^{s(x)}$ for some choice of $s$ and $t$.
The morphisms $\tilde{\rho}$ extend uniquely to morphisms
$\overline{X}_{J_1,N} \to \overline{X}_{J,N}$, for which one finds that the
fibre over $\tilde{C}_E$ is supported at $\tilde{C}_{1,E}$, where $\tilde{C}_1$ is the
cusp of $X_{J_1,N}$ associated to
$T_{O_F,J_1^{-1}}$, with canonical polarisation data
and level structure $\eta(x,y) = t(y)q^{s_1(x)}$ for some choice of $s_1$.

Moreover the corresponding map $\widehat{S}_{\tilde{C},E} \to \widehat{S}_{\tilde{C}_1,E}$ is defined by
$q^m \mapsto q^{\beta_1m}$, and the pullback of the isogeny $\pi$ to
$\widehat{S}_{1,E}$ is the canonical projection
$\tilde{\rho}^*T_{O_F,J^{-1},E} \to T_{O_F,J_1^{-1},E}$ induced by the identity on
$\mathbb{G}_m\otimes \mathfrak{d}^{-1}$ (where $\widehat{S}_1 = 
S_N(\{\sigma_{\alpha_1}^{\tilde{C}_1}\})^\wedge$ is again as in \cite[3.4.2]{C2}
for a suitable cone decomposition $\{\sigma_\alpha^{\tilde{C}}\}$).
In particular, it follows that the morphism (\ref{eqn:verstar}) is compatible with the canonical
trivialisations over $\widehat{S}_{1,E}$, so the resulting map
$$Q_{C_1,E}^{k-2,1} \to \widehat{S}_{C_1,E}  \otimes_{\widehat{S}_{C,E}} Q_{C,E}^{k-2,1}$$
is induced by multiplication by $\mathrm{Nm}_{F/\mathbb{Q}}(JJ_1^{-1})$.

Identifying the pullback of $j_*\mathcal{K}_U$ to $\widehat{S}$
with $\mathrm{Nm}_{F/\mathbb{Q}}(J)^{-1} \otimes \mathcal{O}_{\widehat{S}}$, and similarly for $\widehat{S}_1$ with
$C$ and $J$ replaced by $C_1$ and $J_1$, we find that the trace
$\rho_*\overline{\mathcal{K}}_U \to \overline{\mathcal{K}}_U$ extends (over $\overline{X}_U$)
to the map whose pull-back to $\widehat{S}_E$ is defined by 
$$b \otimes q^m \mapsto \left\{\begin{array}{ll}   \mathrm{Nm}_{F/\mathbb{Q}}(J_1J^{-1})b  \otimes 
q^{\beta_1^{-1}m},&\mbox{if $m \in vJ_1^{-1}$;} \\
0,&\mbox{otherwise.}\end{array}\right.$$
It follows from \cite[(1.1.20)]{K} that the same holds for the
corresponding map $\rho_*\overline{\mathcal{L}}_U^{2,-1} \to \overline{\mathcal{L}}_U^{2,-1}$,
so this formula describes the resulting map 
$Q_{C_1,E}^{2,-1} \to Q_{C,E}^{2,-1}$, giving the proposition.
\epf

One easily sees directly from the definitions that the $T_v$ for $v|p$ 
on $M_{k,0}(U_1(\mathfrak{n});E)$  commute with the
$S_v$ for $v \nmid p\mathfrak{n}$ (assuming all $k_\tau \ge 2$), and it follows from
Propositions~\ref{prop:Tvonqexps},
\ref{prop:Uvonqexps} and~\ref{prop:Uponqexps}
that they commute with each other as well as the $T_v$ for all $v\nmid p$. 
(In fact one can check directly from the definitions that the $T_v$ commute with each other and
the action of the group
$\{\,g\in \mathrm{GL}_2(\mathbb{A}_F^\infty)\,|\, g_p \in \mathrm{GL}_2(O_{F,p})\,\}$ on
$M_{k,0}(E)$.)

\subsection{Partial Frobenius operators}
We also define operators $\Phi_v$ for $v|p$ in the case $R=E$, $l=0$, generalising the
classical $V$-operator.  We maintain the notation from the definition of $T_v$ in \S\ref{subsec:Uv}, except that
we no longer assume $k_\tau \ge 2$ for all $\tau$.  Writing $s_1':A_1' \to \overline{Y}_{J_1,N}$
and  $\tilde{\rho}^* s_* \Omega^1_{A/\overline{Y}_{J,N}} \simeq (s_1')^* \Omega^1_{A_1'/\overline{Y}_{J_1,N}}
= \bigoplus \overline{\omega}_\tau'$, we find the isogenies $A_1^{(p)}  \to A_1' \to A_1$
induce isomorphisms 
\begin{equation}\label{eqn:Phivdef} 
\tilde{\rho}^* \overline{\omega}_\tau \simeq \overline{\omega}_\tau' \simeq \left\{ \begin{array}{ll}
\overline{\omega}_{\mathrm{Fr}^{-1}\circ\tau}^{\otimes p},&\mbox{if $\tau \in\Sigma_v$,}\\
\overline{\omega}_\tau,&\mbox{if $\tau \not\in\Sigma_v$,}\end{array}\right.
\end{equation}
on $\overline{Y}_{J_1,N}$ whose unions over $J_1$ descend to $\overline{Y}_U$.
For $k \in \mathbb{Z}^\Sigma$, define $k'$ by $k'_\tau = pk_{\mathrm{Fr}\circ\tau}$ if $\tau \in \Sigma_v$
and $k'_\tau = k_\tau$ if $\tau\not\in \Sigma_v$, and $\Phi_v:M_{k,0}(U;E) \to M_{k',0}(U;E)$
as the composite
$$M_{k, 0}(U; E)=H^0(\overline{Y}_U, \overline{\mathcal{L}}_U^{k,0})
\to H^0(\overline{Y}_U, \rho^*\overline{\mathcal{L}}_U^{k,0})
\to H^0(\overline{Y}_U, \overline{\mathcal{L}}_U^{k',0})=M_{k', 0}(U; E),$$
where the first map is pull-back and the second is induced by the above isomorphisms.

It is clear from the definition that $\Phi_v$ is injective, and straightforward to check
the operators $\Phi_v$ commute with each other and the action of the groups
$\{\,g\in \mathrm{GL}_2(\mathbb{A}_F^\infty)\,|\, g_p \in \mathrm{GL}_2(O_{F,p})\,\}$ on
$M_{k,0}(E)$ and $M_{k',0}(E)$.  In particular, $\Phi_v$ commutes with the operators
$T_w$ for all $w\nmid p$ and $S_w$ for all $w\nmid \mathfrak{n}p$.
Its effect on $q$-expansions is given by the following:

\begin{proposition}  \label{prop:Phivonqexps} Suppose that $v|p$ and $l=0$.
If $f \in M_{k,l}(U_1(\mathfrak{n});E)$ and $m \in J_+^{-1} \cup \{0\}$, then
$$r_m^J(\Phi_vf) = r_{\beta_	2 m}^{J_2}(f),$$
where the $J_2$ and $\beta_2 \in F_+$ are such that $v^{-1}J = \beta_2J_2$
(interpreting $r_{\beta_2 m}$ as $0$ if $\beta_2m \not\in J_2^{-1}$, i.e., $m\not\in vJ^{-1}$).
\end{proposition}
\begpf  The completion of $\rho$ at the cusps is already computed in
the course of the proof of Proposition~\ref{prop:Uponqexps}.  One finds also
that the pull-back of the isomorphisms of (\ref{eqn:Phivdef}) to $\widehat{S}_{1,E}$
are compatible with the canonical trivialisations of the pushforwards of the cotangent bundles
of the Tate semi-HBAVs $T_{O_F,J^{-1},E}$ and $T_{O_F,J_1^{-1},E}$.  It follows that the map
\begin{equation}\label{eqn:Phivq}
Q_{C,E}^{k,0} \to \widehat{S}_{C_1,E}  \otimes_{\widehat{S}_{C,E}} Q_{C,E}^{k,0}
   \simeq Q_{C_1,E}^{k',0}\end{equation}
induced by $\Phi_v$ is defined by $q^m \mapsto q^{\beta_1m}$.
The desired formula follows on relabelling $J$ as $J_2$ and $J_1$ as $J$,
and taking $\beta_2 = \beta_1^{-1}$.
\epf

The proposition gives an alternative proof (for $U= U_1(\mathfrak{n})$) that the $\Phi_v$ commute
with each other and the $T_w$ for all $w\nmid p$ (after checking that $\Phi_v$ commutes with $S_w$
and applying Proposition~\ref{prop:Tvonqexps} for $w\nmid \mathfrak{n}p$,
and Proposition~\ref{prop:Uvonqexps} for $w|\mathfrak{n}$).
Note however that $\Phi_v$ does not commute with $T_v$ (when the latter
is defined, i.e., $k_\tau \ge 2$ for all $\tau$).

Note that it is immediate from Proposition~\ref{prop:Thetaq} that the kernel
of the operator $\Theta_\tau$ depends only on the prime $v$ such that $\tau \in \Sigma_v$.
Moreover if $\nu(f) = k$ (in the notation of \S\ref{sec:Hasse}), then Theorem~\ref{thm:theta}
implies that $k_\tau$ is divisible by $p$ for all $\tau \in \Sigma_v$.  We will show in
Theorem~\ref{thm:Phiv} that (assuming $k$ is of this form and $l=0$) this kernel is in fact
the image of $\Phi_v$, generalising a result of Katz in Section II of \cite{K2}.  

We need to introduce one more operator:  we define $k^\varphi$
by $k^\varphi_\tau = k_{\mathrm{Fr}^{-1}\circ\tau}$ and 
$\varphi: M_{k,0}(U;E)  \to M_{k^\varphi,0}(U;E)$ as the composite:
$$H^0(\overline{Y}_U, \overline{\mathcal{L}}_U^{k,0})
\to H^0(\overline{Y}_U, \mathrm{Fr}_E^*\overline{\mathcal{L}}_U^{k,0})
\to H^0(\overline{Y}_U, \overline{\mathcal{L}}_U^{k^\varphi,0}),$$
where the first map is pull-back by the automorphism induced by $\mathrm{Fr}_E$ on $\overline{Y}_U$
and the second is given by the
canonical isomorphisms $\mathrm{Fr}_E^*\overline{\omega}_\tau \simeq \overline{\omega}_{\mathrm{Fr}\circ\tau}$.
(Note that we could similarly define $\varphi: M_{k,l}(U;E)  \to M_{k^\varphi,l^\varphi}(U;E)$.)
Its effect on $q$-expansions of $f \in M_{k,0}(U_1(\mathfrak{n});E)$ is given by
$r_m^J(\varphi f) = \mathrm{Fr}_E(r_m^J(f)) = (r_m^J(f))^p$.

\begin{theorem}  \label{thm:Phiv}   Suppose that $k \in \mathbb{Z}^\Sigma$,
$\mathfrak{n}$ is an ideal of $O_F$ prime to $p$, $v$ is a prime dividing $p$
and $\tau \in \Sigma_v$.  Then the image of 
$$\Phi_v : M_{k,0}(U_1(\mathfrak{n});E) \to M_{k',0}(U_1(\mathfrak{n});E)$$
is the kernel of $\Theta_\tau$.
\end{theorem}
\begpf From Proposition~\ref{prop:Thetaq} we see that $f \in M_{k',0}(U_1(\mathfrak{n});E)$ is in the
kernel of $\Theta_\tau$ if and only if $r_m^J(f) = 0$ for all $m,J$ such that $m \not\in vJ^{-1}_+$
It is therefore immediate from Proposition~\ref{prop:Phivonqexps} that 
$\mathrm{image}(\Phi_v) \subset \ker(\Theta_\tau)$.

For the opposite inclusion, first note that we can assume $\mathfrak{n}$ is sufficiently small.
For each cusp $C \in \mathcal{S}$, let $N_C$ denote the stalk
$(j_*\overline{\mathcal{L}}_U^{k,0})_{\overline{C}}$, where as usual $j$ is the inclusion
$\overline{Y}_U \to \overline{X}_U$ and $\overline{C} = C_E$.  Similarly let 
$$N_C' =  (j_*\overline{\mathcal{L}}_U^{k',0})_{\overline{C}_1} =
(j_*\rho_*\overline{\mathcal{L}}_U^{k',0})_{\overline{C}}$$
and consider the $R_C := \mathcal{O}_{\overline{X}_U,\overline{C}}$-linear map 
$\phi_C: N_C \to N_C'$ of finitely generated $R_C$-modules induced by the morphisms
in the definition of $\Phi_v$.   Letting $\mathcal{F}_U$ denote the sheaf of total fractions
on $\overline{Y}_U$, we similarly have a map
$$\widetilde{\Phi}_v:  H^0(\overline{Y}_U, \overline{\mathcal{L}}_U^{k,0}
\otimes_{\mathcal{O}_{\overline{Y}_U}}  \mathcal{F}_{U}) \to H^0(\overline{Y}_U, \overline{\mathcal{L}}_U^{k',0}
\otimes_{\mathcal{O}_{\overline{Y}_U}}  \mathcal{F}_{U}),$$
and thus a commutative diagram of injective maps:
$$\begin{array}{ccccc}
M_{k,0}(U_1(\mathfrak{n});E)  & \to & \bigoplus_{C\in \mathcal{S}} N_C  & \to & H^0(\overline{Y}_U, \overline{\mathcal{L}}_U^{k,0}
\otimes_{\mathcal{O}_{\overline{Y}_U}}  \mathcal{F}_{U}) \\
\downarrow&&\downarrow&&\downarrow\\
M_{k',0}(U_1(\mathfrak{n});E)  & \to & \bigoplus_{C\in \mathcal{S}} N'_C  & \to & H^0(\overline{Y}_U, \overline{\mathcal{L}}_U^{k',0}
\otimes_{\mathcal{O}_{\overline{Y}_U}}  \mathcal{F}_{U})\end{array}$$
where the horizontal maps are the natural inclusions.

The completion $\widehat{\phi}_C$ of $\phi_C$ is precisely the $\widehat{S}_{\overline{C}}$-linear
map $Q^{k,0}_{\overline{C}} \to Q^{k',0}_{\overline{C}_1}$ of (\ref{eqn:Phivq}) (where
$\widehat{S}_{\overline{C}}$ acts on the target via the map to $\widehat{S}_{\overline{C}_1}$
induced by $\rho$).  If $f \in \ker(\Theta_\tau)$, then $r_m^{J_1}(f) = 0$ for all $m \not\in vJ^{-1}_1$,
so the $q$-expansion of $f$ at $C_1$ is in the image of $\widehat{\phi}_C$ for each $C \in \mathcal{S}$.
Since $\widehat{S}_{\overline{C}}$ is faithfully flat over $R_C$, it follows that
$f$ is in the image $\phi_C$ for each $C \in \mathcal{S}$, so there
exists $g \in \bigoplus_{C} N_C \subset H^0(\overline{Y}_U, \overline{\mathcal{L}}_U^{k,0}
\otimes_{\mathcal{O}_{\overline{Y}_U}}  \mathcal{F}_U)$ such that $\widetilde{\Phi}_v(g) = f$.

It remains to prove that $g \in M_{k,0}(U_1(\mathfrak{n});E)$,   Since $\overline{Y}_U$ is smooth
and $\overline{\mathcal{L}}_U^{k,0}$ is invertible, it suffices to prove that $\mathrm{ord}_z (g) \ge 0$
for all prime divisors $z$ on $\overline{Y}_U$.  For this, we note that the map $\varphi$ similarly
extends to a map $\widetilde{\varphi}$, and one checks that
$$\widetilde{\varphi} \circ \prod_{w|p} \widetilde{\Phi}_w : 
H^0(\overline{Y}_U, \overline{\mathcal{L}}_U^{k,0}
\otimes_{\mathcal{O}_{\overline{Y}_U}}  \mathcal{F}_{U}) \to H^0(\overline{Y}_U, \overline{\mathcal{L}}_U^{pk,0}
\otimes_{\mathcal{O}_{\overline{Y}_U}}  \mathcal{F}_{U})$$
is simply the map $g \mapsto g^p$.  Therefore
$$g^p = \left(\widetilde{\varphi}\circ  \prod_{w\neq v} \widetilde{\Phi}_w \right) (f)
           = \left(\varphi \circ \prod_{w\neq v} \Phi_w\right) (f)  \in M_{pk,0}(U_1(\mathfrak{n});E),$$
so that $p \,\mathrm{ord}_z (g) = \mathrm{ord}_z( g^p )\ge 0$, and hence 
$\mathrm{ord}_z (g) \ge 0$.  \epf

\begin{remark} One can also check that the relation $\mathrm{image}(\Phi_v) = \ker(\Theta_\tau)$ holds
for arbitrary $U$ using the same argument as in the proof of the theorem and a straightforward
generalisation of the formula in Proposition~\ref{prop:Phivonqexps} (see the next section
for similar computations of the effect of operators on $q$-expansions at more general cusps).\end{remark}

\section{Normalised eigenforms}
\label{sec:eigenforms}

We will prove that if $\rho$ is irreducible and geometrically modular of weight $(k,l)$, then
in fact $\rho$ is associated to an eigenform $f \in M_{k,l}(U_1(\mathfrak{n});E)$
for some $\mathfrak{n}$ prime to $p$, allowing us to pin down $q$-expansions
of forms giving rise to $\rho$.   We will also use partial $\Theta$-operators to
study the behaviour of minimal weights as $l$ varies, and prove that if an eigenform
is ordinary at a prime over $v$, then so is the associated Galois representation.

\subsection{Preliminaries}
First note that, by definition, if $\rho$ is geometrically
modular of weight $(k,l)$, then $\rho$ is associated to an eigenform
$f \in M_{k,l}(U(\mathfrak{n});E)$ for some $\mathfrak{n}$ prime to $p$.
One approach to replacing $U(\mathfrak{n})$ by $U_1(\mathfrak{n}')$
for some $\mathfrak{n}'$ would be to use the space $M_{k,l}(E)$ to
associate to $\rho$ a representation of $\mathrm{GL}_2(F_v)$ for each $v|\mathfrak{n}$.
One then chooses an irreducible subrepresentation $\pi_v$, whose existence is given
by~\cite[II, 5.10]{V2}, and shows, using the irreducibility of $\rho$,
that $\pi_v$ does not factor through $\det$.  It then follows from~\cite{V} that 
$\pi_v$ has a vector invariant under $U_1(v^{c_v})$ for some exponent $c_v$,
and one can take $\mathfrak{n}' = \prod_c v^{c_v}$.
We shall instead give a more constructive argument that develops some tools
we will need anyway.  In particular we define certain twisting operators
on forms of level $U(\mathfrak{n})$.

We let $U = U(\mathfrak{n})$ and index the components of $Y_U$ by pairs $(J,w)$ where $J$, as usual,
runs through strict ideal class representatives, and $w$ runs through a set $W \subset (O_F/NO_F)^\times$
of representatives for $(O_F/\mathfrak{n})^\times/O_{F,+}^\times$.  More precisely,  choose as before
isomorphisms $s: O_F/NO_F \simeq (NJ)^{-1}/J^{-1}$ (for each $J$) and 
$t:O_F/NO_F \simeq \mu_N \otimes \mathfrak{d}^{-1}$.   Then $s$ determines an
isomorphism $J/NJ \simeq O_F/NO_F$ whose composite with $t$ defines a component
of $Y_{J,N}$, hence of $Y_U$, and we associate to $(J,w)$ the component so defined
with $s$ replaced by $ws$. One easily checks that this defines a bijection between
$Z_U(\mathcal{O})$ and the set of such pairs.  Moreover, there is a unique cusp on
each component of $X_U$ mapping to a cusp at $\infty$ on $X_{U_1(\mathfrak{n})}$,
namely the one associated to the Tate semi-HBAV $T_{O_F,J^{-1}}$ with canonical polarisation
and level $N$ structure $(x,y) \mapsto t(y)q^{ws(x)}$.  For $f \in M_{k,l}(U;R)$ and
$m \in (\mathfrak{n}^{-1}J)_+\cup\{0\}$, we write $r_m^{J,w}(f)$ for the corresponding
$q$-expansion coefficient of $f$.

A computation similar to the proof of Proposition~\ref{prop:Tvonqexps} shows that
the effect of $T_v$ on $q$-expansions of forms in $M_{k,l}(U;R)$ is given by the formula:
\begin{equation}\label{eqn:TvUn} r_m^{J,w}(T_vf) = \beta_1^l  r_{\beta_1 m}^{J_1,w_1}(f)  + \mathrm{Nm}_{F/\mathbb{Q}}(v) \beta_2^l r_{\beta_2 m}^{J_2,w_2}(S_vf)\end{equation}
for $m \in (\mathfrak{n}J)_+^{-1}\cup\{0\}$, where the $J_i$, $\beta_i$ are as before,
with $w_i \in W$ satisfying
$$\begin{array}{rcl}   \beta_1\varpi_v^{-1}ws(1)& \equiv &w_1s_1(1) \bmod \mathfrak{n}(NJ_1)^{-1}\\
\mbox{and}\quad \beta_2\varpi_vws(1)& \equiv& w_2s_2(1) \bmod \mathfrak{n}(NJ_2)^{-1},\end{array}$$
where we view $\beta_1\varpi_v^{-1}$ as inducing an isomorphism $J^{-1}\widehat{O}_F \simeq J_1^{-1}\widehat{O}_F$,
hence $(NJ)^{-1}/J^{-1} \simeq (NJ_1^{-1})/J_1^{-1}$, and similarly $\beta_2\varpi_v$ as inducing
$(NJ)^{-1}/J^{-1} \simeq (NJ_2)^{-1}/J_2^{-1}$.
We note the following consequence:
\begin{lemma} \label{lem:Eisenstein}   If $f \in M_{k,l}(U(\mathfrak{n});E)$ is an eigenform for the operators
$T_v$ and $S_v$ for all but finitely many $v$, and the associated Galois representation $\rho_f$ is
absolutely irreducible, then $r_0^{J,w} = 0$ for all pairs $(J,w)$.
\end{lemma}
\begpf\footnote{Alternatively, this can be proved by revisiting the construction in
Theorem~\ref{thm:galois} and observing that if $r_0^{J,w} \neq 0$ for some $(J,w)$, then
the lift $\tilde{f}$ is non-cuspidal.  One then deduces that
the Galois representation $\rho_{\tilde{f}}$ is reducible, and hence so is $\rho_f$ (possibly after extending
scalars in the case $p=2$).} If $v$ is trivial in the strict class group of conductor $\mathfrak{n}p$, then it follows
from the definitions that $S_v$ acts trivially on $M_{k,l}(U;E)$.  Moreover in the formula (\ref{eqn:TvUn})
for such $v$, we have $J=J_1=J_2$, $w=w_1=w_2$,  $\beta_1 \equiv \beta_2 \equiv 1 \bmod pO_{F,p}$
and $\mathrm{Nm}_{F/\mathbb{Q}}(v) \equiv 1 \bmod p$, so that
$r_0^{J,w}(T_vf) = 2r_0^{J,w}(f)$.  Therefore if $r_0^{J,w}(f) \neq 0$ for some $(J,w)$, then
$\rho_f(\mathrm{Frob}_v)$ has characteristic polynomial $(X-1)^2$ for such $v$.  By the Cebotarev
Density Theorem (and class field theory) it follows that $\rho_f(g)$ has characteristic polynomial
$(X-1)^2$ for all $g \in G_K$, where $K$ is the strict ray class field over $F$ of conductor $\mathfrak{n}p$,
so by the Brauer--Nesbitt Theorem, $\rho_f|_{G_K}$ has trivial semi-simplification.  Since $K$ is abelian
over $F$, this contradicts the absolute irreducibility of $\rho_f$ 
\epf

We continue to assume $U = U(\mathfrak{n})$ and
define an action of  the group $(O_F/\mathfrak{n})^\times$ on $M_{k,l}(U;R)$ via its
isomorphism with the subgroup of $\mathrm{GL}_2(O_F/\mathfrak{n})$ consisting of matrices of the
form $\left(\begin{array}{cc} a&0 \\ 0 & 1 \end{array}\right)$.  Thus $a \in (O_F/\mathfrak{n})^\times$
acts on  $M_{k,l}(U;R)$ by the operator $[UgU]$ for any $g \in \mathrm{GL}_2(\widehat{O}_F)$
congruent to $\left(\begin{array}{cc} a&0 \\ 0 & 1 \end{array}\right) \bmod \mathfrak{n}$; we denote
the operator by $\langle a \rangle$.  It is straightforward to check that its effect on $q$-expansions
is given by the formula:
\begin{equation} \label{eqn:diamond} r_m^{J,w}(\langle a \rangle f) = \nu^l  r_{\nu m}^{J,w'}(f),\end{equation}
where $\nu \in O_{F,+}^\times$ and $w' \in W$ are such that $\nu w \equiv a w' \bmod \mathfrak{n}$.

We also define the  operator $T_v = [U \left(\begin{array}{cc}\varpi_v&0\\0&1 \end{array}\right)U]$
on $M_{k,l}(U;R)$ for $v|\mathfrak{n}$ and a choice of uniformiser $\varpi_v$ for $F_v$.  Another computation
similar to Proposition~\ref{prop:Tvonqexps} (or more precisely, Proposition~\ref{prop:Uvonqexps})
shows that its effect on $q$-expansions is given by:
\begin{equation} \label{eqn:UvUn}  r_m^{J,w}(T_v f) = \beta_1^l  r_{\beta_1 m}^{J_1,w_1}(f) \end{equation}
with notation as in (\ref{eqn:TvUn}).
Note that $T_v$ depends on the choice of $\varpi_v$: replacing $\varpi_v$ by $u\varpi_v$ for
$u \in O_{F,v}^\times$ replaces $T_v$ with $\langle \pi(u) \rangle T_v$ where $\pi$ is
the natural map $O_{F,v}^\times \to (O_F/\mathfrak{n})^\times$.  We see directly from the
definitions that the operators $T_v$ for $v|\mathfrak{n}$ commute with the $T_v$ and $S_v$
for $v\nmid p\mathfrak{n}$ (and each other), as well as the action of $(O_F/\mathfrak{n})^\times$.

Suppose that $\xi:(O_F/\mathfrak{n})^\times \to R^\times$ is a character of conductor 
$\mathfrak{m}|\mathfrak{n}$.  Choose an element $c \in O_{F,\mathfrak{n}} = \prod_{v|\mathfrak{n}} O_{F,v}$ generating $\mathfrak{nm}^{-1}O_{F,\mathfrak{n}}$, and define a {\em twisting operator} on $M_{k,l}(U;R)$ by the formula:
\begin{equation} \label{def:twisting} \Theta_\xi \,\,\,= \!\!\! \sum_{b \in (O_F/\mathfrak{m})^\times}  \!\!\! \xi(b)^{-1}  [U g_b U]
=  \!\!\! \sum_{b \in (O_F/\mathfrak{m})^\times} \!\!\! \xi(b)^{-1} g_b,\end{equation}
where $g_b \equiv \left(\begin{array}{cc} 1 & bc \\ 0 & 1 \end{array}\right) \bmod \mathfrak{n}$.
The operator  $\Theta_\xi$ commutes with the operators $T_v$ and $S_v$ for $v\nmid \mathfrak{n}p$,
and it is straightforward to check that 
\begin{equation}  \label{eqn:twist}  \langle a \rangle \circ \Theta_\xi = \xi(a)\Theta_\xi \circ \langle a \rangle \end{equation}
(Note also the dependence on $c$: replacing $c$ by $uc$ for $u \in O_{F,\mathfrak{n}}^\times$
replaces $\Theta_\xi$ by $\xi(u)\Theta_\xi$.)  One finds the effect on $q$-expansions
is given by:
\begin{equation}  \label{eqn:gauss}  r_m^{J,w}(\Theta_\xi(f))  =  G_J(\xi,w^{-1}cm) r_m^{J,w}(f), \end{equation}
where $G_J(\xi,m) = \displaystyle\sum_{b \in (O_F/\mathfrak{m})^\times} \!\!\! \xi(b)^{-1} \zeta(-bm)$ 
for $m \in (\mathfrak{m}J)_+^{-1}\cup\{0\}$.  (Recall that $\zeta$ is the homomorphism $(NJ)^{-1}/J^{-1} \to \mu_N$
induced by the trace and our choices of $s$ and $t$; see the discussion before Proposition~\ref{prop:QCkl}.)
Standard results on Gauss sums show that  $G_J(\xi,am) = \xi(a)G_J(\xi,m)$
for all $a \in O_F$, $m \in (\mathfrak{m}J)^{-1}$, where $\xi$ as viewed as a function $O_F \to R$ by
setting $\xi(a) = 0$ for $a$ not prime to $\mathfrak{m}$.  One deduces that if $m$ generates
$(\mathfrak{m}J)^{-1}/J^{-1}$, then
$$G_J(\xi,m)G_J(\xi^{-1},-m) = \mathrm{Nm}_{F/\mathbb{Q}}(\mathfrak{m})$$
(in particular, $G_J(\xi,m) \in R^\times$),  and otherwise $G_J(\xi,m) = 0$. 

\subsection{Eigenforms of level $U_1(\mathfrak{n})$}
\begin{lemma} \label{lem:levelU1}  If $\rho:G_F \to \mathrm{GL}_2(\overline{\mathbb{F}}_p)$ is irreducible
and geometrically modular of weight $(k,l)$, then
$\rho$ arises from an eigenform of weight $(k,l)$ and level $U_1(\mathfrak{n})$ for some $\mathfrak{n}$ prime to $p$; i.e., 
there exist $\mathfrak{n}$ prime to $p$, a  field $E$ and an eigenform $f \in M_{k,l}(U_1(\mathfrak{n});E)$
for $S_v$ and $T_v$ for all $v \nmid \mathfrak{n}p$ such that
$\rho \simeq \rho_f$.
\end{lemma}
\begpf By assumption, there exist $\mathfrak{n}$ (prime to $p$), $E$ and
$f \in M_{k,l}(U(\mathfrak{n});E)$, an eigenform  for all $S_v$ and $T_v$ 
with $v \nmid \mathfrak{n}p$, such that
$\rho \simeq \rho_f$.   Since the action of $(O_F/\mathfrak{n})^\times$ commutes
with the operators $S_v$ and $T_v$, we can further assume that $f$ is an
eigenform for this action, i.e., that there is a character $\xi:(O_F/\mathfrak{n})^\times \to E^\times$
(enlarging $E$ if necessary) such that $\langle a \rangle f = \xi(a) f$ for all $a \in (O_F/\mathfrak{n})^\times$.

By Lemmas~\ref{lem:qexp} and~\ref{lem:Eisenstein}, we must have $r_m^{J,w}(f) \neq 0$
for some $J$, some (hence all by (\ref{eqn:diamond})) $w\in W$ and some 
$m \in (\mathfrak{n}J)_+^{-1}$ (i.e., $m\neq 0$).  Letting $e_v = \mathrm{ord}_v(m\mathfrak{n}J)$
for $v|\mathfrak{n}$ and $f' = \prod_{v|\mathfrak{n}}  T_v^{e_v} f$, formula (\ref{eqn:UvUn}) implies that
$r_{m'}^{J',w'}(f') \neq 0$ for some $(J',w')$ and $m' \in (\mathfrak{n}J')_+^{-1}$ with
$m'\mathfrak{n}J'$  prime to $\mathfrak{n}$ (choose $J'$ equivalent to $\displaystyle\prod_{v|\mathfrak{n}} v^{-e_v}J$
and let $m' = \displaystyle\prod_{v|\mathfrak{n}} \beta_{1,v}^{-e_v}m$).
Replacing $f$ by $f'$, we now have
$\langle a \rangle f = \xi(a) f$ for all $a \in (O_F/\mathfrak{n})^\times$, and
$r_m^{J,w}(f) \neq 0$ for some $(J,w)$ and $m \in (\mathfrak{n}J)_+^{-1}$ generating
$(\mathfrak{n}J)^{-1}/J^{-1}$.  

Now replace $f$ by $\Theta_{\xi^{-1}}(f)$ where $\Theta_{\xi^{-1}}$ is the twisting operator associated to $\xi^{-1}$ as defined in (\ref{def:twisting}). Since $cm$ generates $(\mathfrak{m}J)^{-1}/J^{-1}$,
we have $G_J(\xi^{-1},cm) \neq 0$, so formula (\ref{eqn:gauss}) shows that $f \neq 0$.
By formula (\ref{eqn:twist}), $f$ is invariant under the action of $(O_F/\mathfrak{n})^\times$,
hence under the action of the open compact subgroup
$$U'  = \left\{\,\left.\left(\begin{array}{cc} a&b \\ c & d \end{array}\right) \in U_1(\mathfrak{n}) \,\right|\,b \in \mathfrak{n}\widehat{O}_F\,\right\}.$$
Now let $g = \displaystyle\prod_{v|\mathfrak{n}}  \left(\begin{array}{cc} \varpi_v^{-e_v} & 0 \\ 0 & 1 \end{array}\right)$ where $e_v = \mathrm{ord}_v(\mathfrak{n})$.  Then $g^{-1}U_1(\mathfrak{n}^2)g \subset U'$,
so the lemma follows with $f$ replaced by $[U_1(\mathfrak{n}^2)gU']f$ and $\mathfrak{n}$ replaced by $\mathfrak{n}^2$.
\epf

\subsection{Twisting eigenforms}
Suppose now that $k,l,l' \in \mathbb{Z}^{\Sigma}$, and that $\mathfrak{m},\mathfrak{n}$ are ideals of $O_F$
with $\mathfrak{m}|\mathfrak{n}$ and $\mathfrak{n}$ prime to $p$, and let $V_{\mathfrak{m}} \subset \widehat{O}_F^\times$
denote the kernel of the natural projection to $(O_F/\mathfrak{m})^\times$.  
\begin{definition}\label{def:charwt}  We say a character
$$\xi: \{\,a \in (\mathbb{A}_F^\infty)^\times \,|\, a_p \in O_{F,p}^\times\,\} / V_{\mathfrak{m}}
 \to E^\times$$
is a {\em character of weight} $l'$ if $\xi(\alpha) = \overline{\alpha}^{l'}$ for all 
$\alpha \in F_+^\times \cap O_{F,p}^\times$.  \end{definition}

Suppose that
$f \in M_{k,l}(U_1(\mathfrak{n});E)$ and  $\xi$ is a character of weight $l'$ and conductor $\mathfrak{m}$.
Recall from \S\ref{subsec:exi} that we can associate to $\xi$ a form
$e_\xi \in M_{0,l'}(U(\mathfrak{m});E)$, and hence
a form $e_\xi \otimes f \in M_{k,l+l'}(U(\mathfrak{n});E)$.  Choosing $c = (\varpi_v^{e_v - d_v})_v$ where
$d_v = \mathrm{ord}_v{m}$ and $e_v = \mathrm{ord}_v\mathfrak{n}$,  and applying the following (normalised)
composite of operators from the proof of Lemma~\ref{lem:levelU1}:
$$ \mathrm{Nm}_{F/\mathbb{Q}}(\mathfrak{n})^{-1} [U_1(\mathfrak{n}^2)gU']
 \circ \Theta_{\xi^{-1}} \circ \prod_{v|n}  T_v^{e_v}$$
to $e_\xi \otimes f$ then yields a form in $M_{k,l+l'}(U_1(\mathfrak{n}^2);E)$ which we denote $f_\xi'$.
It is straightforward to check that in fact 
$$f_\xi' = e_\xi \otimes \!\!\! \sum_{b\in (O_F/\mathfrak{m})^\times} \!\!\! \xi(b) \left(\begin{array}{cc} 1 & \tilde{b}c' \\ 0 & 1 \end{array}\right)f  \quad \in M_{k,l+l'}(U_1(\mathfrak{mn});E),$$
where $\tilde{b}$ is any lift of $b$ to $O_{F,\mathfrak{n}}^\times$ and $c' = (\varpi_v^{-d_v})_v$.

We now relate the $q$-expansions of $f$ and $f_\xi'$.  Firstly, the form $e_\xi$ has constant $q$-expansions
satisfying the formula
$$\xi(a) r_0^{J_0,w_0}(e_\xi) =  ||a||^{-1} r_0^{J_1,w_1}\left(\left(\begin{array}{cc} a& 0 \\ 0 & 1 \end{array}\right)e_\xi\right)
   =   \overline{\beta}_1^\ell r_0^{J_1,w_1}(e_\xi),$$
for $a \in (\mathbb{A}_F^\infty)^\times$, $\beta_1 \in F_+^\times$, $w_0,w_1 \in W$ such that
$a_p \in O_{F,p}^\times$, $\beta_1J_1 = (a)J$ and $\beta_1a^{-1}w_0s_0(1) \equiv w_1s_1(1) \bmod \mathfrak{m}(NJ_1)^{-1}$.
Assume for simplicity that $1 \in W$, $O_F$ is chosen as the representative for the trivial ideal class, 
$s_0(1) = N^{-1}$ for $J_0 =  O_F$,
and $e_\xi$ is normalised so that $r_0^{O_F,1}(e_\xi) = 1$.  We then have
$$r_0^{J,w}(e_\xi) = \xi(tw^{-1})$$
where $t$ is chosen so that $J = (t)$ and $t^{-1} \equiv s(1) \bmod \mathfrak{m}(NJ)^{-1}$.
Applying (\ref{eqn:UvUn}) and (\ref{eqn:gauss}) with $l$ replaced by $l+l'$, and
$$ r_m^{J}([U_1(\mathfrak{n}^2)gU'] f) = \mathrm{Nm}_{F/\mathbb{Q}}(\mathfrak{n}) \overline{\beta}_2^{l+l'}
    r_{\beta_2m}^{J_2,w_2}(f)$$
for $f \in M_{k,l+l'}(U';E)$ (where $g$ and $U'$ are as in the proof of Lemma~\ref{lem:levelU1}),
$m \in J_+^{-1}\cup\{0\}$, $\mathfrak{n}^{-1}J = \beta_2J_2$ and
$\left(\beta_2\prod_{v|\mathfrak{n}}\varpi_v^{e_v}\right)s(1) \equiv w_2s_2(1) \bmod \mathfrak{m}NJ_2^{-1}$
then gives the formula\footnote{This also follows more directly from the alternative description of $f_\xi'$
and a formula analogous to (\ref{eqn:gauss}).}
$$r_m^J(f_\xi') =  \xi(t)G_J(\xi^{-1},c'm) r_m^J(f).$$
Let $\xi^{-1}_{\mathfrak{m}}$ denote the character of $(O_F/\mathfrak{m})^\times$ induced by $\xi$,
extended to a map $\widehat{O}_F \to O_F/\mathfrak{m} \to E$ by setting $\xi_\mathfrak{m}^{-1}(a) = 0$
if $(a)$ is not prime $\mathfrak{m}$.  We then have
$G_J(\xi^{-1},c'm) = G_{O_F}(\xi^{-1},c'tm) = \xi_{\mathfrak{m}}^{-1}(tm)G_{O_F}(\xi^{-1},c)$,
so setting \begin{equation} f_\xi = G_{O_F}(\xi^{-1},c')^{-1}f'_\xi\end{equation} gives $f_\xi \in M_{k,l+l'}(U_1(\mathfrak{mn}^2);E)$
satisfying
\begin{equation}  \label{eqn:twistqexp} r_m^J(f_\xi) =  \xi(t)\xi_{\mathfrak{m}}^{-1}(tm) r_m^J(f).\end{equation}
Note that this is independent of the choice of $t$ such that $J = (t)$.  Furthermore if we choose
$t$ so that $t_p = 1$, then $\xi(t) = \xi'(t)$ where $\xi': \mathbb{A}_F^\times/F^\times F_{\infty,+}^\times V_{\mathfrak{m}p}
\to E^\times$ is the character in the proof of Theorem~\ref{thm:galois}.  (Recall that $\xi'$ is defined by
$\xi'(\alpha z a) = \xi(a)\overline{a}_p^{-l'}$ for $\alpha \in F^\times$, $z \in F_{\infty,+}^\times$
and $a \in (\mathbb{A}_F^\infty)^\times$ with $a_p \in O_{F,p}^\times$.)  Since $\xi'(m) = 1$, we then have
$$\xi(t)\xi_{\mathfrak{m}}^{-1}(tm) = \left\{\begin{array}{ll}
\xi'((tm)^{(\mathfrak{m})}),&\mbox{if $(tm)$ is prime to $\mathfrak{m}$;}\\
0,&\mbox{otherwise;}\end{array}\right.$$
where $(tm)^{(\mathfrak{m})}$ denotes the projection of $tm$ to the components prime to $\mathfrak{m}$.

We record the above construction:
\begin{lemma} \label{lem:ontwist} If $f \in M_{k,l}(U_1(\mathfrak{n});E)$ and $\xi$ is a character of weight $l'$ and
conductor $\mathfrak{m}$,  then $f_\xi \in M_{k,l+l'}(U_1(\mathfrak{nm}^2);E)$ has $q$-expansion coefficients
defined by (\ref{eqn:twistqexp}).  In particular if $r_m^J(f) \neq 0$ for some $m \in J_+^{-1}$ with
$mJ$ prime to $\mathfrak{m}$, then $f_\xi \neq 0$, in which case if $f$ is an eigenform, then
so is $f_\xi$, and $\rho_{f_\xi} \simeq \rho_{\xi'} \otimes \rho_f$.
\end{lemma}

\subsection{$\Theta$-operators on eigenforms}

\label{subsec:consequences}

Recall from Corollary~\ref{cor:theta} that $\Theta_\tau$ defines a map
$$M_{k,l}(U_1(\mathfrak{n});E) \to M_{k',l'}(U_1(\mathfrak{n});E),$$
where $k'$ and $l'$ are defined in Definition \ref{def:igusatheta} (in particular $l_{\tau'}' = l_\tau - \delta_{\tau,\tau'}$).
Moreover $\Theta_\tau$ commutes with the operators $T_v$ (for all $v\nmid p$) and
$S_v$ (for all $v\nmid \mathfrak{n}p$).

\begin{lemma}  \label{lem:U1better}
With notation as in Lemma~\ref{lem:levelU1}, we can take the eigenform $f$ in the conclusion
so that $\Theta_\tau(f) \neq 0$ for all $\tau \in \Sigma$.
\end{lemma}
\begpf Let $v$ be a prime dividing $p$, and suppose $\tau \in \Sigma_v$.
Let $f$ be an eigenform in $M_{k,l}(U_1(\mathfrak{n});E)$ giving rise to $\rho$,
and let $m,J$ be such that $r_m^J(f) \neq 0$ (so $m \in J_+^{-1})$.   We wish
to prove that we can choose $f$ with
$r_{m'}^J(f) \neq 0$ for some $m' \not\in vJ_+^{-1}$, so that
$f\not\in \ker(\Theta_\tau)$.

By Chevalley's Theorem, we can (enlarging $E$ if necessary)
choose a character $\xi$ of weight $-l$ and conductor
$\mathfrak{m}$ for some $\mathfrak{m}$ prime to $pmJ$. 
Thus $f_{\xi} \in M_{k,0}(U_1(\mathfrak{n}');E)$ where $\mathfrak{n}'=\mathfrak{nm}^2$,
$r_m^J(f_\xi) \neq 0$, and $f_\xi$ is an eigenform giving rise to 
$\rho' = \rho_{\xi'}\otimes\rho$.  

For eigenforms $g \in M_{k,0}(U_1(\mathfrak{n}');E)$ giving rise to $\rho'$, 
define $\delta_v(g)$ to be the least $d \ge 0$ such that $r_m^J(g) \neq 0$
for some $m,J$ such that $m \not\in v^dJ^{-1}$.   Thus $\delta_v(g) = 0$
if and only if $g \not\in \ker(\Theta_\tau)$.  We claim that if $\delta_v(g) > 0$,
then $\rho'$ arises from some $h$ with $\delta_v(h) = \delta_v(g) - 1$;
moreover if $r_m^J(g) = 0$ for all $mJ$ not prime to $\mathfrak{m}$,
then the same is true for $h$.  

To prove the claim, recall that if $g \in \ker(\Theta_\tau)$ and $k' = \nu(g)$,
then $p|k'_\tau$ for all $\tau \in \Sigma_v$.  (Recall that $\nu(g)$ is defined
in \S\ref{subsec:minimal} and that $\nu(g) \in \mathbb{Z}_{\ge 0}$ by \cite{DK}.)
Writing $g = g' \prod_{\tau' \in \Sigma} \mathrm{Ha}_{\tau'}^{n_\tau'}$ for some
$g' \in M_{k',0}(U_1(\mathfrak{n}');E)$ and $n' \in \mathbb{Z}^\Sigma_{\ge 0}$,
we have $g' \in \ker(\Theta_\tau)$.
By Theorem~\ref{thm:Phiv}, we have $g' = \Phi_v(g'')$ for some
$g'' \in M_{k'',0}(U_1(\mathfrak{n}');E)$ where $k''_{\tau'} = k'_{\tau'}$ for
$\tau' \not\in \Sigma_v$, and  $k''_{\tau'} = p^{-1}k'_{\mathrm{Fr}^{-1}\circ\tau'}$
for $\tau' \in \Sigma_v$.   Now 
$$h  := g'' \prod_{\tau' \in \Sigma_v} \mathrm{Ha}_{\tau'}^{k_\tau''}
\prod_{\tau' \in \Sigma} \mathrm{Ha}_{\tau'}^{n_\tau'}$$
is an eigenform in $M_{k,0}(U_1(\mathfrak{n}');E)$ giving rise to $\rho'$,
and Proposition~\ref{prop:Phivonqexps} immediately gives that 
$\delta_v(h) = \delta_v(g) - 1$, and if $r_m^J(g) = 0$
 for all $mJ$ not prime to $\mathfrak{m}$,
then the same is true for $h$.

Starting with $f_\xi$ and applying the claim inductively, we conclude that
$\rho'$ arises from an eigenform $g \in M_{k,0}(U_1(\mathfrak{n}');E)$
such that $r_m^J(g) \neq 0$ for some $m,J$ with $mJ$ prime to
$v\mathfrak{m}$.  Therefore  $g_{\xi^{-1}}$ is an 
eigenform in $M_{k,l}(U_1(\mathfrak{n}'\mathfrak{m}^2);E)$
giving rise to $\rho$, and $g_{\xi^{-1}} \not\in \ker(\Theta_\tau)$.

An elementary linear algebra argument then shows that, after
possibly further shrinking $\mathfrak{n}$ and enlarging $E$,
there is an eigenform $f$ which satisfies the conclusion
simultaneously for all $\tau \in \Sigma$.
\epf

We now have the following immediate consequences of Theorem~\ref{thm:theta}:
\begin{theorem} \label{thm:thetarho}
Suppose that $\rho$ is irreducible and geometrically modular of weight $(k,l)$. %
Then $\rho$ is geometrically modular of weight $(k',l')$, and in fact of weight $(k' - k_{\mathrm{Ha}_\tau},l')$
if $p|k_\tau$ (where $k'$ is as in Definition \ref{def:igusatheta}).
\end{theorem}

\begin{corollary}  \label{cor:thetakmin}
Suppose that $\rho$ is irreducible and $l ,l'\in \mathbb{Z}^\Sigma$ are such that 
$l'_\tau = l_\tau - \delta_{\tau,\tau'}$.
Suppose further that there exist $k = k_{\mathrm{min}}(\rho,l)$ and $k_{\mathrm{min}}(\rho,l')$ as in part 1) of
Conjecture~\ref{conj:geomweights}. Then
$$k_{\mathrm{min}}(\rho,l') \le_{\mathrm{Ha}}\,\,\,  \left\{\begin{array}{ll}  k',& \mbox{if $p\nmid k_\tau$} \\ 
k' - k_{\mathrm{Ha}_\tau} & \mbox{if $p|k_\tau$} \end{array}\right.$$
\end{corollary}

\begin{remark} We remark that we expect equality to hold in the corollary in the case that $p \nmid k_\tau$.
We caution however that the analogous strengthening of Theorem~\ref{thm:theta} is false:
i.e., it is possible for $\Theta_\tau(f)$ to be divisible by $\mathrm{Ha}_{\tau'}$ for some 
$\tau' \neq \tau$ even if $p \nmid k_\tau$ and $f$ is not divisible by $\mathrm{Ha}_{\tau'}$.\end{remark}

We also have:
\begin{corollary} 
Suppose that $\rho$ is irreducible and geometrically modular of some weight $(k_0,l_0)$.
Then for every $l \in \mathbb{Z}^\Sigma$, there exist $k \in \mathbb{Z}^\Sigma$ 
such that $\rho$ is geometrically modular of weight $(k,l)$.
\end{corollary}
\begpf Note that if $\rho$ is geometrically modular of some weight $(k_0,l_0)$,
then multiplying by the constant section $e_1$ of weight $(0,n(p-1))$, we can
replace $l_0$ by $l_0 + n(p-1)$ for any $n\in \mathbb{Z}$ and hence assume
$l_{0,\tau} \ge l_\tau$ for all $\tau \in \Sigma$.   The corollary
then follows from Theorem~\ref{thm:thetarho} by induction on
$\sum_\tau (l_{0,\tau} - l_\tau)$.
\epf

\subsection{Normalised eigenforms}

We continue to assume for simplicity that $J=O_F$ is chosen as an ideal class representative.

\begin{definition}
Suppose that $(k,l)$ is an algebraic weight (i.e., $k_\tau \ge 2$ for all $\tau \in \Sigma$).
We say that $f \in M_{k,l}(U_1(\mathfrak{n});E)$ is a {\em normalised eigenform} if the following hold:
\begin{itemize}
\item $r_1^{O_F}(f) = 1$, 
\item $f$ eigenform for  $T_v$ for all $v\nmid p$ and $S_v$ for all $v\nmid \mathfrak{n}p$, and
\item $f_\xi \in M_{k,0}(U_1(\mathfrak{nm}^2);E')$ is an eigenform for $T_v$ for all $v|p$ and all characters 
$\xi: \{\,a \in (\mathbb{A}_F^\infty)^\times \,|\, a_p \in O_{F,p}^\times\,\} / V_{\mathfrak{m}} \to (E')^\times$
of weight $-l$, conductor $\mathfrak{m}$ prime to $p$, and values in extensions $E'$ of $E$.
\end{itemize}

\end{definition}

It is straightforward to check that if $f_\xi$ is an eigenform for $T_v$ (where $v|p$ and $\xi$ has weight $-l$
and conductor prime to $p$), then so is $f_{\xi_1,\xi_2}$ for any characters $\xi_1,\xi_2$
such that $\xi_1\xi_2 = \xi$ (where the $\xi_i$ have conductors $\mathfrak{m}_i$ prime to $p$ and weights $-l_i$
such that $l=l_1+l_2$).  In particular it follows that if $f$ is a normalised eigenform in 
$M_{k,l}(U_1(\mathfrak{n});E)$, then $f_{\xi_1}$ is a normalised eigenform in 
$M_{k,l+l_1}(U_1(\mathfrak{nm}_1^2);E)$ (enlarging $E$ if necessary).

We have the following strengthening of Lemma~\ref{lem:U1better} for algebraic weights:
\begin{proposition} \label{prop:normalised} If $\rho$ is irreducible and geometrically modular of weight $(k,l)$ with
$k_\tau \ge 2$ for all $\tau$, then $\rho$ arises from a normalised eigenform of weight $(k,l)$ and level 
$U_1(\mathfrak{n})$ for some $\mathfrak{n}$ prime to $p$.
\end{proposition}
\begpf Suppose first that $l=0$.  By Lemma~\ref{lem:levelU1}, $\rho$ arises from an eigenform
$f \in M_{k,0}(U_1(\mathfrak{n});E)$ for some $\mathfrak{n}$ prime to $p$ (and some $E$).  Recall that in this
case ($l=0$ and all $k_\tau \ge 2$),
we have defined Hecke operators $T_v$ for all primes $v|\mathfrak{n}p$, commuting with each other
and the operators $T_v$ and $S_v$ for $v\nmid\mathfrak{n}p$, so we may further assume that $f$ is an
eigenform for $T_v$ for all $v$ and $S_v$ for all $v\nmid\mathfrak{n}p$.  
It suffices to prove that $r_1^{O_F}(f) \neq 0$.  

Suppose that $r_1^{O_F}(f) = 0$; we will show that $f= 0$, yielding
a contradiction.   Recall from Lemma~\ref{lem:Eisenstein} that the absolute irreducibility of $\rho \simeq \rho_f$ 
implies that $r_0^J(f) = 0$ for all $J$.  We will prove that $r_m^J(f) = 0$ for all $J$ and $m \in J^{-1}_+$
by induction on $n=\mathrm{Nm}_{F/\mathbb{Q}}(mJ)$.

If $n=1$, then $mJ = O_F$, so $J=O_F$, $m \in O_{F,+}^\times$, and $r_m^{O_F}(f) = r_1^{O_F}(f) = 0$.

Now suppose that $n > 1$ and $r_m^J(f) = 0$ for all $m,J$ with $\mathrm{Nm}_{F/\mathbb{Q}}(mJ)  < n$,
and let $m_1$, $J_1$ be such that $\mathrm{Nm}_{F/\mathbb{Q}}(m_1J_1)  = n$.  Let $v$ be any prime dividing
$m_1J_1$.  If $v^2\nmid m_1J_1$ or $v|\mathfrak{n}p$, then Propositions~\ref{prop:Tvonqexps},~\ref{prop:Uvonqexps}
and~\ref{prop:Uponqexps} give
$$r_{m_1}^{J_1}(f) =  r_m^J(T_vf) = a_v r_m^J(f),$$
where $m_1J_1 = vmJ$ and $a_v$ is the eigenvalue of $T_v$ on $f$.
We have $r_m^J(f) = 0$ by the induction hypothesis, and hence $r_{m_1}^{J_1}(f) = 0$.
If $v^2 | m_1J_1$ and $v\nmid\mathfrak{n}p$, then Proposition~\ref{prop:Tvonqexps} gives
$$r_{m_1}^{J_1}(f) =  r_m^J(T_vf)  - \mathrm{Nm}_{F/\mathbb{Q}}(v) r_{m_2}^{J_2}(S_vf)
= a_vr_m^J(f) - d_v\mathrm{Nm}_{F/\mathbb{Q}}(v)r_{m_2}^{J_2}(f),$$
where $m_1J_1 = vmJ = v^2m_2J_2$ and $a_v$ (resp.~$d_v$) is the eigenvalue of
$T_v$ (resp.~$S_v$) on $f$.  By the induction hypothesis, we have
$r_m^J(f) = r_{m_2}^{J_2}(f) = 0$, so again it follows that $r_{m_1}^{J_1}(f) = 0$.
This completes the proof of the proposition in the case $l=0$.

Now consider the case of arbitrary $l$.  Let $\mathfrak{m}$ (prime to $p$) be such that there is a character
$$\xi: \{\,a \in (\mathbb{A}_F^\infty)^\times \,|\, a_p \in O_{F,p}^\times\,\} / V_{\mathfrak{m}}
 \to E^\times$$
 of conductor $\mathfrak{m}$ satisfying $\xi(\alpha) = \overline{\alpha}^{-l}$ for all $\alpha \in F_+^\times
 \cap O_{F,p}^\times$.  Then $\rho_{\xi'}\otimes \rho$ is geometrically modular of weight $(k,0)$,
 and therefore arises from a normalised eigenform $f \in M_{k,0}(U_1(\mathfrak{n});E)$ for some
 $\mathfrak{n}$ prime to $p$.  Furthermore we may assume $\mathfrak{m}|\mathfrak{n}$ (for
 example by replacing $\mathfrak{n}$ by $\mathfrak{mn}$).  Then $f_{\xi^{-1}}$ is a normalised
 eigenform in $M_{k,l}(U_1(\mathfrak{mn});E)$ giving rise to $\rho$.
\epf

\subsection{Stabilised eigenforms}

We assume for the rest of the section that the weight $(k,l)$ is algebraic.

\begin{definition}
We say that a normalised eigenform  $f \in M_{k,l}(U_1(\mathfrak{n});E)$ is {\em stabilised}
if $r_m^J(f) = 0$ for all $(m,J)$ such that $m \in J_+^{-1}$ and $mJ$ is
not prime to $\mathfrak{n}$.   Note that this is equivalent to the condition that
$T_vf = 0$ for all $v|\mathfrak{n}$.\end{definition}

\begin{lemma} \label{lem:stable}  If $\rho$ arises from a normalised eigenform
in $M_{k,l}(U_1(\mathfrak{m});E)$, and $\mathfrak{n} \subset \mathfrak{m}$ is an ideal
prime to $p$, then $\rho$ arises from a normalised eigenform
in $M_{k,l}(U_1(\mathfrak{n});E)$ (enlarging $E$ if necessary).
Moreover if $\mathfrak{m}$ and $\mathfrak{n}$ satisfy
\begin{itemize}
\item $\mathrm{ord}_v(\mathfrak{nm}^{-1}) \ge 1$ for all $v|\mathfrak{m}$,
\item $\mathrm{ord}_v(\mathfrak{nm}^{-1}) \neq 1$ for all $v\nmid\mathfrak{m}$,
\end{itemize}
then $\rho$ arises from a stabilised eigenform in $M_{k,l}(U_1(\mathfrak{n});E)$.
\end{lemma}
\begpf The first assertion immediately reduces to the case $\mathfrak{n} = \mathfrak{m}v$
where $v$ is a prime not dividing $\mathfrak{m}p$.  Suppose that $f \in M_{k,l}(U_1(\mathfrak{m});E)$
is a normalised eigenform giving rise to $\rho$, and let $\alpha \in E$ (enlarging $E$ is necessary)
be a root of $X^2 - a_vX + d_v\mathrm{Nm}_{F/\mathbb{Q}}v$, the characteristic polynomial of
$\rho(\mathrm{Frob}_v)$, so $a_v$ (resp.~$d_v$) is the
eigenvalue of $T_v$ (resp.~$S_v$) on $f$.  A standard calculation then shows that
$$f' = f - (\mathrm{Nm}_{F/\mathbb{Q}}v)^{-1}\left(\begin{array}{cc}\varpi_v^{-1}&0\\0&1\end{array}\right)\alpha f$$
is a normalised eigenform in $M_{k,l}(U_1(\mathfrak{n});E)$.  Moreover $f'$ has the same eigenvalues as $f$,
except that $T_vf' = (a_v - \alpha)f'$. and therefore $\rho_{f'} \simeq \rho_f$.

In view of the first assertion, the second immediately reduces to the case 
$\mathfrak{n} = \mathfrak{m}\prod_{v|\mathfrak{m}} v$.  So suppose that $f \in M_{k,l}(U_1(\mathfrak{m});E)$
is a normalised eigenform giving rise to $\rho$, and for each $v|\mathfrak{m}$, let $\beta_v$ be the
eigenvalue of $T_v$ on $f$.  A similar standard calculation then shows that
$$f' = \,\,\, \prod_{v|\mathfrak{m}} \left( 1 - (\mathrm{Nm}_{F/\mathbb{Q}}v)^{-1}
\left(\begin{array}{cc}\varpi_v^{-1}&0\\0&1\end{array}\right)\beta_v\right)f.$$
is a normalised eigenform in $M_{k,l}(U_1(\mathfrak{n});E)$.  Moreover $f'$ has the same eigenvalues as $f$,
except that $T_vf' = 0$ for all $v|\mathfrak{n}$.  Therefore $f'$ is stabilised and gives rise to $\rho$.
\epf

\begin{remark} We remark that a more careful analysis easily shows that the first assertion of the lemma
requires at most a quadratic extension of $E$, and the second holds over the original field $E$.\end{remark}

\begin{definition}
We say that a stabilised eigenform  $f \in M_{k,l}(U_1(\mathfrak{n});E)$ is {\em strongly stabilised}
if $r_m^J(f) = 0$ for all $(m,J)$ such that $m \in J_+^{-1}\cup\{0\}$ and $mJ$ is
not prime to $p$.\footnote{Note that our conventions allow a stabilised eigenform
to have $r_0^J(f) \neq 0$ in the case $\mathfrak{n} = O_F$, but a strongly stabilised
eigenform necessarily has $r_0^J(f) = 0$.} \end{definition}
Thus a stabilised eigenform is strongly stabilised if and only if
$T_vf_\xi = 0$ for all $v|p$ and characters $\xi$ of weight $-l$.  (Note that given $m \in J_+^{-1} \cup \{0\}$,
we can always choose $\xi$ of weight $-l$ and conductor prime to $mJ$ unless $m = 0$
and $\overline{\nu}^l \neq 1$ for some $\nu \in O_{F,+}^\times$, in which case we automatically
have $r_m^J(f) = 0$.)

\begin{lemma}  \label{lem:strong} There is at most one strongly stabilised eigenform $f \in M_{k,l}(U_1(\mathfrak{n});E)$
giving rise to $\rho$.
\end{lemma}
\begpf If $\rho$ arises from $f$, then $T_vf = a_vf$ and $S_vf = d_vf$ for all $v\nmid\mathfrak{n}p$,
where $a_v = \mathrm{tr}(\rho(\mathrm{Frob}_v))$ and  $d_v = \mathrm{Nm}_{F/\mathbb{Q}}(v)^{-1}\mathrm{det}(\rho(\mathrm{Frob}_v))$.  

Suppose then that $f$ and $f'$ are strongly stabilised eigenforms giving rise to $\rho$, and let $f'' = f-f'$.
It suffices to prove that $r_m^J(f'') = 0$ for all $(m,J)$ with $m \in J_+^{-1} \cup \{0\}$.
Since $f$ and $f'$ are strongly stabilised, we have $r_m^J(f'') =0$ whenever
 $mJ$ is not prime to $\mathfrak{n}p$, so we can assume $mJ$ is prime to $\mathfrak{n}p$.
We then proceed as in the proof of Proposition~\ref{prop:normalised} by induction on 
$n = \mathrm{Nm}_{F/\mathbb{Q}}(mJ)$.

If $n =1$, then $mJ = O_F$, so $J=O_F$, $m \in O_{F,+}^\times$, and $r_m^{O_F}(f'') = m^{-l}r_1^{O_F}(f'') = 0$
since $r_1^{O_F}(f) = r_1^{O_F}(f') = 1$.

Now suppose that $n > 1$ and $r_m^J(f'') = 0$ for all $m,J$ with $\mathrm{Nm}_{F/\mathbb{Q}}(mJ)  < n$,
and let $m_1$, $J_1$ be such that $m_1J_1$ is prime to $\mathfrak{n}p$ and $\mathrm{Nm}_{F/\mathbb{Q}}(m_1J_1)  = n$. 
Let $v$ be any prime dividing $m_1J_1$.  If $v^2\nmid m_1J_1$, then Proposition~\ref{prop:Tvonqexps} gives
$$r_{m_1}^{J_1}(f'') =  m_1^{-l} m^l a_v r_m^J(f'')$$
where $m_1J_1 = vmJ$. so the induction hypothesis implies that $r_{m_1}^{J_1}(f'') = 0$.
If $v\nmid\mathfrak{n}p$, then we get instead
$$r_{m_1}^{J_1}(f'') = m_1^{-l}m^la_vr_m^J(f'') - m_1^{-l}m_2^ld_v\mathrm{Nm}_{F/\mathbb{Q}}(v)r_{m_2}^{J_2}(f''),$$
where $m_1J_1 = vmJ = v^2m_2J_2$, and again the induction hypothesis implies that $r_{m_1}^{J_1}(f'') = 0$.
\epf

\begin{remark} Note that if $f$ is a normalised (resp.~stabilised, strongly stabilised) eigenform, then 
the same is true for both $\mathrm{Ha}_\tau f$ and $\Theta_\tau f$
for any $\tau$ (assuming $k_\tau \ge 3$ if $\tau \neq \mathrm{Fr}\circ\tau$
in the case of $\mathrm{Ha}_\tau f$).\end{remark}

\begin{remark} We remark that if $\rho$ is geometrically modular of weight $(k,l)$, then it does not necessarily
arise from a strongly stabilised eigenform of weight $(k,l)$ (for any level $\mathfrak{n}$); for example,
there may be a prime $v|p$ such that $r_m^J(f) \neq 0$ whenever $mJ = v$.  We do however
have the following two ways of establishing the existence of strongly stabilised eigenforms.
One is to apply partial $\Theta$-operators to a stabilised eigenform (hence changing the weight);
the other is to use Theorem~\ref{thm:ordinary} below, or more precisely its corollary. \end{remark}

\subsection{Ordinariness}
The forthcoming Theorem~\ref{thm:ordinary} can be viewed as stating that if an eigenform is ordinary in
a suitable sense, then so is the associated Galois representation.
For the proof, we need to verify certain
compatibility properties for the operators $T_v$ for $v|p$ (assuming $k$ algebraic
and $l=0$), which we shall do using their
effect on $q$-expansions at more general cusps than the ones used above.

Fix sets of ideal class representatives $\{\mathfrak{a}\}$
and coset representatives $\{g\}$ for $P_{\mathfrak{n}}\backslash\mathrm{SL}_2(O_F/\mathfrak{n})$.
For consistency with previous computations, choose $\mathfrak{a} = O_F$ and $g=1$ for the trivial classes.
Also fix choices of 
$t:O_F/NO_F \simeq\mu_N  \otimes (\mathfrak{ad})^{-1} $ for each $\mathfrak{a}$ (independent of $J$)
and $s:O_F/NO_F \simeq (NJ)^{-1}\mathfrak{a}/J^{-1}\mathfrak{a}$ for each $\mathfrak{a}, J$.
The cusps of $X_U$ for $U=U(\mathfrak{n})$ are then in bijection with the quadruples $(J,\mathfrak{a},w,g)$,
where the corresponding cusp is the one associated to the Tate semi-HBAV $T_{\mathfrak{a},\mathfrak{b}}$,
where $\mathfrak{b} = \mathfrak{a}J^{-1}$, with canonical polarisation and level structure
$\eta_w\circ r_{g^{-1}}$, where $\eta_w(x,y) = t(y)q^{ws(x)}$.
Then for $v|p$, we find (by the same proof as for Proposition~\ref{prop:Uponqexps})
that the effect on $q$-expansions of the action of $T_v$ on $M_{k,0}(U;E)$ is given by
\begin{equation}\label{eqn:UpUn}  r_m^C(T_vf) = r_{\beta_1m}^{C_1}(f),\end{equation}
where the $q$-expansion coefficients lie in $\overline{D}^{k,0}$, and
if $C$ is the cusp corresponding to $(J,\mathfrak{a},w,g)$,
then $C_1$ corresponds to $(J_1,\mathfrak{a},w_1,g)$ for $J_1,w_1,\beta_1\in F_+$
such that $vJ = \beta_1J_1$ and  $\beta_1\varpi_v^{-1}ws(1) \equiv w_1s_1(1) \bmod N^{-1}\mathfrak{nb}_1$
(where $\mathfrak{b}_1 = \mathfrak{a}J_1^{-1}$ and $s_1$ is the chosen isomorphism).

Finally we need to consider the action of $T_v$ on $M_{k,0}(U';L)$ for $v|p$, where $U' = U \cap U_1(p)$
and $k \in \mathbb{Z}_{\ge 2}$.   Note that this may be defined in the usual way as the operator 
$\left[U'\left(\begin{array}{cc}\varpi_v&0\\0&1\end{array}\right)U'\right]$ on forms in characteristic zero,
making it compatible with the action of $T_v$ on the space of automorphic forms $A_{k,0}(U')$.
Recall from \S\ref{subsec:XU} that $X_{U'}$ denotes the minimal compactification of
$Y_{U'}$, and that its cusps are in bijection with triples $(C,\mathfrak{f},P)$ where $C$ is
a cusp of $X_U$, $\mathfrak{f}$ is an ideal such that $pO_F \subset \mathfrak{f} \subset O_F$
and $P$ is a generator of $\mathfrak{b}/\mathfrak{bf}$, and the corresponding cusp may
be identified with the $\mathcal{O}$-scheme $\mathrm{Spec}\,\mathcal{O}'_{\mathfrak{f}}$ representing
generators of $\mu_p \otimes \mathfrak{f}(\mathfrak{ad})^{-1}/p(\mathfrak{ad})^{-1}$
(where $\mathfrak{a}$ and $\mathfrak{b}$ are as in the description of $C$).
We only need to consider those cusps for which $\mathfrak{f} = O_F$: for each cusp $C$ of $X_U$,
we write $C'$ for the unique such cusp of $X_{U'}$ lying over it.
We assume  $L$ contains the $p$th roots of unity, so that the components
of $C'_L$ are copies of $\mathrm{Spec}\,L$ in bijection with the generators $\zeta_p$ of
$\mu_p(L)\otimes  (\mathfrak{ad})^{-1}$.
We may then compute the effect of $T_v$ on the completion at each component of $C'_L$
exactly as in Proposition~\ref{prop:Uvonqexps} (see also (\ref{eqn:UvUn})) to conclude that if
$f\in M_{k,0}(U';L)$, then 
\begin{equation} \label{eqn:UpUn2} r_m^{C'}(T_vf) = r_{\beta_1m}^{C'_1}(f),\end{equation}
where the notation is as in  (\ref{eqn:UpUn}), except that the $q$-expansion coefficients lie in
the fibre of $j'_*\mathcal{L}_{U'}^{k,0}$ at $C'_L$, which we may identify with
$\oplus_{\zeta_p} (D^{k,0}\otimes_{\mathcal{O}} L)$ (where $\zeta_p$ runs over generators
of $\mu_p(L)\otimes  (\mathfrak{ad})^{-1}$).

\begin{theorem} \label{thm:ordinary} Suppose that $k\in \mathbb{Z}^{\Sigma}$ with $k_\tau \ge 2$ for all $\tau$,
$U$ is an open compact subgroup of $\mathrm{GL}_2(\widehat{O}_F)$
containing $\mathrm{GL}_2(O_{F,p})$, $Q$ is a finite set of primes containing all $v|p$ and
all $v$ such that $\mathrm{GL}_2(O_{F,v}) \not\subset U$, and $v_0$ is a prime over $p$. 
Suppose that $f \in M_{k,0}(U;E)$ is an eigenform for $T_v$ and $S_v$ for all $v \not\in Q$ and 
that $T_{v_0}f = a_{v_0}f$ for some $a_{v_0} \neq 0$.
Then (possibly after enlarging $E$ and semi-simplifying $\rho_f$)
$$\rho_f|_{G_{F_{v_0}}} \simeq \left(\begin{array}{cc}  \chi_1 & * \\ 0 & \chi_2 \end{array}\right)$$
where $\chi_1$ unramified character, $\chi_1(\mathrm{Frob}_{v_0}) = a_{v_0}$,
and $\chi_2|_{I_{F_{v_0}}} = \prod_{\tau\in\Sigma_{v_0}} \epsilon_\tau^{1-k_\tau}$.
\end{theorem}
\begpf  We may assume that $U= U(\mathfrak{n})$ for some sufficiently small $\mathfrak{n}$
prime to $p$ and that $\mathcal{O}$ is sufficiently large; in particular, we assume 
$\mu_{Np}(\overline{\mathbb{Q}}) \subset \mathcal{O}$ for some $N \in \mathfrak{n}$.

Recall that the proof of Theorem~\ref{thm:galois} in \S\ref{subsec:proof} yields injections
$$M_{k,-1}(U;E) \to M_{k',-1}(U':E) \to M_{m+2,-1}(U';\mathcal{O}) \otimes_{\mathcal{O}} E$$
which are compatible with $T_v$ and $S_v$ for $v\nmid \mathfrak{n}p$, where $U=U(\mathfrak{n})$,
$U' = U(\mathfrak{n})\cap U_1(p)$, $k'$ is nearly parallel and $m$ is a (sufficiently large)
positive integer.   Tensoring with the (pull-backs of the) canonical section 
$e_1 \in H^0(Y_{U}, \mathcal{L}_{U}^{0,1})$, we may replace
$l=-1$ by $l=0$.   Since the first injection is defined by multiplication by partial Hasse invariants,
which have $q$-expansions equal to $1$ at every cusp, we see from (\ref{eqn:UpUn}) that
it is also compatible with $T_v$ for $v|p$.  We may therefore replace $k$ by $k'$ and assume
that $k$ is nearly parallel.

Recall that for a cusp $C$ of $X_U$, we let $C'$ denote the unique cusp of $X_{U'}$ with
$\mathfrak{f} = O_F$.    For $\mathcal{O}$-algebras $R$, let $Q_{C',R}^{m+2,0}$
denote the completion at $C'_R$ of 
$j'_*(\mathcal{K}_{U',R}\otimes_{\mathcal{O}_{Y_{U',R}}} \mathcal{L}_{U',R}^{m,1})$
(in the notation of \S\ref{sec:galois}, and as usual omitting subscripts if $R =\mathcal{O}$
and using $\overline{\,\cdot\,}$ in the case $R=E$).  From the description of
$j'_*\mathcal{K}_{U'}$ in \S\ref{subsec:proof}, we see that
$Q_{C',R}^{m+2,0}$ is canonically isomorphic to
$$D^{m+2,0} \otimes_{\mathcal{O}} \mathrm{Hom}_{\widehat{S}_{C,R}}(\widehat{S}_{C',R},\widehat{S}_{C,R})
    \cong \mathrm{Hom}_{\mathcal{O}}( \mathcal{O}'_{O_F}, D^{m+2,0} ) 
\otimes_{\mathcal{O}} \widehat{S}_{C,R}$$
as a module over $\widehat{S}_{C',R} \cong \mathcal{O}'_{O_F} \otimes_{\mathcal{O}} \widehat{S}_{C,R}$.

Letting $\mathcal{S}$ denote the set of cusps of $X_U$,
we have natural $q$-expansion maps:
$$M_{m+2,0}(U';R) \to \bigoplus_{C\in \mathcal{S}}  Q_{C',R}^{m+2,0},$$
which are injective if $R = L$ (and hence $R=\mathcal{O}$) 
since $\displaystyle\coprod_{C \in \mathcal{S}} C'_L$ includes cusps on every
connected component of $X_{U',L}$.  We define 
$\widetilde{M}_{m+2,0}(U';\mathcal{O})$ to be the preimage of $\displaystyle\bigoplus_{C\in \mathcal{S}}  Q_{C'}^{m+2,0}$
in $M_{m+2,0}(U';L)$ under the $q$-expansion map to $\displaystyle\bigoplus_{C\in \mathcal{S}}  Q_{C',L}^{m+2,0}$.
We thus have an inclusion
$M_{m+2,0}(U';\mathcal{O}) \subset \widetilde{M}_{m+2,0}(U';\mathcal{O})$
 with finite index, so in particular
$\widetilde{M}_{m+2,0}(U';\mathcal{O})$ is finitely generated over $\mathcal{O}$.

We see directly from the definition that, for $v\nmid p\mathfrak{n}$,  the Hecke operators
$T_v$ and $S_v$ on $M_{m+2,0}(U';\mathcal{O})$ also act on the modules
$Q_{C'}^{m+2,0}$ compatibly with the $q$-expansion map, from which it follows that
the operators preserve $\widetilde{M}_{m+2,0}(U';\mathcal{O})$.  
Furthermore, from (\ref{eqn:UpUn2}) and the fact that the isomorphism
$$Q_{C',L}^{m+2,0} \cong D^{m+2,0}  \otimes_{\mathcal{O}} \widehat{S}_{C',L}$$
induced by the Kodaira--Spencer isomorphism $\mathcal{K}_{U',L} \cong \mathcal{L}_{U',L}^{2,-1}$
is the same as the one induced by the canonical isomorphisms
$$\mathrm{Hom}_{\mathcal{O}}(\mathcal{O}_{O_F}',L)
 \,\cong \,\oplus_{\zeta_p} L\,\cong\, \mathcal{O}_{O_F}'\otimes_{\mathcal{O}}L,$$
we see that  $\widetilde{M}_{m+2,0}(U';\mathcal{O})$
is also stable under $T_v$ for $v|p$, with the action on $q$-expansions being defined
by the same formula, but now with coefficients in the fibre
$\mathrm{Hom}_{\mathcal{O}}( \mathcal{O}'_{O_F}, D^{m+2,0} )$
(and a reconciliation of the duplicate use of $m$).

Now consider the commutative diagram:
\begin{equation}\label{eqn:qexpcd}  \begin{array}{ccccc}
M_{k,0}(U;E) &  \to& M_{m+2,0}(U';\mathcal{O})\otimes_{\mathcal{O}}E 
&  \to& \widetilde{M}_{m+2,0}(U';\mathcal{O})\otimes_{\mathcal{O}}E \\
\downarrow&&&&\downarrow \\
\displaystyle\bigoplus_{C\in \mathcal{S}} Q_{\overline{C}}^{k,0} & \to &
 \displaystyle\bigoplus_{C\in \mathcal{S}} Q_{\overline{C}'}^{m+2,0} &\leftarrow  &
 \displaystyle\bigoplus_{C\in \mathcal{S}} Q_{C'}^{m+2,0} \otimes_{\mathcal{O}}E.\end{array}\end{equation}
The proof of Theorem~\ref{thm:galois} shows that the first map on the bottom
row is injective and the second is an isomorphism.  Since the left vertical
arrow is injective by the $q$-expansion principle, it follows that the top
composite is also injective.  Furthermore the right vertical arrow is injective
since $\widetilde{M}_{m+2,0}(U';\mathcal{O}) \to  \displaystyle\bigoplus_{C\in \mathcal{S}} Q_{C'}^{m+2,0}$
is injective with torsion-free cokernel (by construction).  

We already saw in the proof of Theorem~\ref{thm:galois} that the map
$M_{k,0}(U;E)   \to M_{m+2,0}(U';\mathcal{O})\otimes_{\mathcal{O}}E$ is compatible
with the operators $T_v$ and $S_v$ for $v\nmid p\mathfrak{n}$, and it follows that
the same holds for the composite on the top row of (\ref{eqn:qexpcd}).  We claim that 
this composite is also
compatible with the operators $T_v$ for $v|p$.  To see this, note that we have
actions of these operators on $\bigoplus_{C\in \mathcal{S}} Q_{\overline{C}}^{k,0}$
and  $\bigoplus_{C\in \mathcal{S}} Q_{\overline{C}'}^{m+2,0}
\cong \bigoplus_{C\in \mathcal{S}} Q_{C'}^{m+2,0} \otimes_{\mathcal{O}}E$ which are compatible
with the vertical maps of  (\ref{eqn:qexpcd}).  Since these maps are injective, the
claim will follow from the compatibility of the bottom row of (\ref{eqn:qexpcd}) with
these operators.  The desired compatibility then follows from the fact
that the first map on the bottom row is induced by the $\mathcal{O}_S$-dual of the
pull-back to $S = \mathrm{Spec}\,\widehat{S}_{\overline{C}} - \overline{C}$
of the isomorphism
$\overline{\pi}_*i_*\mathcal{O}_{Y^\mu_U} \cong  \oplus_{\kappa}   \overline{\mathcal{L}}_{U}^{\kappa,0}$
constructed in the proof of Theorem~\ref{thm:galois}, which is given with respect to the canonical
trivialisations over $S$ by an isomorphism
$$\mathcal{O}'_{O_F} \otimes_{\mathcal{O}} E \, \cong \, \bigoplus_{\kappa}  \overline{D}^{\kappa,0}.$$

We have now shown that the top row of (\ref{eqn:qexpcd}) defines an injective homomorphism
$$M_{k,0}(U;E)   \to \widetilde{M}_{m+2,0}(U';\mathcal{O})\otimes_{\mathcal{O}}E,$$
compatible with the operators $T_v$ for $v\nmid \mathfrak{n}$ and $S_v$ for $v\nmid\mathfrak{n}p$.
It is therefore a homomorphism 
of $\mathbb{T}$-modules, where $\mathbb{T}$ is the (commutative) $\mathcal{O}$-algebra of
endomorphisms of $\widetilde{M}_{m+2,0}(U';\mathcal{O})$ generated by $T_{v_0}$ and
the operators $T_v$ and $S_v$ for $v\not\in Q$.
The same (standard) argument as at the end of the proof of Theorem~\ref{thm:galois} now shows that
(after enlarging $\mathcal{O}$, $L$ and $E$ if necessary), there
is an eigenform $\tilde{f} \in M_{m+2,0}(U';L)$  for the operators $T \in \mathbb{T}$
such that the eigenvalues are lifts of the corresponding ones for $f$.
In particular $\overline{\rho}_{\tilde{f}}$ and $\rho_f$ have isomorphic semi-simplifications,
and $T_{v_0} \tilde{f} = \tilde{a}_{v_0} \tilde{f}$ for some $\tilde{a}_{v_0} \in \mathcal{O}^\times$.

We now deduce that $\rho_f|_{G_{F_{v_0}}}$ has the desired form from the analogous
fact for the characteristic zero modular Galois representation $\rho_{\tilde{f}}$, which is a special case of local-global compatibility at $v_0$ for
the corresponding automorphic and Galois representations.  More precisely, suppose first
that $\tilde{f}$ is cuspidal and view it as a vector fixed by $U'$ in the associated automorphic
representation $\Pi$, so we have that $\tilde{a}_{v_0}$ is an eigenvalue for $T_{v_0}$ on
$\Pi_{v_0}^{U'_{v_0}}$, where $\Pi_{v_0}$ is the local factor of $\Pi$ at $v_0$ and
$U'_{v_0} = U_1(v_0) \cap \mathrm{GL}_2(O_{F_{v_0}})$.  We may assume for simplicity
that $m > 0$, so that $\Pi_{v_0}$ must be a principal series representation\footnote{Permitting $m=0$
would allow the possibility that $\Pi_{v_0}$ be an unramified twist of the Steinberg representation, which
could anyway have been treated similarly.} of the
form $I(\psi_1|\cdot|^{1/2},\psi_2|\cdot|^{1/2})$ where $\psi_1,\psi_2$ are characters
$F_{v_0}^\times \to \overline{\mathbb{Q}}^\times$ such that $\psi_1$ is unramified with
$\psi_1(\varpi_{v_0}) = \tilde{a}_{v_0}$ (and $\psi_2$ is at most tamely ramified with
$\psi_2(\varpi_{v_0})(\mathrm{Nm}_{F/\mathbb{Q}}(v_0))^{-m-1} \in \mathcal{O}^\times$).
The main theorem of \cite{Sk} (adapted to our conventions) then implies that
$\rho_{\tilde{f}}|_{G_{F_{v_0}}}$ is potentially crystalline with labelled Hodge--Tate
weights $(m+1,0)$ and associated Weil--Deligne representation $\psi_1 \oplus \psi_2$
(writing $\psi_i$ also for the representations of $W_{F_{v_0}}$ to which they correspond
by local class field theory).  A standard exercise in $p$-adic Hodge theory then shows
that $\rho_{\tilde{f}}|_{G_{F_{v_0}}}$ must be of the form:
$$\left(\begin{array}{cc}  \widetilde{\chi}_1 & * \\ 0 & \widetilde{\chi}_2 \end{array}\right)$$
for some $\widetilde{\chi}_1,\widetilde{\chi}_2 : G_{F_{v_0}} \to L^\times$ with
$\widetilde{\chi}_1$  unramified and $\widetilde{\chi}_1(\mathrm{Frob}_{v_0}) = \tilde{a}_{v_0}$ (and
$\widetilde{\chi}_2\chi_{\mathrm{cyc}}^{m+1}$ at most tamely ramified).  The theorem then
follows in this case from the fact that $\rho_f$ is (up to semi-simplification) the reduction
mod $\pi$ of $\rho_{\tilde{f}}$, together with the description of $\det(\rho_f)$ in
Remark~\ref{rmk:det}.

Suppose on the other hand that $\tilde{f}$ is not cuspidal, in which case its eigenvalues
for $T_v$ and $S_v$ for $v\not\in Q$ are the same as those of an Eisenstein series
associated to a pair of Hecke characters $\psi_1,\psi_2$ such that $\psi_1(x) = 1$
and $\psi_2(x) = x^{-m-1}$ for $x \in F_{\infty,+}^\times$.  Moreover
$\tilde{a}_{v_0} = \psi_i(\varpi_{v_0})$ for some $i$ such that $\psi_i$ is
unramified at $v_0$, and we must have $i=1$ since  $\tilde{a}_{v_0} \in \mathcal{O}^\times$.
  In this case the (semi-simplification of the) associated Galois representation
$\rho_{\tilde{f}}$ is $\widetilde{\chi}_1 \oplus \widetilde{\chi}_2$, where $\widetilde{\chi}_1$
(resp.~$\widetilde{\chi}_2\chi_{\mathrm{cyc}}^{m+1}$) is the character associated 
to $\psi_1$ (resp.~$\psi_2|\cdot|^{m+1}$) by class field theory.  The theorem thus follows
as before on reduction mod $\pi$.  \epf

\begin{corollary} \label{cor:ordinary} Let $f \in M_{k,l}(U_1(\mathfrak{n});E)$ 
be a normalised eigenform and ${v_0}$ a prime of $F$ over $p$.  If $T_{v_0}f_\xi \neq 0$
for some character $\xi$ of weight $-l$, then (possibly after enlarging $E$ and semi-simplifying $\rho_f$)
$$\rho_f|_{G_{F_{v_0}}} \simeq \left(\begin{array}{cc}  \chi_1 & * \\ 0 & \chi_2 \end{array}\right)$$
for some characters $\chi_1,\chi_2 : G_{F_{v_0}} \to E^\times$ such that
$\chi_1|_{I_{F_{v_0}}} = \prod_{\tau\in\Sigma_{v_0}} \epsilon_\tau^{-l_\tau}$ and
$\chi_2|_{I_{F_{v_0}}} = \prod_{\tau\in\Sigma_{v_0}} \epsilon_\tau^{1-k_\tau-l_\tau}$.
\end{corollary}
\begpf  Since $\rho_{f_\xi} \simeq \rho_f \otimes \rho_{\xi'}$ and 
$\rho_{\xi'}|_{I_{F_{v_0}}} = \prod_{\tau\in\Sigma_{v_0}} \epsilon_\tau^{l_\tau}$,
we may reduce to the case $l=0$ and $f=f_\xi$, which is
immediate from Theorem~\ref{thm:ordinary}.  \epf

\section{The inert quadratic case}

\label{sec:quadratic}
We now specialise to the  inert quadratic case,  with a focus on non-algebraic weights,
and in particular the case of ``partial weight one''  since it exhibits phenomena
not present in the classical case.  We provide evidence and an approach to
Conjectures~\ref{conj:geomweights} and~\ref{conj:geomweights2}
by deducing results in this setting  from ones in the case of algebraic weights.

\subsection{Notation}

For the rest of the paper, we let $F$ be a real quadratic field in which $p$ is inert,
and we let $\mathfrak{p} = pO_F$ and $K = F_{\mathfrak{p}}$, so $K$ is the unramified
quadratic extension of $\mathbb{Q}_p$.   Fix an embedding $\tau_0:F \to \overline{\mathbb{Q}}_p$
and write $\Sigma = \{\tau_0,\tau_1\}$.  We identify $\Sigma$ with $\Sigma_K$ and hence
with the set of embeddings $O_F/pO_F \to \overline{\mathbb{F}}_p$.
We shall write weights $k \in \mathbb{Z}^{\Sigma}$
in the form $(k_0,k_1)$ where $k_i = k_{\tau_i}$ for $i=0,1$.
Recall that our conventions for Hodge--Tate types and weights of crystalline lifts
of two-dimensional representations are given in \S\ref{subsec:cryslift}.

\subsection{$p$-adic Hodge theory lemmas}
Let $\chi:G_K \to \overline{\mathbb{F}}_p^\times$ be a character such that 
$\chi|_{I_K} = \epsilon_{\tau_0}^i$ with $1 \le i \le p-1$.  Then
 $H^1(G_K,\overline{\mathbb{F}}_p(\chi))$ is two-dimensional, and
we recall  from \cite[\S3]{BDJ} the definition of a certain one-dimensional subspace.
Note that $\chi|_{I_K} = \epsilon_{\tau_0}^{i-1}\epsilon_{\tau_1}^p$, so $\chi$ has a
crystalline lift $\widetilde{\chi}:G_K \to \mathcal{O}^\times$ with Hodge--Tate
type $(1-i,-p) \in  \mathbb{Z}^{\Sigma}$ (where $\mathcal{O}$ is assumed to be
sufficiently large that $\chi$ takes values in $E^\times$).  Such lifts are unique up to twist by unramified characters
with trivial reduction, and we choose\footnote{By Remark~7.13 of \cite{CD}, or more generally the proof of Theorem~9.1 of \cite{GLS},
the subspace $V_\chi$ turns out to be independent of the choice of unramified twist, but we fix it for clarity
and consistency with \cite{BDJ}.  Similarly the proof of Lemma~\ref{lem:lines} below shows the same
holds for $V'_\chi$.} the one such that if $g$ corresponds via local class field theory to $p \in K^\times$,
then $\widetilde{\chi}(g)$ is the Teichm\"uller lift of $\chi(g)$.  A standard argument
shows that the space $H^1_f(G_K,L(\widetilde{\chi}))$ classifying crystalline extensions
is one-dimensional over $L$, with preimage $V_{\widetilde{\chi}} \subset 
H^1(G_K,\mathcal{O}(\widetilde{\chi}))$ free of rank one over $\mathcal{O}$.  We then define
$V_{\chi} =  V_{\widetilde{\chi}} \otimes_{\mathcal{O}} \overline{\mathbb{F}}_p$.
Similarly, $\chi$ has a crystalline lift $\widetilde{\chi}':G_K \to \mathcal{O}^\times$ with Hodge--Tate
type $(-i,0)$, unique up to unramified twist, and
we choose the one sending $g$ (corresponding to $p$) to the Teichm\"uller lift of $\chi(g)$.
We again have that $H^1_f(G_K,L(\widetilde{\chi}'))$ is one-dimensional, with preimage $V_{\widetilde{\chi}'} \subset 
H^1(G_K,\mathcal{O}(\widetilde{\chi}))$ free of rank one over $\mathcal{O}$, and we define
$V_{\chi}' =  V_{\widetilde{\chi}'} \otimes_{\mathcal{O}} \overline{\mathbb{F}}_p$.

\begin{lemma} \label{lem:lines}  With the above notation, $V_\chi = V'_\chi$.
\end{lemma}
\begpf We use the description of $V_\chi$ obtained in \cite{CD} together with
a similar analysis of $V'_\chi$.  All references in this proof are to \cite{CD}.

In the notation of \cite{CD}, the $(\phi,\Gamma)$-module, corresponding to the one-dimensional $E$-vector space $E(\chi)$ equipped with $G_K$ action by $\chi$, has the form $M_{C\vec{c}}$
with $\vec{c} = (p-1-i,p-1)$, and Proposition~5.11 (for $p>2$), Proposition~6.11 (for $p=2$) and Theorem~7.12
show that $V_\chi$ is the subspace corresponding to the span of the class $[B_0] \in \mathrm{Ext}^1(M_{\vec{0}},M_{C\vec{c}})$.

We may analyze $V_\chi'$ similarly as follows.  We can write the $(\phi,\Gamma)$-module corresponding
to $E(\chi^{-1})$ in the form $M_{A\vec{a}}$ where $\vec{a} = (i,0)$ and $A=C^{-1}$, and consider the
subspace of bounded extensions 
$$\mathrm{Ext}_{\mathrm{bdd}}^1(M_{A\vec{a}},M_{\vec{0}}) \subset \mathrm{Ext}^1(M_{A\vec{a}},M_{\vec{0}})$$
defined exactly in Definition~5.1 (dropping the assumption that
one of $a_i$ or $b_i$ is non-zero for each $i$).   As in \S5.1, we have an isomorphism
$$\iota: \mathrm{Ext}^1(M_{A\vec{a}},M_{\vec{0}}) \cong  \mathrm{Ext}^1(M_{\vec{0}},M_{C\vec{c}}).$$
A straightforward adaptation of part (2) of the proof of Proposition~5.11\footnote{Strictly speaking, this is Proposition~6.11 in
the case $p=2$, but the proof there is omitted since it is essentially the same as that of Proposition~5.11, using
the cocycles constructed in \S6.3.} then shows that the image
of $\mathrm{Ext}_{\mathrm{bdd}}^1(M_{A\vec{a}},M_{\vec{0}})$ under $\iota$ is again spanned by $[B_0]$.

By the same argument as in the proof of Theorem~7.8
(with the appeal to Lemma~7.6 replaced by a direct application of Proposition~7.4), one finds that
$V_\chi'$ is contained in (the extension of scalars to $\overline{\mathbb{F}}_p$ of) the image of 
$\mathrm{Ext}_{\mathrm{bdd}}^1(M_{A\vec{a}},M_{\vec{0}})$.  Therefore $V_\chi' \subset V_\chi$,
and equality follows on comparing dimensions.
\epf

\begin{remark}  In \cite[Chapter 5]{HWPhD}, Wiersema gives an alternative proof of the preceding lemma
by first generalising methods and results of \cite{GLS, GLS2} to the setting of two-dimensional crystalline
representations of $G_K$ whose $\tau$-labelled weights $w_{1,\tau}$, $w_{2,\tau}$ need not be distinct.
One can then appeal to those results instead of the ones in \cite{CD} (see \cite[Lemma~5.3.2]{HWPhD}).
Similarly the following lemma can be proved using those results instead of Fontaine--Laffaille theory
(see \cite[Lemma~5.4.7]{HWPhD}).
\end{remark}

\begin{lemma}  \label{lem:pwt1explicit} Suppose that $2 \le k_0 \le p$.  A
representation $\sigma:G_K \to \mathrm{GL}_2(\overline{\mathbb{F}}_p)$
has a crystalline lift of weight $((k_0,1),(0,0))$ if and only if either:
\begin{itemize}
\item $\sigma \simeq \left(\begin{array}{cc}  \chi_1 & * \\ 0 & \chi_2 \end{array}\right)$
with $\chi_1$ unramified, $\chi_2|_{I_K} = \epsilon_{\tau_0}^{1-k_0}$ and associated
extension class in $V_{\chi_1\chi_2^{-1}}$, or
\item  $\sigma\simeq \mathrm{Ind}_{G_{K'}}^{G_K} \xi$ where $K'$ is the unramified quadratic
extension of $K$ and $\xi|_{I_{K'}} = \epsilon_{\tau_0'}^{1-k_0}$ for some extension $\tau_0'$ of
$\tau_0$ to the residue field of $K'$.
\end{itemize}
\end{lemma}
\begpf  For the ``if'' direction, in the first case, let $\widetilde{\chi}_1$ be an unramified
lift of $\chi_1$ and let $\widetilde{\chi}_2 = \widetilde{\chi}_1(\widetilde{\chi}')^{-1}$, where $\chi = \chi_1\chi_2^{-1}$.
By Lemma~\ref{lem:lines} and the definition of $V'_\chi$, the representation $\chi_2^{-1}\otimes\sigma$ is isomorphic to
the reduction of an $\mathcal{O}[G_K]$-module $T$ associated to an extension class
$$0 \to \mathcal{O}(\widetilde{\chi}') \to  T \to \mathcal{O} \to 0$$
such that $T\otimes_{\mathcal{O}}L$ is crystalline.  It follows that $\sigma$ has a crystalline
lift $T\otimes_{\mathcal{O}}L(\widetilde{\chi}_2)$ with $\tau_0$-labelled weights $(k_0-1,0)$
and $\tau_1$-labelled weights $(0,0)$, as required.

In the second case, note that $\xi$ has a crystalline lift $\widetilde{\xi}$ of Hodge--Tate
type $(k_0-1,0,0,0)$ (where the first coordinate corresponds to $\tau'_0 \in \Sigma_{K'}$),
so that $ \mathrm{Ind}_{G_{K'}}^{G_K} \widetilde{\xi}$ is a crystalline lift of $\sigma$ with
the required labelled weights.

The other direction can be proved as follows using Fontaine--Laffaille theory.
We follow the notation of \cite{FL}; in particular see \S0.9 for the definition of the
category $\underline{\mathrm{MF}}^{f,p'}_{\mathrm{tor}}$,
but we take (their) $E$ to be $\mathbb{Q}_p$ and consider objects with an action of (our) $E$.
The results of \S7 and 8 of \cite{FL} imply that there is an
object $M$ of the category $\underline{\mathrm{MF}}^{f,p'}_{\mathrm{tor}}$ and
an embedding $E \to \mathrm{End}(M)$ (for large enough $E$), such that 
$$\sigma \simeq \mathrm{Hom}_E(\underline{\mathrm{U}}_\mathrm{S} (M),\overline{\mathbb{F}}_p)$$
as representations of $G_K$;
moreover decomposing $M = M_0 \oplus M_1$ (according to the idemopotents of $O_K \otimes E$
corresponding to $\tau_0,\tau_1$), each component is two-dimensional over $E$ and
$$\mathrm{Fil}^j M = \left\{\begin{array}{ll}
 M,&\mbox{if $j \le 0$}; \\ Ex_0,&\mbox{if $0 < j < k_0$;}\\ 0,&\mbox{if $j \ge k_0$,}
 \end{array}\right.$$
for some non-zero $x_0 \in M_0$.

We claim that bases $(x_i,y_i)$ for $M_i$ over $E$ can be chosen so that
$\mathrm{Fil}^1M = Ex_0$ as above and the $O_K\otimes E$-linear morphisms
$\phi^j  :  \mathrm{Fr}^*\mathrm{Fil}^j M \to M$
are defined so that $\phi^{k_0-1}x_0 = x_1$, $\phi^0y_0=  y_1$, and
one of the following holds:
\begin{itemize}
\item   $\phi^0x_1 = a x_0 +  by_0$, $\phi^0y_1 = cy_0$ for some $a,c\in E^\times$, $b\in \{0,1\}$,
\item   $\phi^0x_1 = y_0$, $\phi^0y_1 = ax_0$ for some $a \in E^\times$.
\end{itemize}
Indeed choose any basis $x_0$ for $\mathrm{Fil}^1M$ and let $x_1 = \phi^{k_0-1}x_0$.
If $\phi^0x_1 \in \mathrm{Fil}^1M$, then we are in the first case with $b=0$.
Otherwise let $z_0 = \phi^0x_1$, $y_1 = \phi^0z_0$, write $\phi^0y_1 = \alpha x_0 + \beta z_0$
for some $\alpha,\beta \in E$, and note that $\alpha \neq 0$.  If $\beta = 0$, then let
$y_0 = z_0$, giving the second case; otherwise let $y_0 = z_0+ \beta^{-1}\alpha x_0$,
giving the first case with $b = 1$.

In the first case, $M$ is reducible (as an object of $\underline{\mathrm{MF}}^{f,p'}_{\mathrm{tor}}$
with $E$-action), fitting in an exact sequence $0 \to M' \to M \to M'' \to 0$, 
where $M' = Ey_0 \oplus Ey_1$. 
It follows that $\sigma$ has the form $\left(\begin{array}{cc}  \chi_1 & * \\ 0 & \chi_2 \end{array}\right)$ where
$\chi_1$ (resp.~$\chi_2$) is obtained by applying the functor 
$\mathrm{Hom}_E(\underline{\mathrm{U}}_\mathrm{S} (\cdot),\overline{\mathbb{F}}_p)$ to $M'$
(resp.~$M''$).  Moreover $\chi_1$ (resp.~$\chi_2$) has a crystalline lift of Hodge--Tate type $(0,0)$
(resp.~$(k_0-1,0)$) and the subspace of $H^1(G_K,\overline{\mathbb{F}}_p(\chi_1\chi_2^{-1}))$
obtained from such extensions is one-dimensional.  Therefore $\chi_1$ is unramified, 
$\chi_2|_{I_K} = \epsilon_{\tau_0}^{1-k_0}$, and since the subspace 
must contain $V_{\chi_1\chi_2^{-1}}' = V_{\chi_1\chi_2^{-1}}$, these subspaces in fact
coincide, so $\sigma$ has the required form.

In the second case, consider $\sigma|_{G_{K'}}$, which (in view of the compatibility noted at the end
of \S3 of \cite{FL}) is isomorphic to $\mathrm{Hom}_E(\underline{\mathrm{U}}'_\mathrm{S} (M'),\overline{\mathbb{F}}_p)$,
where $M' = M\otimes_{O_K}O_{K'}$ and $\underline{\mathrm{U}}'_\mathrm{S}$ is the Fontaine--Laffaille
functor on the category $\underline{\mathrm{MF}}^{f,p'}_{\mathrm{tor}}$ defined with $K$ replaced by $K'$.
Assuming $E$ is chosen sufficiently large (in particular containing the residue field of $K'$), we may
decompose $M' = \oplus M_i'$ according to the embeddings $\tau_i' =  \mathrm{Fr}^i\circ\tau_0'$
for $i=0,1,2,3$, where $\tau_0'$ is a choice of extension of $\tau_0$.  Writing $x_i',x_{i+2}'$ (resp.~$y_i',y_{i+2}'$) for
the image of $x_i\otimes 1$ (resp.~$y_i\otimes 1$) in the corresponding component,
observe that $M'$ decomposes as
$$\left(Ex_0' \oplus Ey_1'\oplus Ey_2' \oplus Ex_3'\right) \quad
\bigoplus \quad \left( Ey_0' \oplus Ex_1'\oplus Ex_2' \oplus Ey_3'\right).$$
It follows that $\sigma|_{G_{K'}} \simeq \xi \oplus \xi'$ where $\xi$
has a crystalline lift of Hodge--Tate type $(k_0-1,0,0,0)$, so that
$\xi|_{I_{K'}} = \epsilon_{\tau_0'}^{1-k_0}$ (note that similarly $\xi' = \epsilon_{\tau_2'}^{1-k_0}$) 
and $\sigma$ has the required form.
\epf

\begin{remark} For completeness we note that $\sigma$ has a crystalline lift of weight $((1,1),(0,0))$ if and
only if it is unramified.\end{remark}

\begin{remark} We remark that the non-semisimple representations of $G_K$ occurring in the statement of the lemma
are precisely those which are gently (but not tamely) ramified in the terminology of \cite[\S3.3]{DDR}.
This is a special case of Conjecture~7.2 of \cite{DDR}, proved in \cite{CEGM}. \end{remark}

\begin{lemma}  \label{lem:pwt1shift} Suppose that $2 \le k_0 \le p$.  A
representation $\sigma:G_K \to \mathrm{GL}_2(\overline{\mathbb{F}}_p)$
has a crystalline lift of weight $((k_0,1),(0,0))$ if and only if all of the following hold:
\begin{enumerate}
\item $\sigma$ has a crystalline lift of weight $((k_0-1,p+1),(0,0))$ if $k_0>2$, and 
of weight $((p+1,p),(0,0))$ if $k_0=2$,
\item $\sigma$ has a crystalline lift of weight $((k_0+1,p+1),(-1,0))$,
\item and $\sigma$ is not of the form $\left(\begin{array}{cc}  \chi_1 & * \\ 0 & \chi_2 \end{array}\right)$
where $\chi_1|_{I_K} = \epsilon_{\tau_0}$.
\end{enumerate}
\end{lemma}
\begpf We firstly prove the ``only if" direction. Suppose that $\sigma$ has a crystalline lift of weight $((k_0,1),(0,0))$.

First consider the case that $\sigma$ is reducible, so by Lemma~\ref{lem:pwt1explicit}, 
it is an unramified twist of a representation 
of the form $\left(\begin{array}{cc}  1 & * \\ 0 & \chi^{-1}  \end{array}\right)$
with $\chi|_{I_K} = \epsilon_{\tau_0}^{k_0-1}$ and associated
extension class in $V_{\chi}$. 

For 1), note that $\chi|_{I_K} = \epsilon_{\tau_0}^{k_0-2}\epsilon_{\tau_1}^{p}$
(resp.~$\epsilon_{\tau_0}^{p}\epsilon_{\tau_1}^{p-1}$) if $k_0>2$  (resp.~$k_0=2$), so that
$\chi$ has a crystalline lift $\widetilde{\chi}''$ of Hodge--Tate type $(2-k_0,-p)$ (resp.~$(-p,1-p)$).
Since $H_f^1(G_K,L(\widetilde{\chi}'')) = H^1(G_K,L(\widetilde{\chi}''))$ and $H^1(G_K,\mathcal{O}(\widetilde{\chi}''))$
maps surjectively to $H^1(G_K,E(\chi))$, it follows as in the proof of Lemma~\ref{lem:pwt1explicit}
that $\sigma$ has a crystalline lift of the required weight.

For 2), we instead write $\chi|_{I_K} = \epsilon_{\tau_0}^{k_0}\epsilon_{\tau_1}^{-p}$ and use the
lift $\widetilde{\chi}$ in the definition of $V_\chi$.  Since the extension class associated to $\sigma$ lies
in $V_\chi$, it follows that $\sigma$ has a crystalline lift with $\tau_0$-labelled weights $(k_0,0)$
and $\tau_1$-labelled weights $(0,-p)$.  Twisting by a crystalline character of Hodge-Tate type
$(-1,p)$ and trivial reduction, we conclude that $\sigma$ has a crystalline lift of the required weight.

Finally 3) is clear since $\epsilon_{\tau_0}$ is not $\epsilon_{\tau_0}^{1-k_0}$ or the trivial character.

Now suppose that $\sigma$ is irreducible, so  $\sigma\simeq \mathrm{Ind}_{G_{K'}}^{G_K} \xi$
where $\xi|_{I_{K'}} = \epsilon_{\tau_0'}^{1-k_0}$ for some extension $\tau_0'$ of
$\tau_0$.  Writing $\xi|_{I_K'} = \epsilon_{\tau_0'}^{2-k_0}\epsilon_{\tau_3'}^{-p}$
(resp.~$\epsilon_{\tau_2'}^{-p}\epsilon_{\tau_3'}^{1-p}$) if $k_0>2$ (resp.~$k_0=2$),
we see that $\xi$ has a lift $\tilde{\xi}$ of Hodge--Tate type $(k_0-2,0,0,p)$ (resp.~$(0,0,p,p-1)$),
and $\mathrm{Ind}_{G_{K'}}^{G_K}\tilde{\xi}$  is a crystalline lift of $\sigma$ of the
required weight for 1).    

For 2), we proceed similarly by writing 
$\xi|_{I_{K'}} = \epsilon_{\tau_0'}^{1-k_0}\epsilon_{\tau_1'}^{-p}\epsilon_{\tau_2'}$ to
see that $\xi$ has a crystalline lift of Hodge--Tate type $(k_0-1,p,-1,0)$ whose induction
to $G_K$ has the required weight.

Finally 3) is clear since $\sigma$ is irreducible.

We now prove the ``if direction". Suppose that 1), 2) and 3) all hold.  We will use the results of \cite{GLS}
and their extension to $p=2$ in \cite{W}, which show that if $\sigma$ has a crystalline lift of weight $(k,l)$
with $2 \le k_\tau \le p+1$ for all $\tau$, then $\sigma$ is of the form prescribed in \cite{BDJ} for the
corresponding Serre weight (i.e., that $W^{\mathrm{cris}}(\sigma) \subset W^{\mathrm{explicit}}(\sigma)$)
in the terminology of  \cite{GLS2}, but note that the conventions for Hodge--Tate weights in \cite{GLS}
and \cite{GLS2} are opposite to ours).

First suppose that $\sigma$ is reducible and write
$\sigma \simeq \left(\begin{array}{cc}  \chi_1 & * \\ 0 & \chi_2 \end{array}\right)$.
From condition 1) and \cite[Thm.~9.1]{GLS}
(extended to $p=2$ in \cite{W}), it follows that if $k_0 > 2$, then $\sigma$ is isomorphic to the reduction
of a lattice in a crystalline representation of the form
$$ \left(\begin{array}{cc}  \widetilde{\chi}_1 & * \\ 0 &  \widetilde{\chi}_2 \end{array}\right),$$
where $\widetilde{\chi}_1$ and $\widetilde{\chi}_2$ are of Hodge--Tate types $(k_0-2,p)$
and $(0,0)$ (in either order) or of Hodge--Tate types $(k_0-2,0)$ and $(0,p)$ (again in
either order).  Furthermore in the first case, if $\widetilde{\chi}_1$ has Hodge--Tate type
$(k_0-2,p)$, then the representation is necessarily decomposable, so we may exchange
$\widetilde{\chi}_1$ and $\widetilde{\chi}_2$.  We therefore conclude that
$\sigma|_{I_K}$ is of the form:
$$ \left(\begin{array}{cc}  1 & * \\ 0 &  \epsilon_{\tau_0}^{1-k_0}\end{array}\right),
\quad  \left(\begin{array}{cc}  \epsilon_{\tau_0}^{-1} & * \\ 0 & \epsilon_{\tau_0}^{2-k_0}\end{array}\right)
\quad\mbox{or}\quad
  \left(\begin{array}{cc}  \epsilon_{\tau_0}^{2-k_0} & * \\ 0 & \epsilon_{\tau_0}^{-1}\end{array}\right).$$
If $k_0=2$, then the same reasoning shows that $\sigma|_{I_K}$ is of the form:
$$\left(\begin{array}{cc}  1 & * \\ 0 &  \epsilon_{\tau_0}^{-1}\end{array}\right),
\quad  \left(\begin{array}{cc}  \epsilon_{\tau_0}^{p-1} & * \\ 0 & \epsilon_{\tau_0}^{-p}\end{array}\right)
\quad\mbox{or}\quad
  \left(\begin{array}{cc}  \epsilon_{\tau_0}^{-p} & * \\ 0 & \epsilon_{\tau_0}^{p-1}\end{array}\right).$$

Similarly, from condition 2), we find that $\sigma|_{I_K}$ is of the form:
$$ \left(\begin{array}{cc}  \epsilon_{\tau_0} & * \\ 0 &  \epsilon_{\tau_0}^{-k_0}\end{array}\right),
\quad  \left(\begin{array}{cc} 1 & * \\ 0 & \epsilon_{\tau_0}^{1-k_0}\end{array}\right)
\quad\mbox{or}\quad
  \left(\begin{array}{cc}  \epsilon_{\tau_0}^{1-k_0} & * \\ 0 & 1\end{array}\right).$$
Moreover in the second case, the associated extension class lies in $V_\chi$
(where we exchange $\chi_1$ and $\chi_2$ if necessary if $\sigma$ splits, and 
use the fact that $V_\chi$ is independent of the choice of unramified twist
in its definition).   That $\sigma$ has the required form is then immediate on
comparing the possibilities resulting from 1) and 2), taking 3) into account
in the case $k_0=p=2$, and applying Lemma~\ref{lem:pwt1explicit}.

Finally suppose that $\sigma$ is irreducible.  Then condition 1) and \cite[Thm.~10.1]{GLS}
(extended to $p=2$) imply that $\sigma \simeq  \mathrm{Ind}_{G_{K'}}^{G_K} \xi$ 
for some $\xi$ with $\xi|_{I_{K'}}$ of the form 
$\epsilon_{\tau_0'}^{2-k_0}\epsilon_{\tau_3}^{-p} = \epsilon_{\tau_0'}^{1-k_0}$ 
(if the balanced subset $J$ in \cite[Thm.~10.1]{GLS} is $\{\tau'_0,\tau'_3\}$ or its complement) or
$\epsilon_{\tau_0'}^{2-k_0} \epsilon_{\tau_1'}^{-p} = \epsilon_{\tau_0'}^{2-k_0-p^2}$
(if $J = \{\tau'_0,\tau_1'\}$ or its complement), with the latter possibility replaced by $\epsilon_{\tau_0'}^{p-p^2-p^3}$
if $k_0=2$.  Similarly condition 2) implies that $\sigma \simeq  \mathrm{Ind}_{G_{K'}}^{G_K} \xi'$ 
for some $\xi'$ with $\xi'|_{I_K'}$ of the form $\epsilon_{\tau_0'}^{1-k_0}$ or
$\epsilon_{\tau_0'}^{p^2-k_0}$.  Since neither $\epsilon_{\tau_0'}^{p^2-k_0}$ nor
its conjugate $\epsilon_{\tau_0'}^{1-p^2k_0}$ agrees with any of the possibilities resulting
from condition 1), we deduce that $\xi|_{I_K'} = \epsilon_{\tau_0'}^{1-k_0}$, and the desired conclusion
again follows from Lemma~\ref{lem:pwt1explicit}.
\epf

\begin{remark} \label{rmk:123}
Note that we only needed to use condition 3) in the case $k_0=p=2$, so it is otherwise implied by 1) and 2). \end{remark}

 \subsection{Weight shifting}  
 
 We now prove an analogue of Lemma \ref{lem:pwt1shift} in the context of geometric modularity.
 (See~\cite[Chapter~4]{HWPhD} for generalisations to the setting of arbitrary totally real $F$
 in which $p$ is unramified.)

\begin{lemma}  \label{lem:pwt1shift2} Suppose that $2 \le k_0 \le p$ and that 
$\rho:G_F \to \mathrm{GL}_2(\overline{\mathbb{F}}_p)$ is irreducible.
If $\rho$ is geometrically modular of weight $((k_0,1),(0,0))$, then
\begin{enumerate}
\item $\rho$ is geometrically modular of weight $((k_0-1,p+1),(0,0))$ if $k_0>2$, and 
of weight $((p+1,p),(0,0))$ if $k_0=2$, and
\item $\rho$ is geometrically modular of weight $((k_0+1,p+1),(-1,0))$.
\end{enumerate}
Moreover the converse holds if we assume in addition that
\begin{enumerate}
\item[(3)] $\rho|_{G_K}$ is not of the form $\left(\begin{array}{cc}  \chi_1 & * \\ 0 & \chi_2 \end{array}\right)$
where $\chi_1|_{I_K} = \epsilon_{\tau_0}$.
\end{enumerate}
\end{lemma}
\begpf Suppose first that $\rho$ is geometrically modular of weight $((k_0,1),(0,0))$, i.e. that
$\rho$ is equivalent to $\rho_f$ for some eigenform $f \in M_{(k_0,1),(0,0)}(U;E)$.  Multiplying $f$ by $\mathrm{Ha}_{\tau_0}$
(resp.~$\mathrm{Ha}_{\tau_0}\mathrm{Ha}_{\tau_1}$) if $k_0>2$ (resp.~$k_0=2$) yields an eigenform giving
rise to $\rho$ of the weight required for 1).   Conclusion 2) is immediate from Theorem~\ref{thm:thetarho}.

Conversely suppose 1), 2) and 3) all hold.   First consider the case $k_0 > 2$.
By Proposition~\ref{prop:normalised}, hypotheses 1) and 2) imply that $\rho$ arises from
normalised eigenforms in $M_{(k_0-1,p+1),(0,0)}(U_1(\mathfrak{m}_1);E)$ and $M_{(k_0+1,p+1),(-1,0)}(U_1(\mathfrak{m}_2);E)$
for some ideals $\mathfrak{m}_1$, $\mathfrak{m}_2$ prime to $p$ (and sufficiently large $E$).
 We may then choose $\mathfrak{n}$ satisfying the
conditions in Lemma~\ref{lem:stable} with $\mathfrak{m} = \mathfrak{m}_i$ for $i=1,2$ (for example
take $\mathfrak{n} = \mathfrak{m}_1^2\mathfrak{m}_2^2$) to deduce that $\rho$ arises from
stabilised eigenforms in 
$f_1 \in M_{(k_0-1,p+1),(0,0)}(U_1(\mathfrak{n});E)$ and $f_2 \in M_{(k_0+1,p+1),(-1,0)}(U_1(\mathfrak{n});E)$.
By Proposition~\ref{prop:Thetaq}, $\Theta_{\tau_0}(f_1)$ is a strongly stabilised eigenform in 
$M_{(k_0,2p+1),(-1,0)}(U_1(\mathfrak{n});E)$.   By Corollary~\ref{cor:ordinary} and hypothesis 3),
so is $\mathrm{Ha}_{\tau_0}f_2$.  Lemma~\ref{lem:strong} then implies that $\Theta_{\tau_0}(f_1)
= \mathrm{Ha}_{\tau_0}f_2$, and now it follows from Theorem~\ref{thm:theta} that
$f_1 = \mathrm{Ha}_{\tau_0}f$ for some $f \in M_{(k_0,1),(0,0)}(U_1(\mathfrak{n});E)$,
so $\rho$ is geometrically modular of weight $((k_0,1),(0,0))$.

The case $k_0 =2$ is similar, but instead one has $f_1 \in M_{(p+1,p),(0,0)}(U_1(\mathfrak{n});E)$,
and obtains $f \in M_{(p+2,0),(0,0)}(U_1(\mathfrak{n});E)$.  Theorem~1.1 of \cite{DK} now implies that
$f$ is divisible by $\mathrm{Ha}_{\tau_1}$, so that $\rho$ is geometrically modular of weight $((2,1),(0,0))$.
\epf

\subsection{Geometric modularity in partial weight one}

\begin{theorem}  \label{thm:quadratic1} Suppose that $2 \le k_0 \le p$ and that 
$\rho:G_F \to \mathrm{GL}_2(\overline{\mathbb{F}}_p)$ is irreducible and modular.  
Suppose that Conjecture~3.14 of \cite{BDJ} and Conjecture~\ref{conj:algvsgeomweights}
hold for $\rho$.  Then $\rho$ is geometrically modular of weight $((k_0,1),(0,0))$ if and only if
$\rho|_{G_K}$ has a crystalline lift of weight $((k_0,1),(0,0))$.
\end{theorem}
\begpf Suppose first that $\rho|_{G_K}$ has a crystalline lift of weight $((k_0,1),(0,0))$.
Lemma~\ref{lem:pwt1shift} implies that $\rho|_{G_K}$ has crystalline lifts of weight
$((k_0-1,p+1),(0,0))$ (resp.~$((p+1,p),(0,0))$) if $k_0 > 2$ (resp.~$k_0=2$) and
$((k_0+1,p+1),(-1,0))$, and that $\rho|_{G_K}$ has no subrepresentation on which
$I_K$ acts as $\epsilon_{\tau_0}$.
Conjecture~3.14 of~\cite{BDJ} then implies that $\rho$ is algebraically modular of
weights of the two indicated weights, and then Conjecture~\ref{conj:algvsgeomweights}
implies it is geometrically modular of those weights.  It then follows from Lemma~\ref{lem:pwt1shift2}
that $\rho$ is geometrically modular of weight $((k_0,1),(0,0))$.

Now suppose that $\rho$ is geometrically modular of weight $((k_0,1),(0,0))$.
We can then reverse the argument to conclude that $\rho|_{G_K}$ has crystalline lifts of weight
$((k_0-1,p+1),(0,0))$ (resp.~$((p+1,p),(0,0))$) if $k_0 > 2$ (resp.~$k_0=2$) and
$((k_0+1,p+1),(-1,0))$.  If $p>2$, then as noted in Remark~\ref{rmk:123},
this already implies that $\rho|_{G_K}$ has a crystalline lift of weight $((k_0,1),(0,0))$.
To conclude, we can assume $k_0=p=2$, and we just need to rule out the possibility that 
$\rho|_{G_K} \simeq\left(\begin{array}{cc}  \chi_1 & * \\ 0 & \chi_2 \end{array}\right)$
where $\chi_1|_{I_K} = \epsilon_{\tau_0}$.  We do this by an ad hoc argument.

It is more convenient to work with $\rho' = \rho\otimes\rho_{\xi'}$, where $\xi$ is a character of
weight $(1,0)$.  Then $\rho'$ is geometrically modular of weight $((2,1),(1,0))$,
Conjecture~\ref{conj:algvsgeomweights} holds as well for $\rho'$, and we assume
for the sake of contradiction that $\rho'|_{G_K}$ has an unramified subrepresentation.
We let $v_0 = 2O_F$.

By Lemma~\ref{lem:levelU1}, $\rho$ arises from an eigenform $f_0 \in M_{(2,1),(1,0)}(U_1(\mathfrak{m}_1);E)$
for some $\mathfrak{m}_1$ and $E$.  Moreover by Lemma~\ref{lem:U1better}, we can assume that
 $\Theta_{\tau_1}(f_0) \neq 0$, i.e. that $r_m^J(f_0) \neq 0$ for some $m,J$ such that $m \not\in 2J^{-1}$.
The same argument as in the proof of Proposition~\ref{prop:normalised} then shows that we may assume
$f_0$ satisfies the first two conditions in the definition of a normalised eigenform.  (With regard to the third
condition, note that we have not defined $T_{v_0}$ in this context.)  Therefore $\Theta_{\tau_0}(f_0)$
is a normalised eigenform in $M_{(3,3),(0,0)}(U_1(\mathfrak{m}_1);E)$; note that it is an eigenform for
$T_{v_0}$ since $r_m^J(\Theta_{\tau_1}(f_0)) = 0$ if $m \in 2J^{-1}$.  By Theorem~\ref{thm:theta},
we have $\nu(\Theta_{\tau_0}(f_0)) \le_{\mathrm{Ha}} (4,1)$ (where the notation is as in \S\ref{sec:Hasse}),
and  Theorem~1.1 of \cite{DK} then implies that
$\nu(\Theta_{\tau_0}(f_0)) \le_{\mathrm{Ha}} (2,2)$.   We may therefore write $\Theta_{\tau_0}(f_0) =
\mathrm{Ha}_{\tau_0}\mathrm{Ha}_{\tau_1} f_1$ for a normalised eigenform $f_1 \in 
M_{(2,2),(0,0)}(U_1(\mathfrak{m}_1);E)$ with $r_m^J(f_1) = 0$ for all $m \in 2J^{-1}$.

We have shown in particular that $\rho'$ is geometrically modular of weight $((2,2),(0,0))$, hence
algebraically modular of weight $((2,2),(0,0))$ by our supposition of Conjecture~\ref{conj:algvsgeomweights}.
Therefore (for example by \cite[Prop.~2.5]{BDJ}), $\rho' \simeq \overline{\rho}_{\tilde{f}}$ for a
characteristic zero eigenform $\tilde{f}$ of weight $((2,2),(0,0))$; we may further assume that
$\tilde{f}$ is a newform in $M_{(2,2),(0,0)}(U_1(\mathfrak{m}_2);\mathcal{O})$ for some $\mathfrak{m}_2$,
enlarging $L$ if necessary, so it is a normalised eigenform for $T_v$ for all primes $v$, and 
$S_v$ for all $v \nmid \mathfrak{m}_2$.   By local-global compatibility, $\rho_{\tilde{f}}|_{G_K}$
is crystalline with $\tau_i$-labelled weights $(1,0)$ for $i=0,1$, and the characteristic polynomial
of  $\phi^2$ on $D_{\mathrm{cris}}(\rho_{\tilde{f}}|_{G_K})$ is $X^2 - \tilde{a} X + 4\tilde{d}$
where $\tilde{a}$ is the eigenvalue of $T_{v_0}$ on $\tilde{f}$ and $\tilde{d} \in \mathcal{O}^\times$
is the eigenvalue of $S_{v_0}$. Using for example that $\rho_{\tilde{f}}|_{G_K}$ is dual to a representation
arising from a $2$-divisible group over $O_K$, we see from the form of $\rho'|_{G_K}$ that
$$\rho_{\tilde{f}}|_{G_K}  \simeq \left(\begin{array}{cc}  \widetilde{\chi}_1 & * \\ 0 & \widetilde{\chi}_2 \end{array}\right)$$
with $\widetilde{\chi}_1$ unramified and $\widetilde{\chi}_1(\mathrm{Frob}_{v_0}) = \tilde{a} \in \mathcal{O}^\times$
(and $\chi_{\mathrm{cyc}}\widetilde{\chi}_2$ is unramified with $\chi_{\mathrm{cyc}}\widetilde{\chi}_2(\mathrm{Frob}_{v_0}) = \tilde{a}^{-1}\tilde{d}$).
The reduction of $\tilde{f}$ is thus a normalised eigenform $f_2 \in M_{(2,2),(0,0)}(U_1(\mathfrak{m}_2);E)$ giving
rise to $\rho'$, with the property that the eigenvalue of $T_{v_0}$ on $f_2$ is non-zero.

As in the proof of Lemma~\ref{lem:pwt1shift2}, we can choose $\mathfrak{n}$ so that the
conditions in Lemma~\ref{lem:stable} are satisfied for $\mathfrak{m}_1$ and $\mathfrak{m}_2$,
and the proof of the lemma then yields eigenforms $g_1,g_2\in M_{(2,2),(0,0)}(U_1(\mathfrak{n});E)$
such that $g_1$ is strongly stabilised, whereas $g_2$ is stabilised and satisfies $T_{v_0} g_2 = a g_2$
for some $a \in E^\times$.  Now consider the form $f_3 = a^{-1}(g_2-g_1)$; its $q$-expansion coefficients
are given by $r_m^J(f_3) = 0$ unless $m \in 2J^{-1}_+$, in which case $r_m^J(f_3) = r_{m/2}^J(g_2)$.
In particular $f_3 \in \ker(\Theta_{\tau_0})$, so $f_3 =  \Phi_{v_0}(g_3)$ for some
$g_3 \in M_{(1,1),(0,0)}(U_1(\mathfrak{n});E)$.  By Proposition~\ref{prop:Phivonqexps},
we have $r_m^J(g_3) = r_m^J(g_2)$ for all $m,J$, so
$g_2 = \mathrm{Ha}_{\tau_0}\mathrm{Ha}_{\tau_1}g_3$.

Furthermore, note that $\nu(g_3) = (1,1)$; otherwise Corollary~1.2 of \cite{DK} would
force $\nu(g_3) = (0,0)$, making $g_3$ locally constant and contradicting the irreducibility of $\rho$.
Now consider $\Theta_{\tau_1}(g_3) \in M_{(3,2),(0,-1)}(U_1(\mathfrak{n});E)$.
By Theorem~\ref{thm:theta}, $\Theta_{\tau_1}(g_3)$ is not divisible by $\mathrm{Ha}_{\tau_1}$.
We claim that $\Theta_{\tau_1}(g_3)$ is not divisible by $\mathrm{Ha}_{\tau_0}$ either.  Indeed
if it were, then we would have $\nu(\Theta_{\tau_1}(g_3)) \le_{\mathrm{Ha}} (4,0)$, and
Theorem~1.1 of \cite{DK} would imply divisibility by $\mathrm{Ha}_{\tau_1}$.
Therefore Theorem~\ref{thm:theta} implies that $\Theta_{\tau_0}\Theta_{\tau_1}(g_3)$
is not divisible by $\mathrm{Ha}_{\tau_0}$ (and in fact a similar argument gives
$\nu(\Theta_{\tau_0}\Theta_{\tau_1}(g_3)) = (4,4)$).   Note that 
$\Theta_{\tau_0}\Theta_{\tau_1}(g_3) \in M_{(4,4),(-1,-1)}(U_1(\mathfrak{n});E)$
is a strongly stabilised eigenform giving rise to $\rho$.  However so is
$e_1\mathrm{Ha}^2_{\tau_0}\mathrm{Ha}^2_{\tau_1}g_1$, where $e_1$
is the constant section in $M_{(0,0),(-1,-1)}(U_1(\mathfrak{n});E)$ with value $1$.
We therefore conclude that
$\Theta_{\tau_0}\Theta_{\tau_1}(g_3) = e_1\mathrm{Ha}^2_{\tau_0}\mathrm{Ha}^2_{\tau_1}g_1$
is divisible by $\mathrm{Ha}_{\tau_0}$, yielding the desired contradiction.
\epf

\begin{remark} Note that the theorem holds just as well for weights of the form $((k_0,1),l)$
and $((1,k_0),l)$ for any $l \in \mathbb{Z}^\Sigma$.\end{remark}

Recall from Proposition~\ref{prop:algvsgeomweights}
that one direction of Conjecture~\ref{conj:algvsgeomweights} holds if $k$ is paritious in the sense of Definition \ref{def:paritious}.
Recall also that Conjecture~3.14 of \cite{BDJ} has been proved under mild technical
hypotheses by Gee and collaborators (see especially~\cite{GLS,GKi}), with an
alternative to part due to Newton~\cite{N}.  
In particular it holds under the assumptions that $p>2$, $\rho|_{G_{F(\zeta_p)}}$
is irreducible, and if $p=5$, then $\rho|_{G_{F(\zeta_5)}}$ does not have projective
image isomorphic to $A_5 \cong \mathrm{PSL}_2(\mathbb{F}_5)$.  It might be possible
to treat this exceptional case with $p=5$ using the same methods along with~\cite[Thm.~3.2.1]{BLGG}, but we only
need to do this in a particular instance in order to obtain one direction of 
Theorem~\ref{thm:quadratic1} for odd $k_0$ unconditionally.

\begin{theorem}  \label{thm:quadratic2} Suppose that $3 \le k_0 \le p$, $k_0$ is odd and that 
$\rho:G_F \to \mathrm{GL}_2(\overline{\mathbb{F}}_p)$ is irreducible and modular.  
If $\rho|_{G_K}$ has a crystalline lift of weight $((k_0,1),(0,0))$, then
$\rho$ is geometrically modular of weight $((k_0,1),(0,0))$.
\end{theorem}
\begin{proof}  
We first show that the local condition at $p$ implies that
$\rho|_{G_{F(\zeta_p)}}$ is irreducible.  Indeed if it is not, then
$\rho$ is induced from $G_{F'}$ for a
quadratic extension $F'/F$ which is
ramified at $p$, and hence $\rho|_{G_K}$ is induced
from $G_{K'}$ for a ramified quadratic extension $K'/K$.
This in turns implies that $\rho|_{I_K} \simeq \chi_1\oplus \chi_2$
for some characters $\chi_1,\chi_2$ such that $\chi_1\chi_2^{-1}$
is quadratic.   However the explicit description of the possibilities
for $\rho|_{I_K}$ from Lemma~\ref{lem:pwt1explicit} shows that
$\chi_1\chi_2^{-1}$ would have the form $\epsilon_{\tau_0}^{\pm(k_0-1)}$
or $\epsilon_{\tau_0'}^{\pm(k_0-1)(p^2-1)}$, which gives a contradiction
since such a character has order $(p^2-1)/i$ or $(p^2+1)/i$ for some
$i \le p-1$.

We may therefore apply Theorem~A of \cite{GLS2} to conclude that
$\rho$ is algebraically modular of weights $((k_0-1,p+1),(0,0))$ 
and $((k_0+1,p+1),(-1,0))$, unless $p=5$ and $\rho|_{G_{F(\zeta_5)}}$
has projective image isomorphic to $\mathrm{PSL}_2(\mathbb{F}_5)$.
Aside from this exceptional case, it follows from Proposition~\ref{prop:algvsgeomweights}
that  $\rho$ is geometrically modular of weights $((k_0-1,p+1),(0,0))$ 
and $((k_0+1,p+1),(-1,0))$, and then from Lemma~\ref{lem:pwt1shift2} that
$\rho$ is geometrically modular of weight $((k_0,1),(0,0))$.

Suppose then that $p=5$ and $\rho|_{G_{F(\zeta_5)}}$
has projective image isomorphic to $\mathrm{PSL}_2(\mathbb{F}_5)$, so
that of $\rho$ is isomorphic to $\mathrm{PSL}_2(\mathbb{F}_5)$ or
$\mathrm{PGL}_2(\mathbb{F}_5)$.  Again using the explicit
descriptions in Lemma~\ref{lem:pwt1explicit}, we see this is only
possible if $k_0 = 5$ and $\rho|_{G_K} \simeq  \chi_1 \oplus \chi_2$
where $\chi_1$ is unramified and $\chi_2|_{I_K} = \epsilon_{\tau_0}^{-4}$
has order $6$.   In this case the conjectural set of Serre weights
for $\rho^\vee = \mathrm{Hom}_{\overline{\mathbf{F}}_5}(\rho,\overline{\mathbb{F}}_5)$
(with the notation of \S\ref{sec:conj}) is:
\[ \{\, V_{(4,6),(0,0)}, \, V_{(2,2),(-1,0)},\, V_{(6,6),(-1,0)},\,V_{(6,4),(4,0)}\,\}.\]
In particular if $\xi$ is a character of weight $(1,0)$, then
$(\rho\otimes\rho_{\xi'})^\vee|_{G_K}$ has a Barsotti--Tate lift (necessarily non-ordinary),
and the argument of \cite[\S3.1]{G2} (using the method of Khare--Wintenberger~\cite{KW})
then shows that $(\rho\otimes\rho_{\xi'})^\vee$ is modular of weight $V_{(2,2),(0,0)}$,
from which it follows that $\rho$ is algebraically modular of weight $((2,2),(-1,0))$.

Similarly $\rho^\vee|_{G_K}$ has a potentially Barsotti--Tate lift of type
$[\epsilon_{\tau_0}^2\epsilon_{\tau_1}^4] \oplus 1$, 
so the same argument (but now using  \cite[Thm.~3.2.1]{BLGG} for the existence of ordinary lifts)
shows that $\rho^\vee$ is modular of some weight in the set of Jordan--H\"older constituents 
$\mathrm{Ind}_{B}^{\mathrm{GL}_2(O_F/p)} \psi$ where
\[. \psi\left(\left(\begin{array}{cc}a&b\\0&d\end{array}\right)\right) = \tau_0(a)^2\tau_1(a)^4,\]
namely
$\{\, V_{(4,6),(0,0)}, \, V_{(3,5),(3,0)},\, V_{(4,2),(2,4)}\,\}$.
Therefore $\rho$ is algebraically modular of weight $((4,6),(0,0))$,
$((3,5),(3,0))$ or $((4,2),(2,4))$.  Since $\rho|_{G_K}$ has no
crystalline lifts of weight $((3,5),(3,0))$ or $((4,2),(2,4))$
(by \cite[Thm.~2.12]{GLS}, but in fact already by \cite{FL}),
these possibilities are ruled out by local-global compatibility
and the discussion before Proposition~\ref{prop:algvsgeomweights}.
Therefore $\rho$ is algebraically modular of weight $((4,6),(0,0))$.

We have now shown that $\rho$ is algebraically modular of weights
$((4,6),(0,0))$ and $((2,2),(-1,0))$, so also geometrically modular
of these weights by  Proposition~\ref{prop:algvsgeomweights}.
Therefore $\rho$ is also geometrically modular of weight $((6,6),(-1,0))$,
and it follows from  Lemma~\ref{lem:pwt1shift2} that
$\rho$ is geometrically modular of weight $((5,1),(0,0))$,
as required.
\end{proof}

\begin{remark} Again the theorem holds also for weights of the form $((k_0,1),l)$
and $((1,k_0),l)$ for any $l \in \mathbb{Z}^\Sigma$.\end{remark}

\begin{remark} We remark that the assumption that $F$ is unramified at $p$ ensures
that the conditions at $p=5$ in the modularity lifting theorems
of \cite{Ki}, \cite{G} and \cite{BLGG} are satisfied.
If $p=5$ is ramified in $F$, then a more exceptional case can arise, and is treated in work of Khare and Thorne~\cite{KT}. \end{remark}

\subsection{An example}
Consider the Galois representation defined in Example IIIb${}_1$
of \cite[\S9]{DDR}, so $F = \mathbb{Q}(\sqrt{5})$, $p=3$ and
$\rho:G_F \to \mathrm{GL}_2(\mathbb{F}_9)$ is absolutely irreducible and has the property that
$\rho|_{G_K} \simeq \left(\begin{array}{cc}  \chi_1 & * \\ 0 & \chi_0 \end{array}\right)$
where $\chi_i|_{I_K} = \epsilon_{\tau_i}$ for appropriately chosen $\tau_i:O_F/p \simeq 
\mathbb{F}_9$.  Setting $\chi = \chi_1\chi_0^{-1}$, we have 
$\chi|_{I_K} = \epsilon_{\tau_0}^{-2}$, and the discussion in \cite{DDR}
shows that the associated extension class lies in the line $V_\chi$
of Lemma~\ref{lem:lines}.  It follows that $\rho$ has a crystalline lift of weight
$((3,1),(0,-1))$.   

The modularity of $\rho$ is strongly indicated by
the data exhibited in \cite[\S10.4]{DDR}.  In particular it follows from the explicit
computations described there  that there is an eigenform
$f \in M_{(2,4),(0,-1)}(U_1(\mathfrak{n});\mathbb{F}_9)$ with $\mathfrak{n} = (10\sqrt{5})$
whose eigenvalue for $T_v$ coincides with $\mathrm{tr}\rho(\mathrm{Frob}_v)$ for all $v\nmid 30$
such that $\mathrm{Nm}_{F/\mathbb{Q}}(v) < 100$, and whose eigenvalue for 
$S_v$ is $1= \mathrm{Nm}_{F/\mathbb{Q}}(v)^{-1}\det\rho(\mathrm{Frob}_v)$
for all $v\nmid 30$.
We assume for the rest of the discussion that it is indeed the case that $\rho_f \simeq \rho$.
It then follows from Theorem~\ref{thm:quadratic1}
that $\rho$ is geometrically modular of weight $((3,1),(0,-1))$, i.e., $\rho\simeq \rho_g$
for some eigenform $g \in M_{(3,1),(0,-1)}(U_1(\mathfrak{n});\mathbb{F}_9)$

Consider also the form $g_\xi$ for a character $\xi$ of conductor $(\sqrt{5})$ and weight $(0,2)$,
in the sense of Definition~\ref{def:charwt}.
(There are two such characters, both of order $4$, differing by the quadratic character corresponding
to the extension $F(\mu_5)$.)  Then we have $g_\xi \in M_{(3,1),(0,1)}(U_1(\mathfrak{n});\mathbb{F}_9)$,
and as the weight $((3,1),(0,1))$ is paritious (in the sense of Definition~\ref{def:paritious}),
it is natural to ask whether $g_\xi$ lifts to a characteristic zero eigenform of partial weight one.

\bibliographystyle{plain} 
\bibliography{refs} 

\bigskip

\end{document}